\documentclass[11pt,reqno]{amsart}
\usepackage{amsmath, amssymb, graphicx, mathrsfs, hyperref}
\usepackage{amsthm}
\usepackage{palatino}
\usepackage{autonum}
\usepackage[makeroom]{cancel}
\usepackage{enumitem}  
\usepackage[usenames,dvipsnames,svgnames]{xcolor}
\usepackage[top=1in, bottom=1in, left=1in, right=1in]{geometry} 

\newcommand{\Gal}{\operatorname{Gal}}
\newcommand{\Res}{\operatorname{Res}}
\renewcommand{\Re}{\operatorname{Re}}
\renewcommand{\Im}{\operatorname{Im}}
\renewcommand{\min}{\operatorname{min}}
\newcommand{\s}{{\sigma}}

\renewcommand{\arctan}{\tan^{-1}}
\renewcommand{\a}{\alpha}
\renewcommand{\b}{\beta}
\newcommand{\e}{\epsilon}
\renewcommand{\d}{{\delta}}

\newcommand{\G}{\Gamma}

\renewcommand{\r}{{\rho}}

\renewcommand{\pmod}[1]{\,(\textup{mod}\,#1)}
\theoremstyle{plain}
\newtheorem{theorem}{Theorem}[section]
\newtheorem{lemma}[theorem]{Lemma}
\newtheorem{corollary}[theorem]{Corollary}

\newtheorem{proposition}[theorem]{Proposition}
\newtheorem{definition}[theorem]{Definition}

\newtheorem*{remark*}{Remark}
\newtheorem{example}[]{Example}

\usepackage{import}
\usepackage{xifthen}
\usepackage{pdfpages}
\usepackage{mathtools}

\numberwithin{equation}{section}
\usepackage{enumitem, hyperref}
\makeatletter
\def\namedlabel#1#2{\begingroup
	#2%
	\def\@currentlabel{#2}%
	\phantomsection\label{#1}\endgroup
}
\makeatother

\usepackage{hyperref} 
\definecolor{labelkey}{rgb}{1,0,0}

\begin{document}
	\title[Improved Omega results for Riesz sums of a class of $L$-functions]{Improved Omega results for Riesz sums of a class of $L$-functions}

\author[Arindam Roy]{Arindam Roy}

\address{Department of Mathematics, University of North Carolina at Charlotte, NC 28223, USA}
	\email{arindam.roy@uncc.edu}

\author[Jagannath Sahoo]{Jagannath Sahoo} 

\address{
	Department of Mathematics, Indraprastha Institute of Information Technology Delhi, New Delhi 110020, India} 

\email{jagannath.sahoo@alumni.iitgn.ac.in}

\author[Akshaa Vatwani]{Akshaa Vatwani}

\address{Department of Mathematics, 
	Indian Institute of Technology Gandhinagar, Gandhinagar, 
	Gujarat 382355, India}

\email{akshaa.vatwani@iitgn.ac.in }
	

\begin{abstract}
We investigate a general class of $L$-functions satisfying functional equations with multiple gamma factors, and establish sharper lower bounds for the error terms associated with the Riesz sums of their coefficients. By integrating the resonance method with Riesz smoothing, we improve upon previously known omega results for the corresponding Riesz sums.  
 We prove these estimates within a unified axiomatic framework under a precise set of hypotheses. Our  general theorem can be applied to several classical arithmetic problems, including the $k$-fold Piltz divisor function $d_{\mathbb{K},k}(n)$ over a number field $\mathbb{K}$, as well as the lattice point counting functions for the circle, ellipse, and sphere problems.  
 In particular, we sharpen the classical $\Omega$-bound for the error
 in counting lattice points in the ellipse $u^2+5v^2\le x$,
 strengthening a result of Nowak. Moreover, for the smoothed error terms corresponding to the divisor and circle problems, we improve upon previously known $\Omega$-bounds by increasing the   exponent of the $\log\log x$ factor and sharpening the $\log\log\log x$ dependence by replacing  exponential decay factors of the form $\exp(-c\sqrt{\log\log\log x})$ in previous results with an explicit polynomial power of $\log\log\log x$. 
\end{abstract}

	\maketitle
	\section{Introduction}
	Two of the most enduring problems in analytic number theory -- the Dirichlet divisor problem and the Gauss circle problem -- concern the oscillatory behavior of lattice point counting functions. Specifically, they address the distribution of the divisor function $d(n)$ and the number of representations of an integer as a sum of two squares, $r(n)$. Since the seminal works of  Gauss (1834) and  Dirichlet (1849), it has been understood that the summatory functions
\begin{align}
	\sum_{n\le x} r(n) \quad \text{and} \quad \sum_{n\le x} d(n),
\end{align}
admit asymptotic expansions with dominant main terms $\pi x$ and $x\log x + (2\gamma-1)x$, respectively, accompanied by highly oscillatory error terms, denoted $P(x)$ and $\Delta(x)$.
Despite nearly two centuries of progress, the optimal upper bounds for these error terms remain unknown. The conjectured exponent $\theta = 1/4$ in bounds of the form
\begin{align}
	P(x),\Delta(x) \ll x^{\theta+\varepsilon},
\end{align}
continues to elude proof. While substantial effort has been devoted to upper bounds, the corresponding lower bounds -- expressed via $\Omega$-estimates -- play an equally fundamental role: they reveal the intrinsic magnitude and irregularity of the arithmetic fluctuations underlying these lattice point problems. 

\normalcolor

A powerful method for probing these oscillations is the introduction of Riesz means, which smooths the sharp cutoff in the summatory functions. For $\rho \ge 0$, define the normalized Riesz means
\begin{align}
	D_\rho(x)=\frac{1}{\Gamma(\rho+1)}\sum_{n\le x} d(n)(x-n)^\rho,
	\qquad
	R_\rho(x)=\frac{1}{\Gamma(\rho+1)}\sum_{n\le x} r(n)(x-n)^\rho.
\end{align}
These weighted sums admit asymptotic expansions whose error terms, $\Delta_\rho(x)$ and $P_\rho(x)$, exhibit behavior that depends delicately on the smoothing parameter $\rho$. The study of their lower bounds traces a rich historical development, reflecting the interplay between analytic smoothing and arithmetic oscillation.

The principal goal of this paper is to establish sharper $\Omega$-estimates  for the error terms of Riesz means associated with a  class $\mathcal{C}$ of $L$-functions satisfying functional equations with multiple Gamma factors. Specifically, we investigate the structural interaction between the resonance method and Riesz smoothing. The weight $(x-n)^\rho$ induces rapid decay in the high-frequency components of the associated Bessel or trigonometric series, effectively acting as a low-pass filter, whereas resonance arguments depend on amplifying  oscillations to produce extreme values.  By  systematically accounting for these two effects,   we improve upon known results of 
Hafner \cite{Hafner}  for Riesz sums associated to this class. Our results demonstrate that despite the damping constraints imposed by analytic smoothing, the resonance method can be successfully adapted to extract  larger arithmetic fluctuations than previously known for this class of functions.

Moreover, we develop a unified axiomatic framework,  relying on a precise set of hypotheses, that yields improved $\Omega$-bounds for the residual error. Mirroring the strategy used in the resonance method for individual classical problems, this framework formalizes  the crucial ingredient, which  is the identification of a ``distinguished set" $\mathcal{M}$ on which the average value of the coefficients is larger than usual  (see \eqref{condition of H4}).  
It is the product of this average  with the size of the set $ \mathcal{M}$ -- in other words,  the summatory function  $\sum_{m \in \mathcal M} b_m$
of the coefficients over $\mathcal{M}$ -- that shows up in
the shape of our  general $\Omega$-bound (see \eqref{def:g(x)} and \eqref{omega bd}). 


As concrete applications of our main general theorem, we improve upon known foundational bounds for several classical arithmetic functions in the critical range of smoothing $0 \le \rho < 1/2$.  Our main results include the following. Here and henceforth, 
$\log_m x$ denotes  the logarithm iterated $m$ times.
\begin{itemize}[leftmargin=*]
	\item \textbf{The Piltz Divisor Problem:} For the Riesz means  of the $k$-fold divisor function $d_{\mathbb{K}, k}(n)$ associated to a number field $\mathbb{K}$, we sharpen known omega results by   improving the exponent of the $\log_2 x$ term as well as refining the $\log_3 x$ component in the bound. 
	More precisely,   previous exponential decay factors of the form $\exp(-c\sqrt{\log_3 x})$ are replaced with an  explicit power of $\log_3 x$.  
	
In particular, this improves the $\Omega$-bound for the smoothed error term $\Delta_\rho(x)$ in the classical divisor problem, in the  range $0<\rho<1/2$, improving upon the work of Hafner \cite{Hafner}. 
When $\rho= 0$, our  exponent for the $\log_2 x$ factor (see \eqref{omega bdd for div fn}) reduces to  the value obtained by Soundararajan, which is thought to be optimal (see \cite[p.~1990]{Sound}).  
Analogously, it is likely that our exponent for the $\log_2 x$ term  captures the true growth of $\Delta_\rho(x)$  across the  range  $0\le \rho \le 1/2$. 

	\item \textbf{The Circle Problem:} We achieve structurally identical improvements for the smoothed error term $P_\rho(x)$ concerning the number of representations of an integer as a sum of two squares, yielding a larger $\log_2 x$ exponent and  a precise power of  $\log_3 x$.
	
	\item \textbf{The Sphere Problem:} We establish new omega results for the Riesz means of $r_3(n)$:  the number of ways of expressing $n$ as a sum of three squares, in the extended range $0 \le \rho < 1$. This extends prior results of Tsang \cite{TsangSphere} for $\rho=0$ and  provides sharper quantitative information compared to the classical limit theorems of Chandrasekharan and Narasimhan \cite{ChandrasekharanV} (see \eqref{CN lim result} below).
	\item \textbf{An Ellipse Problem:}
	For the positive-definite quadratic form
	$\mathcal Q(u,v)=u^2+5v^2,$
	we study the Riesz means associated with 
	$r_{\mathcal Q}(n)$, which counts the number of representations of $n$
	by $\mathcal Q$. We obtain improved $\Omega$-bounds for the corresponding
	error term $P_\rho(x,\mathcal Q)$ in the range $0\leq \rho<1/2$. In particular, when $\rho=0$, our result improves the 
	classical
	$
	\Omega_-\bigl(x^{1/4}(\log x)^{1/4}\bigr)
	$
	bound of Nowak \cite{Nowak1985} for lattice-point discrepancies of smooth convex planar
	domains. 
	\normalcolor
	
\end{itemize}

\section{Main results}
We will  not only sharpen  known omega results for \(\Delta_{\rho}(x)\) and \(P_{\rho}(x)\), but also treat a substantially more general class of functions (see Definition~\ref{def:fn eqn of F}), and explore  applications to the Piltz divisor problem and the sphere problem. We formulate explicit hypotheses (\ref{H1}-\ref{H3}) under which omega results for the error terms of the associated Riesz sums follow uniformly for every function in this class.
More precisely, for \(x>0\) and \(\rho \ge 0\), we consider the Riesz sum
\begin{align} \label{def:A_rho}
	A^{\rho}(x)
	:=
	\frac{1}{\Gamma(\rho+1)}\sideset{}{'}\sum_{\lambda_n \le x} a_n (x-\lambda_n)^{\rho},
\end{align}
associated with the Dirichlet series \(\sum_{n=1}^{\infty} a_n \lambda_n^{-s}\), where \((\lambda_n)\) is a strictly increasing sequence of real numbers tending to infinity. Here, and throughout, the prime on the summation indicates that the last term is weighted by \(1/2\) when \(\rho=0\) and \(x=\lambda_n\) for some \(n\). Recall that for a real valued function $f$ and a positive function $g$, we write $f=\Omega (g)$  if  $\limsup_{x \rightarrow \infty} \frac{|f(x)|}{g(x)} > 0$.
Similarly, we write  $f=\Omega_+ (g)$
(respectively $f=\Omega_- (g)$) if 
$\limsup_{x \rightarrow \infty} \frac{f(x)}{g(x)} > 0$ 
(respectively $\liminf_{x \rightarrow \infty} \frac{f(x)}{g(x)} < 0$).
If one has both  $f=\Omega_+ (g)$ and  $f=\Omega_- (g)$,
then we write 
$f=\Omega_\pm (g)$.

Our main result, Theorem \ref{thm:main theorem} stated later, holds for a general class of functions. We first illustrate this result by giving $\Omega$-bounds for several classical examples in Section \ref{sec: appl}.   Comparison of these bounds with existing results is given in Section \ref{sec: comparison}. 

\subsection{Application to Classical Problems}
\label{sec: appl}
Let $a_n$ be the Piltz divisor function 
$d_{\mathbb{K}, k}(n)$
for a  number field $\mathbb{K}$, defined in terms of the Dedekind zeta function 
$\zeta_{\mathbb{K}}(s) $
by the relation,
\begin{equation}
	\zeta^k_{\mathbb{K}}(s) 
	=
	\sum_{n=1}^{\infty}
	\frac{d_{\mathbb{K},k}(n)}{n^s} \quad (\Re(s) > 1).
\end{equation}
The error term corresponding to the Riesz sum for this function is given by
\begin{align}
	\Delta_{\mathbb{K}, k,\rho}(x)
	= 
	\frac{1}{\Gamma(\rho+1)}\sideset{}{'}\sum_{n \leq x} d_{\mathbb{K},k}(n)(x-n)^\rho 
	- 
	\frac{1}{2 \pi i} \int_{C}\frac{ \zeta^k_{\mathbb{K}}(s)  \Gamma(s)x^{s+\rho}}{\Gamma(s+\rho+1)} ds,   	
\end{align}
where $C$ is a circle of finite radius containing points
$0,1$
and lying inside the strip
$-\epsilon < \Re(s)< 1+ \epsilon$ 
for some 
$\epsilon >0$. 
In order to state the omega result for 
$\Delta_{\mathbb{K}, k,\rho}(x)$, we need the following notation. 

We will denote the degree $[\mathbb K : \mathbb Q]$ by   $m$. 
For  $0 \leq v \leq m$, 
let 
$P_v$ 
denote the set of all rational primes which are unramified in 
$\mathbb{K}$
and which are divisible by exactly 
$v$-many $\mathcal{O}_{\mathbb{K}}$-prime ideals of inertia degree $1$. 
Let 
$\d_v$ 
be the Dirichlet densities of the sets 
$P_v$, which can be calculated as follows. 
Let 
$\mathbb{L}$
be the minimal normal extension of 
$\mathbb{K}$, 
$G= \Gal(\mathbb{L}/\mathbb{Q})$
the Galois group of 
$\mathbb{L}$
over 
$\mathbb{Q}$
and 
$H= \Gal(\mathbb{L}/\mathbb{K})$ 
the corresponding Galois subgroup of $G$. Then, we have
\begin{align} \label{def delta v}
	\delta_v
	= 
	|G|^{-1} 
	\big|\{g \in G: |\{\sigma \in G: g \in \sigma H \sigma^{-1}\}| = v |H|\}\big|.
\end{align} 
Define  the set 
\begin{align}\label{set}
I \coloneqq	\{ 1 \leq v \leq m : \d_v > 0\}. 
\end{align}
We note that the constants $\delta_v$ 
satisfy 
\begin{align}
	\sum_{v \in I} \delta_v v = 1. 
	\label{delta v}
\end{align}
\begin{theorem}  	\label{thm: piltz}
Let $\mathbb K$ be a number field of degree $m$. Let  $r_1$ and $2r_2$ be the number of real and complex embeddings of $\mathbb{K}$ respectively.
Let $\delta_v$ and $I$ be as in \eqref{def delta v} and \eqref{set} respectively and $R= |I|$. Then for any  $ k \in \mathbb N$ satisfying $km>1$ and   
	$0 \leq \rho <  \frac{km-1}{2}$, we have 
	\begin{align}\label{final piltz div bdd}
		\Delta_{\mathbb{K}, k,\rho}(x)
		= 
		\Omega \Bigg(
		x^{\theta}
		(\log x)^{\kappa}
		(\log_2 x)^{\tau'}
		(\log_3 x)^{\frac{R}{2}
			(\kappa-1)}\Bigg), 
	\end{align}
	where 
	\begin{align}\label{tau' in piltz}
		\theta = \frac{(km-1)(2 \rho +1)}{2km}, \quad
		\kappa = \frac{km-1- 2 \rho}{2km}, \quad 
		\tau' = \frac{km+1+2 \rho}{2km} \Big(\sum_{v \in I} \d_v (vk)^{\frac{2km}{km+1+2 \rho}}-1\Big). 
	\end{align}
	Moreover, the omega bound \eqref{final piltz div bdd} is true for   $\Re(\Delta_{\mathbb{K}, k,\rho}(x))$    
	with $\Omega$ replaced by 
	$\Omega_+$  (resp. $\Omega_-$) 
	whenever 
	$kr_1 - 2\rho \equiv 3 \pmod 8$ (resp. 
	$kr_1 - 2\rho  \equiv 7 \pmod 8$).
\end{theorem}
In the case 
$\mathbb{K}= \mathbb{Q}$, $d_{\mathbb{K},k}(n)$
is simply the usual Piltz divisor function 
$d_{k}(n)$, and we have 
\begin{corollary} \label{cor: d_k}
	For $k\ge 2$, $k \in \mathbb N$ and
	$0 \le \rho < \frac{k-1}{2}$, 
	we have
	\begin{align}\label{final omega bdd d_k}
		\Delta_{k,\rho}(x)
		=
		\Omega
		\left(
		x^{\theta}
		(\log x)^{\kappa} 
		(\log_2 x)^{\tau'}
		(\log_3 x)^{-\frac{\r}{2k}-\frac{1}{4k} - \frac{1}{4}}\right),
	\end{align}
	where 
	\begin{align}
		\theta = \frac{(k-1)(2 \rho +1)}{2k},  \:\:
		\kappa = \frac{k-1- 2 \rho}{2k} \: \text{ and } \:  \tau' = \frac{k+1+2 \rho}{2k} \left(k^{\frac{2k}{k+1+ 2 \rho}}-1\right). 
	\end{align}
	Moreover,   the bound  \eqref{final omega bdd d_k} is true for $\Re(\Delta_{k,\rho}(x))$ with $\Omega$ replaced by 
	$\Omega_+$  (resp. $\Omega_-$) 
	whenever 
	$k - 2\rho \equiv 3 \pmod 8$ (resp. 
	$k - 2\rho  \equiv 7 \pmod 8$).  \normalcolor
\end{corollary} 
In particular, for the classical divisor function $d(n)$, we obtain for $0\le \rho < 1/2$,
\begin{align} \label{omega bdd for div fn}
	\Delta_\rho(x)
	=
	\Omega\left(
	x^{\theta}
	(\log x)^{\kappa}
	(\log_2 x)^{\tau'}
	(\log_3 x)^{-\frac{\rho}{4}-\frac{3}{8}}
	\right),
\end{align}
where \(\theta = \tfrac{2\rho+1}{4}\), \(\kappa = \tfrac{1-2\rho}{4}\), and
$
\tau' = \left(\frac{3+2\rho}{4}\right)\left(2^{\frac{4}{3+2\rho}} - 1\right).
$

In the case of the circle problem, 
for $\rho \geq 0$ and $x>0$, we define the error term $P_{\rho}(x)$ by
\begin{align}
	P_{\rho}(x)
	&:=
	A^{\rho}(x)-Q_{\r}(x)
	\\
	&= 
	\frac{1}{\Gamma(\rho+1)}\sideset{}{'}\sum_{n \leq x} r(n)(x-n)^\rho 
	- \frac{1}{2 \pi i} \int_{C}\frac{\phi(s)\Gamma(s)x^{s+\rho}}{\Gamma(s+\rho+1)} ds,	
\end{align}
using \eqref{def:A_rho} and \eqref{def:Q_rho}, 
where $C$ is a circle of finite radius containing points
$0,1$
and lying inside the strip
$-\epsilon < \Re(s)< 1+ \epsilon$ 
for some 
$\epsilon >0$. 
We obtain the following result.
\begin{theorem} \label{thm: r(n)}
	For  $0 \leq \rho < \frac{1}{2}$, we have 
	\begin{align}\label{final omega bdd r(n)}
		P_{\rho}(x)
		=
		\Omega
		\left(
		x^{\theta}
		(\log x)^{\kappa} 
		(\log_2 x)^{\tau'}
		(\log_3 x)^{-\frac{\r}{4}- \frac{3}{8}}\right),
	\end{align}
	where 
	\begin{align}
		\theta = \frac{2 \rho +1}{4},  \:\:
		\kappa = \frac{1- 2 \rho}{4} \: \text{ and } \:  \tau' = \frac{3+2 \rho}{4} \left(2^{\frac{1-2 \rho}{3+ 2 \rho}}-1\right). 
	\end{align} 
\end{theorem}

Let 
$a_n=r_3(n)$,
the number of ways to write $n$ as a sum of three squares and 
$ \phi(s):=\sum_{n=1}^{\infty}r_3(n) n^{-s}$.
For
$\rho \geq 0$ and $x>0$,
we define the error term 
$P_{3, \rho}(x)$ by
\begin{align}
	P_{3,\rho}(x)
	:= 
	\frac{1}{\Gamma(\rho+1)}\sideset{}{'}\sum_{n \leq x} r_3(n)(x-n)^\rho 
	- \frac{1}{2 \pi i} \int_{C}\frac{\phi(s)\Gamma(s)x^{s+\rho}}{\Gamma(s+\rho+1)} ds,	
\end{align}
using \eqref{def:A_rho} and \eqref{def:Q_rho},
where $C$ is a circle of finite radius containing points
$0,\frac{3}{2}$
and lying inside the strip
$-\epsilon < \Re(s)< \frac{3}{2}+ \epsilon$ 
for some 
$\epsilon >0$.
We obtain the following omega result as a consequence of our general result.
\begin{theorem} \label{thm: r_3(n)}
	For  $0 \leq \rho < 1$, we have 
	\begin{align}\label{final omega bdd r_3(n)}
		P_{3,\rho}(x)
		=
		\Omega
		\left(
		x^{\theta}
		(\log x)^{\kappa} 
		\right),
	\end{align}
	where 
	\begin{align}
		\theta = \frac{1+ \rho}{2} \: \text{ and } \:
		\kappa = \frac{1-  \rho}{2}.
	\end{align} 
	Moreover, in the special case $\rho =0$, the bound \eqref{final omega bdd r_3(n)} is true for 
	$\Re(P_{3,\rho}(x))$
	with $\Omega$ replaced by 
	$\Omega_-$.
\end{theorem} 

Let
$\mathcal{Q}(u,v):=u^2+5v^2$
be a positive definite binary quadratic form, and define
\begin{align}
	r_\mathcal{Q}(n):=\#\{(u,v)\in\mathbb Z^2:\mathcal{Q}(u,v)=n\}.
\end{align}
Let $\phi(s)$ be the associated Dirichlet series.
For $\rho\geq 0$, define
\begin{align}
	P_\rho(x,\mathcal {Q})
	:=
	\frac{1}{\Gamma(\rho+1)}
	\sideset{}{'}\sum_{n\leq x}r_\mathcal {Q}(n)(x-n)^\rho-\frac{1}{2\pi i}\int_C
	\frac{\phi(s)\Gamma(s)x^{s+\rho}}{\Gamma(s+\rho+1)}\,ds,
\end{align}
using \eqref{def:A_rho} and \eqref{def:Q_rho}, 
where $C$ is a circle of finite radius containing points
$0,1$
and lying inside the strip
$-\epsilon < \Re(s)< 1+ \epsilon$ 
for some 
$\epsilon >0$.
We have the following $\Omega$-bound for $P_\rho(x,\mathcal {Q})$.
\begin{theorem}\label{thm:r_Q(n)}
	For  $0 \leq \rho < \frac{1}{2}$, we have 
\begin{align}\label{final omega bdd r_Q(n)}
	P_\rho(x,\mathcal {Q})
	=
	\Omega
	\left(
	x^{\theta}
	(\log x)^{\kappa} 
	(\log_2 x)^{\tau'}
	(\log_3 x)^{-\frac{\r}{4}- \frac{3}{8}}\right),
\end{align}
where 
\begin{align}
	\theta = \frac{2 \rho +1}{4},  \:\:
	\kappa = \frac{1- 2 \rho}{4} \: \text{ and } \:  \tau' = \frac{3+2 \rho}{4} \left(2^{\frac{1-2 \rho}{3+ 2 \rho}}-1\right). 
\end{align} 
\end{theorem}

\normalcolor

\subsection{The General Framework} 
\label{sec:notation and result}
In the Riesz sum \eqref{def:A_rho}, the coefficients
$a_n$
will be taken to be coefficients of a Dirichlet series contained in a certain class, defined as follows. 
We consider a class 
$\mathcal{C}$ 
of Dirichlet series satisfying a general functional equation with multiple Gamma factors.
\begin{definition}\label{def:fn eqn of F} 
	Let 
	$(\lambda_n)$ 
	and 
	$(\mu_n)$ 
	be two strictly increasing sequences of positive real numbers tending to infinity and 
	$(a_n)$ 
	and 
	$(b_n)$ 
	be two sequences of complex numbers not identically zero. Consider the functions 
	$\phi(s)$ 
	and 
	$\psi(s)$ 
	representable as Dirichlet series 
	\begin{equation} \label{dirichlet series phi psi}
		\phi(s)=\sum_{n=1}^{\infty} \frac{a_n}{\lambda_n^{s}} \hspace{4mm}\text{and}\hspace{4mm} \psi(s)=\sum_{n=1}^{\infty} \frac{b_n}{\mu_n^{s}},
	\end{equation}
	satisfying the following conditions.
	We say that 
	$\phi(s),\psi(s) \in \mathcal{C}$ if 
	\begin{enumerate}
		\item 
		$	\phi(s)$ 
		and 
		$ \psi(s)$ 
		have finite abscissae of absolute convergence 
		$\sigma_a$ 
		and 
		$\sigma_b$ 
		respectively. 
		\item 
		$\phi(s)$ and $\psi(s)$ 
		satisfy a functional equation  of the type 
		\begin{equation}\label{generalfneqn}
			Q^sF(s)=\omega' Q^{\d-s} \overline{G(\d-\bar{s})},
		\end{equation}
		for some $\delta>0$, 
		with 
		$F(s)= \phi(s)\prod_{i=1}^r \Gamma\left(\alpha_i s+\beta_i \right)$,
		$G(s) = \psi(s)\prod_{i=1}^r \Gamma\left(\alpha_i s+\beta_i \right)$, where 
		$r \in \mathbb{N}$,
		$\alpha_i>0$, 
		$\beta_i \in \mathbb C$  
		with 
		$\Re (\beta_i)\geq 0$,
		$Q>0$ 
		and 
		$\omega' \in \mathbb C$ 
		with 
		$|\omega'|=1$. 
		More precisely, this means that 
		there exists a domain $D$ which is the exterior of a  compact set $S$, such that in $D$, there is  an analytic function $\chi(s)$ with the properties: 
		\begin{enumerate}
			\item $\chi(s)=	Q^sF(s)$ for $\Re(s)> \sigma_a$ and $\chi(s)=	\omega' Q^{\d-s} \overline{G(\d-\bar{s})}$ for $\Re(s)< \delta-\sigma_b$. 
			\item For some constant $\theta<1$, 
			\begin{equation} \label{chigrowthcond}
				\chi(s)=O\left(\exp\left[\exp\left(\frac{\theta \pi|s|}{\s_2-\s_1}\right)\right]\right),
			\end{equation} 
			uniformly in every strip 
			$-\infty<\sigma_1\leq \Re(s) \leq \sigma_2<\infty$. 				
		\end{enumerate}
	\end{enumerate}
\end{definition}		
Given the functional equation \eqref{generalfneqn}, we define
$A \coloneqq \sum_{i=1}^{r}\alpha_i$.	
In the case of functions in the Selberg class,
$2A$ 
is said to be the degree of the function, which is an invariant.  For the purpose of our main theorem, we assume that the $\beta_i$'s appearing in the functional equation \eqref{generalfneqn} are real.

Let 
$a> \max\{0,\s_a,\s_b\}$ 
be sufficiently large so that all the singularities of 
$F(s)= \phi(s)\prod_{i=1}^r \Gamma\left(\alpha_i s+\beta_i \right)$ 
are contained in the strip 
$\d-a<\Re(s)<a$. 
For $\rho \geq 0$,
we  define a residual function $Q_{\rho}$ by 
\begin{align}\label{def:Q_rho}
	Q_\rho(x)
	:= 
	\frac{1}{2 \pi i}\int_{C}\frac{\phi(s)\Gamma(s)x^{s+\rho}}{\Gamma(s+\rho+1)} ds,
\end{align} 
where $C=C_a$ denotes a circle of finite radius, lying inside the strip $\d-a<\Re(s)<a$,  containing all the singularities of $F(s)$, such that all the singularities of $\phi(s)\Gamma(s)x^{s+\rho}$ which lie in this strip are contained inside $C$. It is this function $Q_\rho(x)$ which turns out to be the main term in the evaluation of the Riesz sum \eqref{def:A_rho} (see for instance \cite[Theorem 1.4]{rsv}, stated below as Theorem \ref{thm:FE implies Rieszsum}). In this paper, we will obtain new omega results for the difference $A^{\rho}(x) - Q_\rho(x)$.

In what follows, we will come across functions of the form
\begin{align}\label{conditions on g(x)}
	h(t)
	= 
	\Big(Z(|\log t|)\Big)^{\ell_1} 
	\cdot
	\Big(\log \big(Z(|\log t|)\big) \Big)^{\ell_2}
	\cdots
	\Big(\log_B \big(Z(|\log t|)\big)\Big) ^{\ell_B},
\end{align}
for some fixed $B \in \mathbb{N}$, 
some constants 
$ \ell_1, \ell_2, \hdots,  \ell_B \in \mathbb{R}$, 
such that 
$h(t) \rightarrow \infty$ 
as 
$t \rightarrow \infty$, 
where  $Z(t)$ is a function which satisfies,
\begin{enumerate}
	\item  $Z(ct) \asymp Z(t)$ 
	for any constant 
	$c> 0$. \label{condiion on Z}
	\item  For some $\eta' < \min\left\{\frac{1}{2}, \frac{1}{3A}\right\}$
	and all
	$\tau >0$,
	$$Z(x+\tau) - Z(x) \ll \tau e^{\eta' x},$$
	as $x \rightarrow \infty$. \label{condiion2 on Z}
\end{enumerate}
In particular, \eqref{condiion on Z} means that 
$Z(t)$ has at most polynomial growth and satisfies
$$Z(t) \ll t^{\eta}, \quad \log Z(t) \ll \log t,$$ 
for some 
$\eta > 0$, 
as $t \rightarrow \infty$.

In order to state our main result, we now add certain conditions to our class of Dirichlet series $\mathcal{C}$.
Let 
$\psi(s) \in \mathcal{C}$ 
be as defined in \eqref{dirichlet series phi psi}.
Assume the following additional hypotheses on 
$\psi(s)$.
\begin{itemize}
	\item[\namedlabel{H1}{\bf H1)}] 
	$b_n \geq 0$ for all $n$. 
	
	\item[\namedlabel{H2}{\bf H2)}] The $\beta_i$'s appearing in the functional equation \eqref{generalfneqn} are real.
	
	\item[\namedlabel{H3}{\bf H3)}]  As $x \rightarrow \infty$, we have 
	\begin{align} \label{partial sum b_n}
		\sum_{\mu_n \leq x} b_n \sim c_0 x^{\sigma_b} (\log x)^{\ell-1},
	\end{align}
	where 
	$c_0 > 0$
	and 
	$\ell \in \mathbb{R}$ 
	are some fixed constants and 
	$\sigma_b$
	is the abscissa of absolute convergence of 
	$\psi(s)$.
\end{itemize}

Let 
$Y \geq 1$ 
be  sufficiently large and let 
$c> 1$.
Consider the set
\begin{align}\label{set S(y)}
	S(Y):= \{m \in \mathbb{N}: \mu_m \in [Y,cY]\}.
\end{align}
In what follows, we will choose a subset
$\mathcal{M}_Y$ 
of ``distinguished" integers from $S(Y)$.  A similar set was also considered in $(1.1.8)$  of \cite{Hafner}. 
This means that 
$\mathcal{M}_Y$
is a subset of
$S(Y)$ 
on which the values of $b_n$ are larger than usual, on average. More precisely,  
$\mathcal{M}_Y$ satisfies
\begin{align}\label{condition of H4}
	\frac{1}{|\mathcal{M}_Y|} \sum_{m \in \mathcal{M}_Y} b_m 
	\gg
	G(Y)
	\gg
	\frac{1}{Y} \sum_{m \in S(Y)} b_m,
\end{align}
for some function $G(Y)$.

With these hypotheses in hand, we now state our main result in the general framework. It is worth keeping in  mind that in most examples of interest, the variable $\mu_N$ that appears in the statement of our theorem below is simply $N$. 
\begin{theorem} \label{thm:main theorem}
	Consider 
	$\phi(s), \psi(s) \in \mathcal{C}$,
	as defined in \eqref{dirichlet series phi psi} and recall that
	$A= \sum_{i=1}^{r} \alpha_i$.
	For
	$x>0$, $\rho \geq 0$, 
	consider 
	$A^{\rho}(x),Q_{\r}(x)$ 
	as defined in \eqref{def:A_rho} and \eqref{def:Q_rho} respectively. 
	Suppose that 
	$\phi(s), \psi(s)$
	satisfy hypotheses \ref{H1} to \ref{H3}.
	For  a parameter $N$ chosen suitably (as described below)
	in terms of $x$,
	we consider the set $S(Y)$ as defined in \eqref{set S(y)}, with 
	$$Y= \left(\frac{1}{2}\right)^{2A} \mu_N, \quad c = 3^{2A}.$$
Let
	$\mathcal{M}_Y$ (denoted by $\mathcal{M}_N$ henceforth to ease notation) be 
	any subset of distinguished integers inside 
	$S(Y)$. 
	Let $|\mathcal{M}_N|$ denote the size of 
	$\mathcal{M}_N$ and $G(N)$ denote the function that satisfies \eqref{condition of H4} above for the chosen set $\mathcal{M}_N$. 
	Moreover, the parameter $N=N(x)$ is chosen such that:
	\begin{enumerate}
		\item[\namedlabel{assumption 1}{\rm (1)}] $N \rightarrow \infty$ as $x \rightarrow \infty$, 
		
		\item[\namedlabel{assumption 2}{\rm (2)}]$\mu_N \leq x^{\epsilon A} $ for any $\epsilon > 0$, 
		
		\item[\namedlabel{assumption 3}{\rm (3)}] The function 
		\begin{align} \label{def:g(x)} 
			g(x):= \mu_{N}^{-\frac{\delta}{2}- \frac{\rho}{2A}- \frac{1}{4A}} 
			|\mathcal{M}_N| G({N})
		\end{align} 
		satisfies  $g(x) \asymp h(x)$ as $x\to \infty$, for some function $h(x)$ of the form \eqref{conditions on g(x)},  and 
		\item[\namedlabel{assumption 4}{\rm (4)}]  
		 $2^{|\mathcal{M}_N|} \leq  x^{\epsilon}$
		 as 
		 $x \rightarrow \infty$,
		 for some
		 $0 < \epsilon < \frac{1}{6}$. 
		\normalcolor
	\end{enumerate}
	Then for 
	$0 \leq \r <  (2\s_b-\d)A-\frac{1}{2}$, we have
	\begin{align} \label{omega bd}
		A^{\rho}(x)-Q_{\r}(x) = \Omega\left(x^{\frac{\delta}{2} - \frac{\rho}{2A}- \frac{1}{4A} + \rho } g(x) \right).
	\end{align}
	Moreover, writing
	$\mu=\frac{1}{2}+\sum_{i=1}^{r}\left(\beta_i-\frac{1}{2}\right)$ 
	and 
	$\gamma=-\left(\frac{\d}{2}+\frac{\r}{2A}+\frac{1}{4A}\right)$,
	we have 
	\begin{align}\label{special case 1 thm}
			\Re \bigg(\overline{\omega'}\big(A^{\rho}(x)-Q_{\r}(x)\big)\bigg) \normalcolor
		= 
		\Omega_+\left(x^{\frac{\delta}{2} - \frac{\rho}{2A}- \frac{1}{4A} + \rho } g(x) \right),
	\end{align}
	whenever 
	$A\gamma - \mu \equiv 0 \pmod {2 } $ 
	and 
	\begin{align}\label{special case 2 thm}
		\Re \bigg(\overline{\omega'}\big(A^{\rho}(x)-Q_{\r}(x)\big)\bigg) \normalcolor
		= 
		\Omega_-\left(x^{\frac{\delta}{2} - \frac{\rho}{2A}- \frac{1}{4A} + \rho  } g(x) \right),
	\end{align}
	whenever
	$A\gamma - \mu \equiv 1 \pmod {2}$.
\end{theorem}

\section{Connections to earlier work}

\subsection{Comparison to known Omega results}
\label{sec: comparison}
When $\rho =0$,  Theorem \ref{thm: piltz} allows us to  recover the bound given by Karak and Mahatab \cite[Theorem 1.1]{karakmahatab}.   
Our  result  sharpens  the bounds obtained by Hafner \cite{hafner1983} for  Riesz sums associated with the Piltz divisor problem. Hafner \cite[Theorem 2]{hafner1983} obtained that if
$0 \leq \rho <  \frac{km-1}{2} $,
then there exists a positive number $c$ depending on 
$\mathbb{K}$, $k$ and $\rho$, such that 
\begin{align} \label{hafner piltz bdd}
	\Delta_{\mathbb{K}, k,\rho}(x)
	= 
	\Omega_* \big(x^\theta (\log x)^\kappa (\log_2 x)^\tau \exp(-c\sqrt{\log_3 x})\big),
\end{align}
where
$\theta$
and 
$\kappa$ are as in \eqref{tau' in piltz},  
\begin{align}\label{tau in piltz}
	\tau = \frac{km-1- 2 \rho}{2km}\big(k \log k + k \sum_{v=1}^{m} \delta_v v \log v-k +1\big)+k-1, 
\end{align}
and 
\begin{align}
	\Omega_*
	=
	\begin{cases}
		\begin{array}{ll}
			\Omega_+ 	
			& \text{if}  \quad kr_2 \equiv 0 \pmod 4  \hspace{1mm} \text{and} \hspace{1mm} \rho \geq \frac{km-3}{2},
			\\
			\Omega_+ 	
			& \text{if}  \quad kr_2 \equiv 3 \pmod 4  \hspace{1mm} \text{and} \hspace{1mm} \rho \geq \frac{km-2}{2},
			\\
			\Omega_-	
			& \text{if}  \quad kr_2 \equiv 2 \pmod 4  \hspace{1mm} \text{and} \hspace{1mm} \rho \geq \frac{km-3}{2},
			\\
			\Omega_-	
			& \text{if}  \quad kr_2 \equiv 1 \pmod 4  \hspace{1mm} \text{and} \hspace{1mm} \rho \geq \frac{km-2}{2},
			\\
			\Omega_\pm 
			& \text{otherwise}.
		\end{array}
	\end{cases}
\end{align}
We not only improve \eqref{hafner piltz bdd} due to the 
$\log_3x$  
term in \eqref{final piltz div bdd}, but also since our exponent of 
$\log_2x$
in \eqref{final piltz div bdd} can be shown to be larger than  the exponent 
$\tau$ 
of 
$\log_2x$  
in \eqref{hafner piltz bdd}  in the regime
$0 \leq \rho < \frac{km-1}{2}$ (see Lemma \ref{lem: tau < tau'}).

The bound given by Corollary \ref{cor: d_k} for the Piltz divisor function sharpens a result of  Hafner \cite[p.~71]{Hafner}, which  showed that 
\begin{align} \label{hafner omega bdd for d_kn}
	\Delta_{k,\rho}(x)  =
	\begin{cases}
		\begin{array}{ll}
			\Omega_\pm \big(x^\theta (\log x)^\kappa (\log_2 x)^\tau \exp(-c\sqrt{\log_3 x})\big) ,	
			& \text{if}  \quad \rho < \frac{k-3}{2},
			\\
			\Omega_+ \big(x^\theta (\log x)^\kappa (\log_2 x)^\tau \exp(-c\sqrt{\log_3 x})\big),		
			& 
			\text{if} \quad 
			\frac{k-3}{2} \leq \rho < \frac{k-1}{2},
		\end{array}
	\end{cases}
\end{align}
where $c$ is a positive constant depending on $k$ and $\rho$; 
$\theta$ and $\kappa$ are the same as in \eqref{final omega bdd d_k} and 
\begin{align} 
	\tau = \frac{k-1- 2 \rho}{2k}(k \log k -k +1)+k-1,
\end{align} 
is smaller than $\tau'$ in the domain 
$0 \le \rho < \frac{k-1}{2}$.  Putting $\rho=0$ in  Corollary \ref{cor: d_k} reduces \eqref{final omega bdd d_k} to the $\Omega$-bound obtained recently by Mahatab \cite[Theorem~3]{mahatab2025} and Sourmelidis \cite[Corollary~3.1]{sourmelidis2024}. Similarly, the bound \eqref{omega bdd for div fn} for the classical divisor function coincides with  \cite[Theorem~1]{mahatab2025} and  \cite[Corollary~1.3]{sourmelidis2024} when $\rho=0$. 

Theorem \ref{thm: r(n)} improves upon a result of  Hafner \cite[p. 74]{Hafner}, who showed that for $0 \leq \rho < \frac{1}{2}$, 
\begin{align} \label{hafner omega bdd for r(n)}
	P_{\rho}(x)  = \Omega_- \big(x^\theta (\log x)^\kappa (\log_2 x)^\tau \exp(-c\sqrt{\log_3 x})\big),
\end{align}
where $c$ is a positive constant depending on  $\rho$; 
$\theta$ and $\kappa$ are the same as in \eqref{final omega bdd r(n)} and 
\begin{align} 
	\tau = \frac{(1- 2 \rho)(\log 2)}{4},
\end{align} 
can be seen to be less than the exponent 
$\tau'$ 
of 
$\log_2x$ 
in \eqref{final omega bdd r(n)}  for
$0 \leq \rho < \frac{1}{2}$. Again, in the $\rho=0$ case, Theorem \ref{thm: r(n)} yields \cite[Theorem~2]{mahatab2025} and  \cite[Corollary~1.3]{sourmelidis2024}. 

In the special case $\rho=0$, Theorem~\ref{thm: r_3(n)} recovers a classical result of Tsang~\cite[Theorem~1]{TsangSphere} for the sphere problem, albeit without retaining the directional  sign information of the $\Omega$-bounds. More generally, in their foundational work, Chandrasekharan and Narasimhan~\cite[p.~135]{ChandrasekharanV} proved that for $0 \leq \rho \leq 1$, 
\begin{align} \label{CN lim result}
	\limsup_{x \rightarrow \infty} \frac{P_{3,\rho}(x)}{x^{\frac{1+ \rho}{2}}} = +\infty \quad \text{and} \quad \liminf_{x \rightarrow \infty} \frac{P_{3,\rho}(x)}{x^{\frac{1+ \rho}{2}}} = -\infty.
\end{align}
Theorem~\ref{thm: r_3(n)} thus provides a substantial quantitative strengthening of these classical bounds across the entire range of $\rho$.

Theorem \ref{thm:r_Q(n)}   gives a new result for the lattice point discrepancy associated with an ellipse. 
In the special case $\rho=0$, it improves the result of Nowak \cite{Nowak1985}. Nowak proved an omega bound for the lattice point discrepancy of
compact convex planar domains $\mathcal D$ containing the origin as an interior point, whose boundary $\partial\mathcal D$ is the image of the unit circle under a one-to-one $C^\infty$ map and has finite non-vanishing curvature throughout. For the ellipse
$u^2+5v^2\leq x$, his result implies that
\begin{align}
	P_0(x,\mathcal Q)
	=
	\Omega_-\!\left(x^{1/4}(\log x)^{1/4}\right).
\end{align}
Thus, Theorem \ref{thm:r_Q(n)} improves Nowak's bound in the case $\rho=0$, while also providing new omega bounds for the corresponding Riesz sums throughout the range
$0 \leq \rho < \frac{1}{2}$.
\normalcolor


\subsection{An overview of Omega results and Riesz smoothing}
A decisive breakthrough in the study of lower bounds for error terms occurred with the work of Vorono\"i \cite{VoronoiI, VoronoiII} (1903–1904), who derived explicit formulae expressing $\Delta_\rho(x)$ and $P_\rho(x)$ in terms of Bessel functions. This Vorono\"i summation formula provided a bridge between discrete divisor sums and continuous special function expansions, furnishing the analytic framework that still underpins modern investigations. Building on this foundation, for the classical case $\rho=0$, G. H. Hardy (1916) established 
that
\begin{align}
	\Delta_0(x)=
	  \begin{cases}
		\Omega_+ \left( (x \log x)^{1/4} \log_2 x\right),
		\\
		\Omega_-\left( x^{1/4}\right),
	\end{cases}
\end{align}
with analogous bounds for $P_0(x)$.  These results demonstrated that the error terms attain magnitudes of order at least $x^{1/4}$, thereby showing that the classical $O(\sqrt{x})$ bounds of Gauss and Dirichlet could not approach the conjectural barrier from below. Hardy further identified a fundamental transition at $\rho=\tfrac12$: for $\rho>1/2$, the associated Vorono\"i series converges absolutely and the oscillatory structure becomes more regular, whereas for $0\le\rho<1/2$ it remains only conditionally convergent and is markedly volatile.

Subsequent decades saw efforts to generalize Hardy’s results to arbitrary $\rho\ge 0$.  Early extensions by  Landau and Cramér laid important groundwork, but a comprehensive understanding of the critical range $0\le\rho<1/2$ emerged only later. In a sequence of influential works culminating in 1982, Hafner \cite{Hafner, hafner1983} established sharp $\Omega_+$ and $\Omega_-$ bounds for $\Delta_\rho(x)$ and $P_\rho(x)$, proving that in this regime the error term attains a magnitude of
$
x^{\frac{1}{4}+\frac{\rho}{2}}
$
up to explicit logarithmic factors. Hafner’s results confirmed that smoothing modifies the frequency profile of the oscillations but preserves the structural amplitude dictated by the Vorono\"i expansion.

More recently, attention has turned to refining the logarithmic components of these lower bounds. A major advance was achieved by Soundararajan \cite{Sound}, who introduced the resonance method in 2003. This technique constructs Dirichlet polynomials that “resonate’’ with the oscillatory terms in the Vorono\"i expansion, leading to significantly strengthened $\Omega$-results for the classical case $\rho=0$, including additional powers of iterated logarithms. In particular, Soundararajan obtained 
\begin{align}
	\Delta_0 (x)
	= \Omega \left(
	(x \log x)^{1/4}
	(\log_2 x)^{(3/4)(2^{4/3}-1)}
	(\log_3 x)^{-5/8}
	\right),\label{sound div bdd}
	\\
	P_0(x)
	= \Omega \left(
	(x \log x)^{1/4}
	(\log_2 x)^{(3/4)(2^{1/3}-1)}
	(\log_3 x)^{-5/8}
	\right).
\end{align}
Recently, these bounds were sharpened by Mahatab \cite[Theorem 1]{mahatab2025} and Sourmelidis
\cite[Corollary 1.3]{sourmelidis2024}  
independently, by further refining the resonance method, to obtain 
\begin{align}
	\Delta_0(x)
	= \Omega \left(
	(x \log x)^{1/4}
	(\log_2 x)^{(3/4)(2^{4/3}-1)}
	(\log_3 x)^{-3/8}
	\right),
	\\
	P_0(x)
	= \Omega \left(
	(x \log x)^{1/4}
	(\log_2 x)^{(3/4)(2^{1/3}-1)}
	(\log_3 x)^{-3/8}
	\right).
\end{align}

\normalcolor

 Riesz sums were first studied systematically by Chandrasekharan and Narasimhan \cite{ChandrasekharanV}, who established \(\Omega\)-results for the corresponding error terms and showed their optimality for large values of \(\rho\). For smaller values of \(\rho\), however, these bounds admit further refinement, as demonstrated in subsequent work (see, for example, \cite{BerndtAvgI, BerndtAvgII, Hafner}). In particular, Berndt improved the \(\rho=0\) case in \cite{BerndtAvgI, BerndtAvgII} and extended the analysis to \(\rho \ge 0\) in \cite{BerndtAvg}. Later, Hafner \cite{Hafner} considered a broader class of functions and obtained \(\Omega_\pm\)-results by adapting ideas of Szeg{\"o} \cite{szego1926I} and Szeg{\"o}-Walfisz \cite{SzegoWalI, SzegoWalII}. Indeed, Hafner investigated the same class $\mathcal{C}$ of Dirichlet series (cf. Definition \ref{def:fn eqn of F}), with the sole exception that he required a slightly more restrictive growth condition than our condition \eqref{chigrowthcond}. Specifically, he assumed that
 \begin{align}
 	\lim_{|\Im(s)| \rightarrow \infty} \chi(s) = 0,
 \end{align}
 uniformly in any vertical strip $-\infty < \sigma_1 \leq \Re(s) \leq \sigma_2 < \infty$.
 
In the present work, we strengthen Hafner’s results by employing the resonance method, which applies to trigonometric series of the form
\[
\sum_n f(n)\cos(2\pi \gamma_n x + \beta),
\]
where \((\gamma_n)\) is a nondecreasing sequence of nonnegative real numbers, \(\beta \in \mathbb{R}\), and the coefficients $f(n)$ are nonnegative.

\normalcolor

     \section{Preliminary Lemmas}
   We record Stirling's asymptotic formulae \cite[Corollary 1.4.3]{AndrewsAskeyRoy1999} for the Gamma function below.
\begin{lemma} \label{lem:stirling}
	For 
	$\epsilon>0$ 
	and 
	$|\arg s| \leq \pi -\epsilon$, 
	we have as 
	$s\rightarrow \infty$, 
	\begin{equation} \label{stirling gen formula}
		\G(s)\sim \sqrt{2\pi} e^{-s}s^{s-\frac{1}{2}}. 
	\end{equation}
\end{lemma}
Note that for $s=\sigma+it$, if  $\sigma$ is fixed and 	 $|t|\rightarrow \infty$, this yields 
\begin{equation} \label{stirling's formula}
	|\Gamma(s)|\sim e^{- \frac{\pi}{2}|t|}|t|^{\sigma-\frac{1}{2}}\sqrt{2\pi}. 
\end{equation} 
If we put
$t=0$
in \eqref{stirling gen formula}, as $\sigma \rightarrow \infty$ we have
$\G(\sigma)\sim \sqrt{2\pi} e^{-\sigma}\sigma^{\sigma-\frac{1}{2}}$, which gives
\begin{align}\label{stirling for real}
	\G(\sigma+1)\sim \sqrt{2\pi} e^{-\sigma}\sigma^{\sigma+\frac{1}{2}},
\end{align}
upon using $\Gamma(s+1) = s \Gamma(s)$.
From \eqref{stirling for real} one can derive that for 
$x, \alpha \in \mathbb{R}$, we have
\begin{align} \label{gamma of x plus alpha}
	\Gamma(x+\alpha) \sim x^{\alpha} \Gamma(x),
\end{align}
as $x \rightarrow \infty$.

     
 We will require the following lemma, which is analogous to a result of Hafner \cite[Lemma 2.3.6]{Hafner}. 
	\begin{lemma}\label{lem: asymp for I_k,u}
	For $u>0$ and $x$ sufficiently large, let
	\begin{align}
		I_{x,u}
		=
		\frac{1}{\Gamma(x+1)} 
		\int_{0}^{\infty}
		(ut^A)^{\theta} h(ut^A)
		e^{-t} t^x  dt,
	\end{align}
	where $\theta \in \mathbb{R}$,  $A= \sum_{i=1}^{r}\alpha_i$ and $h(t)$ is a function of the form \eqref{conditions on g(x)}. 
	Fix
	$0 < \epsilon < \frac{A}{2} $ 
	and let
	$u_0$
	be sufficiently small.
	Then as 
	$x \to \infty$, we have 
	\begin{align}
		I_{x,u} 
		\sim
		c 
		u^{\theta}
		x^{ A \theta}
		h(ux^A),
	\end{align}
	for some constant 
	$c>0$,
	provided
	$x^{-\epsilon} \leq u \leq u_0$.
\end{lemma}
The proof of Lemma \ref{lem: asymp for I_k,u} follows exactly in the same way as the proof of Lemma 2.3.6 in \cite{Hafner}. For the sake of completeness, we provide a sketch of the proof, with suitable modifications corresponding to our choice of $h(t)$. 
\begin{proof}
	We change the variable $t$ to $xt$  in the integral in 
	$I_{x,u}$ to get
	\begin{align}\label{1st int in hafners lem}
		I_{x,u}
		=
		\frac{u^{\theta}x^{x+A \theta+1}e^{-x}}{\Gamma(x+1)} 
		\int_{0}^{\infty}
		t^{A\theta} h(ux^At^A)
		(e^{1-t} t)^x  dt,
	\end{align}
	after some simplification.
	Let
	$\frac{1}{3} < \tau < \frac{1}{2}$,
	$y= x^{-\tau}$
	and 
	partition the interval
	$(0, \infty)$
	into 
	$(0, 1-y], (1-y,1+y), [1+y, \infty)$.
	Denote the integral on the right hand side of \eqref{1st int in hafners lem} by 
	$I^*_{x,u}$ 
	and write
	\begin{align}
		I^*_{x,u}
		&=
		\left[
		\int_{1-y}^{1+y}
		+
		\int_{1+y}^{\infty}
		+
		\int_{0}^{1-y}
		\right]
		t^{A\theta} h(ux^At^A)
		(e^{1-t} t)^x  dt
		\\
		&=J_1+J_2+ J_3 \quad \text{(say).} 
	\end{align}
	We proceed to estimate the integral $J_1$ by estimating the integrand in the interval $(1-y, 1+y)$. From  \cite[equation (2.3.17)-(2.3.18)]{Hafner}, for 
	$t \in (1-y,1+y)$,
	we have
	\begin{align}\label{est of t^Atheta}
		t^{A\theta} = 1+ O(y)
	\end{align}
	and 
	\begin{align}
		(e^{1-t}t)^x= 
		e^{-\frac{x(1-t)^2}{2}}\left(1+o(1)\right),
	\end{align}
	as 
	$x \rightarrow \infty$.
	Similar to \eqref{est of t^Atheta}, we have 
	$t^{A} = 1+ O(y)$ and hence we get
	\begin{align}
		\log(ux^At^A)
		&= \log(ux^A)+ \log(1+ O(y))
		\\
		&= \log(ux^A)+O(y), \label{est of logux^At^A}
	\end{align}
	since 
	$y \rightarrow 0$ 
	as $x \rightarrow \infty$.
	Using \eqref{est of logux^At^A} and condition \eqref{condiion2 on Z} on 
	$Z(t)$, we obtain 
	\begin{align}
		Z(|\log(ux^At^A)|)
		&= Z(\log(ux^A))+ O(ye^{\eta'\log (ux^A)})
		\\
		&= Z(\log(ux^A))+ O(yu^{\eta'}x^{\eta' A})
		\\
		&=Z(\log(ux^A))\left(1+O\Big(\frac{x^{-\tau+\eta' A}}{Z(\log(ux^A))}\Big)\right), \label{est of Z(log)}
	\end{align}
	putting 
	$y= x^{-\tau}$ and noting that 
	$u \leq u_0$.
	Using  
	$\eta' < \frac{1}{3A}$, 
	$\tau > \frac{1}{3}$ and 
	$Z(\log(ux^A)) \rightarrow \infty$ 
	as 
	$x \rightarrow \infty$
	in \eqref{est of Z(log)}, we get
	\begin{align}
		Z(|\log(ux^At^A)|) = Z(\log(ux^A))(1+o(1)),\label{final est of Z(log)}
	\end{align}
	which gives
	\begin{align}
		\big(Z(|\log(ux^At^A)|)\big)^{\ell_1}=\big(Z(\log(ux^A))\big)^{\ell_1} (1+o(1)),
	\end{align}
	for any 
	$\ell_1 \in \mathbb{R}$. 
	Taking logarithm on both sides of \eqref{final est of Z(log)}, after some  simplification, we see that
	\begin{align}
		\Big(\log Z(|\log(ux^At^A)|)\Big)^{\ell_2}= \Big(\log\big(Z(\log(ux^A))\big)\Big)^{\ell_2}(1+o(1)),\label{final est log Z}
	\end{align}
	for any 
	$\ell_2 \in \mathbb{R}$. Proceeding in a similar manner, we  obtain
	\begin{align}
		\Big(\log_B Z(|\log(ux^At^A)|)\Big)^{\ell_B}= \Big(\log_B\big(Z(\log(ux^A))\big)\Big)^{\ell_B}(1+o(1)),\label{final est log_B Z}
	\end{align}
	for any
	$B \in \mathbb{N}$
	and 
	$\ell_B \in \mathbb{R}$. 
	Recalling the form of $h(t)$ in \eqref{conditions on g(x)} and combining \eqref{final est of Z(log)}, \eqref{final est log Z}, \eqref{final est log_B Z}, we conclude that
	\begin{align}
		h(ux^At^A)= h(ux^A)(1+o(1)),
	\end{align}
	for $t \in (1-y, 1+y)$.
	Following Hafner's proof of Lemma 2.3.6 after equation (2.3.19), we obtain for some constant 
	$c>0$,
	\begin{align}\label{est of J1}
		J_1= cx^{-\frac{1}{2}} h(ux^A)(1+o(1)), \quad \text{and} \quad  J_2, J_3=o\big(x^{-\frac{1}{2}} h(ux^A)\big).
	\end{align}
	Combining the above estimates for 
	$J_1, J_2, J_3$,
	we obtain
	\begin{align}
		I^*_{x,u}= cx^{-\frac{1}{2}} h(ux^A)(1+o(1)). \label{final est I^*}
	\end{align}
	Substituting \eqref{final est I^*} in \eqref{1st int in hafners lem} and using Stirling's formula (see \eqref{stirling for real}), we conclude that
	\begin{align}
		I_{x,u} 
		\sim
		c 
		u^{\theta}
		x^{ A \theta}
		h(ux^A),
	\end{align}
	where $c>0$ is some constant.
\end{proof}
 
 We will rely on the following result by Roy, Sahoo, and Vatwani \cite[Theorem 1.4]{rsv}, which gives a Riesz sum identity for a class of arithmetic functions in the class $\mathcal{C}$.
 \begin{theorem} \label{thm:FE implies Rieszsum}
 	Let 
 	$\phi(s)$ 
 	and 
 	$\psi(s)$ 
 	be two Dirichlet series satisfying Definition \ref{def:fn eqn of F}. Let 
 	$a> \max\{0,\s_a,\s_b\}$ 
 	be sufficiently large so that all the singularities of $F(s)= \phi(s)\prod_{i=1}^r \Gamma\left(\alpha_i s+\beta_i \right)$ 
 	are contained in the strip 
 	$\d-a<\Re(s)<a$. 
 	Then the functional equation  
 	$$Q^sF (s)=\omega' Q^{\delta-s}\overline{G(\delta-\bar{s})},$$ 
 	where  
 	$G(s) = \psi(s)\prod_{i=1}^r \Gamma\left(\alpha_i s+\beta_i \right)$, 
 	implies the identity
 	\begin{multline}\label{riesz sum}
 		\frac{1}{\Gamma(\rho+1)}\sideset{}{'}\sum_{\lambda_n \leq x} a_{n}(x-\lambda_n)^\rho
 		\\
 		=
 		Q_\rho(x)
 		+
 		\bigg(\frac{\omega' x^{\delta+\rho}}{Q^\delta}\bigg) \sum_{n=1}^{\infty}  \bar{b}_n
 		\frac{1}{2 \pi i}
 		\int_{a-i\infty}^{a+i\infty} 
 		\frac{\Gamma(\delta-s) \prod_{i=1}^{r}\Gamma(\alpha_i s+\bar{\beta}_i)}{\Gamma(1+\delta-s+\rho)\prod_{i=1}^{r}\Gamma(\alpha_i(\delta-s)+\beta_i)}\bigg(\frac{\mu_nx}{Q^2}\bigg)^{-s}ds	, 
 	\end{multline}
 	for $x>0$, and $\r> (2\s_b-\d)A-\frac{1}{2}$. Here,  $A=\sum_{i=1}^{r}\a_i$ and $Q_\rho(x)$ is the residual function defined in \eqref{def:Q_rho}.
 	
 \end{theorem}

Recall the setup of Definition \ref{def:fn eqn of F}. 
For $R>0$ and $a,b\in\mathbb{R}$, let $C_R(a,b)$ denote the contour consisting of the line 
$(a-i\infty, a-iR)$, 
followed by the three sides 
$[a-iR, b-iR]$, 
$[b-iR, b+iR]$, 
$[b+iR, a+iR]$ 
of a rectangle, followed by the line 
$(a+iR, a+i\infty)$.
For  $x>0, \rho \geq 0$
and any non-negative integer  $m$, 
define
$I_{\rho+m}(x)$ 
to be the integral
\begin{align}\label{def:I_rho*}
	I_{\rho+m}(x)
	:= 
	\frac{1}{2 \pi i}
	\int_{C_R(c_\rho, c_\rho+U)} \frac{\Gamma(\delta-s) \prod_{i=1}^{r}\Gamma(\alpha_is+\bar{\beta}_i)}{\Gamma(1+\delta-s+\rho+m)\prod_{i=1}^{r}\Gamma(\alpha_i(\delta-s)+\beta_i)}
	x^{\r+m+\d-s}ds,
\end{align}  
where
\begin{align} \label{c rho}
&c_\r=\frac{A \d+\r}{2A}-\xi, \quad 0<\xi<\min\left\{\frac{1}{4A}, \frac{\delta}{2}\right\},
\\
&R>
\max\left\{
\frac{|\sum_{i=1}^{r} \Im \beta_i |}{A},
\max_{1\le i\le r}
\frac{|\operatorname{Im}\beta_i|}{\alpha_i}
\right\},
\end{align}
and 
$U>0, R$ 
are taken to be such that the contour 
$C_R(c_\rho, c_\rho+U)$ 
does not pass through any singularities of the integrand. \normalcolor
For  $a>0$, we define
$I_{\rho+m}^{(a)}(x)$ 
to be the integral 
\begin{align}\label{def:I_rho}
	I_{\rho+m}^{(a)}(x)
	:=	
	\frac{1}{2 \pi i}
	\int_{a-i\infty}^{a+i\infty} 
	\frac{\Gamma(\delta-s) \prod_{i=1}^{r}\Gamma(\alpha_is+\bar{\beta}_i)}{\Gamma(1+\delta-s+\rho+m)\prod_{i=1}^{r}\Gamma(\alpha_i(\delta-s)+\beta_i)}
	x^{\r+m+\d-s}ds,
\end{align}
such that the line
$\Re(s)=a$ 
does not pass through any singularities of the integrand. Moreover, 
$a, U, R$ 
are chosen in 
\eqref{def:I_rho*} and \eqref{def:I_rho}
in such a way that the region between 
$C_R(c_\rho, c_\rho+U)$ 
and 
$\Re(s)=a$
does not contain any singularities of the integrand.
Note that the path of integration in each of the integrals above is independent of $m$.
The next few lemmas will be concerned with these integrals.

\begin{lemma} \label{lem:I_rho to I*_rho}
     For 
     $x > 0$, 
     $\rho \geq 0$ and any non-negative integer  $m$, 
     consider 
     $I_{\r+m}(x)$,
     $I_{\r+m}^{(a)}(x)$
      defined above. Then for  any 
     $m > (2a-\d)A-\rho - 1$, 
     we have
     \begin{align}
     	 I_{\r+m}^{(a)}(x)= I_{\r+m}(x).
     \end{align}
\end{lemma}

\begin{proof}
		Let us define 
	\begin{equation}\label{def of G_rho}
		G_\r(s)
		:=
		\frac
		{\Gamma(\delta-s) \prod_{j=1}^{r}\Gamma(\alpha_js+\bar{\beta}_j)}
		{\Gamma(1+\delta-s+\rho)\prod_{j=1}^{r}\Gamma(\alpha_j(\delta-s)+\beta_j)}.
	\end{equation}
Recalling the definitions \eqref{def:I_rho*} and \eqref{def:I_rho} of 
$I_{\r+m}(x)$ 
and 
$ I^{(a)}_{\r+m}(x)$ respectively, 
we need to prove that
\begin{equation} \label{eq:G_rho equality}
	\frac{1}{2 \pi i}
	\int_{C_R(c_\rho, c_\rho+U)} 
	G_{\r+m}(s) 
	x^{\r+m+\d-s}ds
	=
	\frac{1}{2 \pi i}
	\int_{a-i\infty}^{a+i \infty}
	G_{\r+m}(s) 
	x^{\r+m+\d-s}ds,
\end{equation} 
for sufficiently large $a$. Take $T$ sufficiently large.
Consider the closed contour 
$\mathcal{C}$
consisting of the line segments 
$[c_{\rho}-iT, a-iT]$, 
$[a-iT,a+iT]$,
$ [a+iT, c_{\rho}+iT]$,
$ [c_{\rho}+iT, c_{\rho}+iR ]$,
$ [c_{\rho}+iR, c_{\rho}+U+iR]$,
$[c_{\rho}+U+iR, c_{\rho}+U-iR]$, 
$[c_{\rho}+U-iR, c_{\rho}-iR]$, 
$[c_{\rho}-iR, c_{\rho}-iT]$.
By construction, this contour does not contain any singularities of 
$G_{\r+m}(s)$. 
Hence by Cauchy's residue theorem, we have
\begin{align} \label{eq:int of G_rho on contour}
    \int_{\mathcal{C}}
	G_{\r+m}(s) x^{\r+m+\d-s}ds
	=0.
\end{align} 
The proof will be complete if we show that the horizontal integrals 
$\int_{c_{\rho}-iT}^{ a-iT}	G_{\r+m}(s) x^{\r+m+\d-s}ds $
and 
$\int_{a+iT}^{c_{\rho}+iT} 	G_{\r+m}(s) x^{\r+m+\d-s}ds$
go to zero as
$T \rightarrow \infty$.
Using the Stirling formula \eqref{stirling's formula},   after some simplification we obtain  the estimate 
\begin{equation}
	G_{\r+m}(s) 	\ll T^{-1-\r-m+(2a-\d)A},
\end{equation} 
for the integrand. 
Since this is $o(1)$  as $T\rightarrow \infty$ whenever $m>(2a-\d)A-\rho - 1$, we obtain \eqref{eq:G_rho equality}.
\end{proof}

\begin{lemma} \label{lem:d/dx I_rho}
	For 
 $x > 0$, 
$\rho \geq 0$ and any non-negative integer  $m$, 
	let 
	$I_{\r+m}(x)$ 
	be the integral defined in \eqref{def:I_rho*}. Then  we have
	\begin{align}
		\frac{d}{dx} I_{\r+m+1}(x) = I_{\r+m}(x).
	\end{align}
\end{lemma}
\begin{proof}
Using \cite[Theorem 2.3]{temmebook}, we may differentiate under the integral sign to get
	\begin{align}
		\frac{d}{dx} 
		I_{\r+m+1}(x)
		=
		\frac{1}{2 \pi i}
		\int\limits_{C_R(c_\rho, c_\rho+U)}\!\! \frac{\Gamma(\delta-s) \prod_{i=1}^{r}\Gamma(\alpha_is+\bar{\beta}_i)}{\Gamma(1+\delta-s+\rho+m+1)\prod_{i=1}^{r}\Gamma(\alpha_i(\delta-s)+\beta_i)}
		\left(\frac{d}{dx} x^{\r+m+1+\d-s} \right)ds,
	\end{align}
as 
$I_{\r+m+1}(x)$ is
uniformly convergent for
$m \geq 0$
and 
$x$
in any compact set. 
  Using the property that $\Gamma(z+1)= z\Gamma(z)$, we obtain from above
 \begin{align}
 \frac{d}{dx} 
 I_{\r+m+1}(x)
 =&
 \frac{1}{2 \pi i}
 \int_{C_R(c_\rho, c_\rho+U)} \frac{\Gamma(\delta-s) \prod_{i=1}^{r}\Gamma(\alpha_is+\bar{\beta}_i)}{\Gamma(1+\delta-s+\rho+m)\prod_{i=1}^{r}\Gamma(\alpha_i(\delta-s)+\beta_i)}
  x^{\r+m+\d-s}ds
 \\
 =&
 I_{\r+m}(x),
 \end{align}
as the curve of integration is independent of $m$.
\end{proof}
The following lemma gives an asymptotic formula for the integral $I_\r(x)$, defined by letting $m=0$ in \eqref{def:I_rho*}.

\begin{lemma} \label{lem:asym for I_rho} 
	Let $\rho\geq0$ and $m$ be a non-negative integer. For a fixed $m$, define
	$$
	d_m:=c_\rho+\frac{m+1}{2A},
	$$
	and suppose that $C_R(c_\rho, d_m)$
	does not pass through any poles of the integrand.
	If
	$c_\rho +U < d_m$,
	  we have as $x\to\infty$,
	\begin{align}\label{formula for I_rho}
		I_\rho(x)
		&=
		\sum_{n=0}^{m}	A_nx^{\frac{\omega-n-2i\nu-\frac12}{2A}}
		\cos\left((h^{-1}x)^{\frac1{2A}}
		+\pi\left(A\gamma+\frac n2+i\nu-\mu\right)\right) 
		\\
		&+ O_m\left(x^{\frac{\omega-m-1}{2A}+\xi}\right)
		+O_m\left(x^{\rho+\delta-c_\rho-U}\right),
	\end{align}
where  $ A=\sum_{i=1}^{r}\a_i$,
$ A_n \in \mathbb{C}$ 
are some constants, and writing $\b_i=b_i+ic_i$, we have
\begin{align} 
	&\nu=\sum_{i=1}^{r}c_i, \hspace{0.5mm} \omega=A\d+(2A-1)\r,\hspace{0.5mm} h=\exp \left( \sum_{i=1}^{r}2\a_i\log\a_i-2A \log2A\right),\hspace{0.5mm} \mu=\frac{1}{2}+\sum_{i=1}^{r}\left(\beta_i-\frac{1}{2}\right),
	\\
	& \gamma=-\left(\frac{\d}{2}+\frac{\r}{2A}+\frac{1}{4A}\right). 
	\label{notations in lem for I_rho}
\end{align}
 If $c_\rho+ U \geq d_m$, then \eqref{formula for I_rho} holds without the
	second error term.
	\end{lemma}

\begin{proof}
We follow the proof of Lemma 2.11 in \cite{rsv}, with a modification in the
contour-shifting argument. Let
\begin{align}
	G_\r(s)
	:=
	\frac
	{\Gamma(\delta-s) \prod_{j=1}^{r}\Gamma(\alpha_js+\bar{\beta}_j)}
	{\Gamma(1+\delta-s+\rho)\prod_{j=1}^{r}\Gamma(\alpha_j(\delta-s)+\beta_j)}.
\end{align}
Following the arguments up to (2.23) of \cite{rsv}, we write
\begin{align}
	G_\rho(s)=\sum_{n=0}^{m}F_n^\rho(s)+H_m(s),
\end{align}
where $F_n^\rho(s)$ are defined in (2.20) of \cite{rsv}, and
\begin{align}
	H_m(s)\ll_m\frac{|F_0^\rho(s)|}{|s|^{m+1}}
\end{align}
as $|\operatorname{Im}(s)|\to\infty$, uniformly in any fixed bounded
vertical strip. 

We first consider the case $c_\rho+U<d_m$ and move the contour of integration from $C_R(c_\rho, c_\rho+U)$  
to $C_R(c_\rho, d_m)$. 
Since there are only finitely many simple poles of 
$G_\rho(s)$ between these two contours, by Cauchy's residue theorem,
\begin{align}
	I_\rho(x)
	&=\frac{1}{2\pi i}\int_{C_R(c_\rho, c_\rho+U)}
	G_\rho(s)x^{\rho+\delta-s}\,ds
	\\
	&=\frac{1}{2\pi i}\int_{C_R(c_\rho, d_m)}
	G_\rho(s)x^{\rho+\delta-s}\,ds
	-
	\sum_{\substack{k\geq 0 \\ c_\rho+U < \delta+k < d_m}}
	\Res_{s=\delta+k}	 	G_\rho(s)x^{\rho+\delta-s}. \label{eq:int contour shift}
\end{align}
Further, since both the contours are fixed for a fixed $m$,
we have 
\begin{align}
\Bigg|	\sum_{\substack{k\geq 0 \\ c_\rho+U < \delta+k < d_m}}
	\Res_{s=\delta+k}	 	G_\rho(s)x^{\rho+\delta-s}\Bigg|
	\ll
	 \sum_{\substack{k\geq 0 \\ c_\rho+U < \delta+k < d_m}}
	 x^{\rho-k} 
	 \ll_m 
	  x^{\rho+\delta-c_\rho-U}.
\end{align}
Substituting the expansion of $G_\rho(s)$, we obtain
\begin{align}\label{eq:int I_rho}
	I_\rho(x)
	=\sum_{n=0}^{m}\frac{1}{2\pi i}\int_{C_R(c_\rho, d_m)}
	F_n^\rho(s)x^{\rho+\delta-s}\,ds
	+\frac{1}{2\pi i}\int_{C_R(c_\rho, d_m)}
	H_m(s)x^{\rho+\delta-s}\,ds
	+O_m\left(x^{\rho+\delta-c_\rho-U}\right).
\end{align}
By the choice of $d_m$ and $R$, all the singularities of $F_n^\rho(s)$ for $n= 0, \hdots, m,$ lie to the left of $C_R(c_\rho,d_m)$.
Then, we proceed as in (2.24)--(2.26) of \cite{rsv}
and obtain
\begin{align}\label{eq:int F_n}
	\frac{1}{2\pi i}\int_{C_R(c_\rho,d_m)}
	F_n^\rho(s)x^{\rho+\delta-s}\,ds
	=A_nx^{\frac{\omega-n-2i\nu-\frac12}{2A}}
	\cos\left((h^{-1}x)^{\frac1{2A}}
	+\pi\left(A\gamma+\frac n2+i\nu-\mu\right)\right).
\end{align}
In order to estimate the integral involving $H_m(s)$, we write 
\begin{align}
	\frac{x^{\rho+\delta}}{2\pi i}\int_{C_R(c_\rho, d_m)}
	H_m(s)x^{-s}\,ds 
	= \frac{x^{\rho+\delta}}{2\pi i}\int_{d_m-i\infty}^{d_m+i\infty}
	H_m(s)x^{-s}\,ds,
\end{align}
since the region between two contours contains no singularities of the integrand. Proceeding as in (2.27)--(2.29) of \cite{rsv}, we obtain
\begin{align}\label{eq:int H_m}
	\frac{1}{2\pi i}\int_{C_R(c_\rho,d_m)}
	H_m(s)x^{\rho+\delta-s}\,ds
	\ll_m x^{\frac{\omega-m-1}{2A}+\xi}.
\end{align}
Combining \eqref{eq:int I_rho}, \eqref{eq:int F_n}, and \eqref{eq:int H_m}, we obtain
\begin{align}
I_\rho(x)
&=
\sum_{n=0}^{m}	A_nx^{\frac{\omega-n-2i\nu-\frac12}{2A}}
\cos\left((h^{-1}x)^{\frac1{2A}}
+\pi\left(A\gamma+\frac n2+i\nu-\mu\right)\right) 
\\
&+O_m\left(x^{\rho+\delta-c_\rho-U}\right)
+ O_m\left(x^{\frac{\omega-m-1}{2A}+\xi}\right).
\end{align}

Finally, suppose that $c_\rho+ U \geq d_m$. In this case, all the poles of
$F_n^\rho(s)$, $0\leq n\leq m$, already lie to the left of
$C_R(c_\rho,c_\rho+U)$. Hence, we do not need to shift the contour as done in \eqref{eq:int contour shift}, instead we directly substitute the expansion of $G_\rho(s)$. The main terms are then evaluated as in \eqref{eq:int F_n}. To estimate the remainder, we shift the contour  of integration from $C_R(c_\rho, c_\rho+U)$  
to $C_R(d_m, c_\rho+U)$ without crossing any singularities and obtain the same bound as in \eqref{eq:int H_m}. This completes the proof.
\end{proof}

\normalcolor

We will need the following versions of results of Szeg{\"o} \cite{szego1926}, also given by Berndt \cite[Lemmas 3.2 and 3.4]{BerndtAvgI}. 
\begin{lemma}  \label{lem:szego's lem}
	Let $0 < \gamma_1 < 1, 0 < \gamma_2 < \frac{1}{3}$, and $\ell > \frac{3 \gamma_2}{2}$. Then, for
	$x^{-\gamma_1} \leq u \leq x^{\gamma_2},$
	we have as 
	$x \rightarrow \infty$,
	 \begin{align}
	 	\int_{0}^{\infty}
	 	(e^{1-t}t)^x
	 	t^{\frac{\omega}{2}- \frac{1}{4}} 
	 	\cos\big(c(xtu)^{\frac{1}{2}}+ \beta \pi\big) dt
	 	=
	 	\sqrt{2 \pi} e^{-c^2u/8} x^{-\frac{1}{2}} 
	 	\cos\big(c(xu)^{\frac{1}{2}}+ \beta \pi\big)
	 	+ O(x^{-1+\ell}),
	 \end{align}
  where $c$ is a positive constant and $\omega$ is as defined in \eqref{notations in lem for I_rho}.
\end{lemma}
\normalcolor
\begin{lemma}\label{lem:bdd for integral I_2}
	Let $\epsilon > 0$.
	Then for 
	$u> x^{\epsilon}$ 
	as 
	$x$ 
	continuously tends to 
	$\infty$, we have
	\begin{align}
		\frac{1}{\Gamma(x+1)} 
		\int_{0}^{\infty}
		e^{-t^2} 
		t^{2x}
		e^{i ut} dt
		= O(u^{-M}),
	\end{align}
where $M$ is any fixed positive number.
	 
\end{lemma}

The following lemma is due to Soundararajan \cite[Lemma 2.1]{Sound}. Let $(f(n))$ be a sequence of non-negative real numbers and $(\lambda_n)$ be a non-decreasing sequence of non-negative real numbers. Suppose that $\sum_{n=1}^{\infty}f(n) < \infty$. Consider the trigonometric series
\begin{align}\label{def of F(x)}
	F(x):= \sum_{n=1}^{\infty}f(n) \cos (2 \pi \lambda_n x + \beta),
\end{align} where $\beta \in \mathbb{R}$.
\begin{lemma} \label{lem:sound's key lemma}
	Let $L \geq 2$ and $N \geq 1$ be integers. Let $\mathcal{M}$ be a set of $M$ integers such that $\lambda_m \in [\lambda_N/2,3 \lambda_N/2]$  for each $m \in \mathcal{M}$. For any $X \geq 2$, there exists a point $x \in [X/2, (6L)^{M+1}X]$ such that 
	\begin{align}\label{eq:abs bdd in sound lem}
		|F(x)| 
		\geq
		\frac{1}{8} \sum_{m \in \mathcal{M}} f(m)
		-
		\frac{1}{L-1} \sum_{\substack{n \\ \lambda_n \leq 2 \lambda_N}} f(n)
		-
		\frac{4}{\pi^2 X \lambda_N} \sum_{n} f(n).
	\end{align}
If $\beta \equiv 0 \pmod {2\pi}$, then there is a point $x \in [X/2, (6L)^{M+1}X]$ such that 
\begin{align}\label{eq:bdd in sound lem}
     F(x)
	\geq
	\frac{1}{8} \sum_{m \in \mathcal{M}} f(m)
	-
	\frac{1}{L-1} \sum_{\substack{n \\ \lambda_n \leq 2 \lambda_N}} f(n)
	-
	\frac{2}{\pi^2 X \lambda_N} \sum_{n} f(n).
\end{align}
If $\beta \equiv \pi \pmod {2 \pi}$, then the conclusion  \eqref{eq:bdd in sound lem} holds with $-F(x)$ in place of $F(x)$.
\end{lemma}

 A refinement of the above lemma was recently obtained  by Mahatab \cite[Theorem 6]{mahatab2025} and Sourmelidis \cite[Theorems 1.1, 1.2]{sourmelidis2024}. While it is possible to incorporate either of these results into our proof to obtain Theorem \ref{thm:main theorem}, at this point Theorems 1.1 and 1.2  of \cite{sourmelidis2024} are more amenable to our setup. We state the relevant result below in terms of the function $F(x)$ of \eqref{def of F(x)}.  
\begin{lemma}\cite[Theorems 1.1 and 1.2]{sourmelidis2024} \label{lem:sourmelidis theorem}
	Let $N,M \geq 1$ be any integers.
	Let $\mathcal{M} \subseteq \mathbb N$ be a set of size $M$ such that $\lambda_m \in [\lambda_N/2,3 \lambda_N/2]$ 
	for each
	$m \in \mathcal{M}$.  
	For any
	$1 \leq Y < X$,
	there exist a point 
	 $x \in [Y/2, 2X]$ 
	 and an explicitly computed absolute constant $C>0$ \normalcolor
	 such that 
	\begin{align}\label{eq:abs bdd in sourmelidis lem}
		|F(x)| 
		\geq 
		C \sum_{m \in \mathcal{M}} f(m)
		+
		O\Bigg(\frac{2^M Y \log X}{X}  \sum\limits_{\substack{n \\ \lambda_n \leq 2 \lambda_N }} f(n)\Bigg)  \normalcolor
		-
		\frac{4}{\pi^2 Y \lambda_N} \sum_{n} f(n).
	\end{align}
	If $\beta \equiv 0 \pmod {2\pi}$, then there exist a point 
	 $x \in [Y, X]$ 
	 and an absolute constant 
	 $C>0$
	 as above,  \normalcolor
	  such that 
	\begin{align}\label{eq:bdd in sourmelidis lem}
		F(x)
		\geq 
		2C \sum_{m \in \mathcal{M}} f(m)
		+
		O\Bigg(\frac{2^M Y \log X}{X}  \sum_{n \geq 1}f(n)\Bigg).
	\end{align}
	If $\beta \equiv \pi \pmod {2 \pi}$, then the conclusion  \eqref{eq:bdd in sourmelidis lem} holds with $-F(x)$ in place of $F(x)$.
\end{lemma}

We will also need the following results.
\begin{lemma} \cite[Section II.6.1, Theorem 4]{tenenbaumbook} \label{lem: size having exactly k fact}
Let 
$\pi(x, k)$ 
denote the number of positive integers up to $x$ that have exactly $k$ distinct prime factors. Let
$A> 0$.
There exist positive constants 
$c_1 = c_1(A)$
 and 
 $c_2= c_2(A)$
 such that, uniformly for 
 $x \geq 3$,
 $1 \leq k \leq A \log_2 x$,
 $N \geq 0$,
 we have 
 \begin{align}
 	\pi(x, k)
 	=
 	\frac{x}{\log x}
 	\left\{ 
 	\sum_{j=0}^{N}
 	\frac{P_{j,k}(\log_2 x)}{(\log x)^j}
 	+ O \Big( \frac{(\log_2 x)^k}{k !} R_N(x) \Big)
 	\right\},
 \end{align}
where
$P_{j,k}(X)$
is a polynomial of degree at most $k-1$ and 
\begin{align}
	R_N(x)=
	e^{-c_1\sqrt{\log x}}
	+ 
	\left(\frac{c_2 N +1}{\log x}\right)^{N+1}.
\end{align}
	In particular, we have 
	\begin{align}
		P_{0,k}(X) = 
		\sum_{m + \ell = k-1} \frac{1}{m! \ell!} \lambda^{(m)}(0) X^{\ell},
	\end{align}
	where 
$
		\lambda(z)= 
	{\Gamma(z+1)}^{-1} 
		\prod_{p} \left( 1+ \frac{z}{p-1}\right) 
		\left(1- \frac{1}{p}\right)^z,
$
and $\lambda^{(m)}(z)$ denotes its $m$th derivative. 
	Moreover, under the same conditions, we have 
	\begin{align}
		\pi(x, k)
		=
		\frac{x}{\log x}
		\frac{(\log_2 x)^{k-1}}{(k-1)!}
		\left\{ \lambda\Big( \frac{k-1}{\log_2 x}\Big) + O \Big(\frac{k}{(\log_2 x)^2}\Big) \right\}.
	\end{align}
\end{lemma}

We record the following lemmas from \cite{BerndtAvgI} and \cite{Hafner}, following the notation of Section \ref{sec:notation and result}.

\begin{lemma}\cite[Lemma 2.2.1]{Hafner}\label{lem:asymp of b_n with log}
	Suppose \ref{H3} holds for the partial sums of coefficients of the  Dirichlet series 
	$\psi(s)=\sum_{n=1}^{\infty} b_n \mu_n^{-s}$. Then for 
	$\tau \leq \sigma_b$, 
	we have
	\begin{align}
		\sum_{\mu_n \leq x} b_n \mu_n^{-\tau} 
		\sim
		\begin{cases}
			\begin{array}{ll}
				 c x^{\sigma_b - \tau} (\log x)^{\ell-1},	
				& \text{if}  \quad \tau < \sigma_b, \label{sum b_n first}
				\\
				c (\log x)^{\ell},		
				& 
				\text{if} \quad \tau = \sigma_b, \label{sum b_n second}
			\end{array}
		\end{cases}
	\end{align}
	as
$x \rightarrow \infty$. In particular, we note that 
\begin{align}
	\sum_{\mu_n \leq x} b_n \mu_n^{-\tau} 
	=
	O\left(x^{\sigma_b - \tau + \epsilon}\right), \label{eq:asymp of b_n}
\end{align}
as $x \rightarrow \infty$.
\end{lemma}

\begin{lemma} \cite[Lemma 3.9] {BerndtAvgI}\label{lem:tail of sum b_n}
	Suppose \ref{H3} holds for the partial sums of coefficients of the  Dirichlet series 
$\psi(s)=\sum_{n=1}^{\infty} b_n \mu_n^{-s}$. Then
	for $\tau > \sigma_b$, 
	as
	$x \rightarrow \infty$,
	we have
	\begin{align}
		\sum_{\mu_n > x}
		b_n 
		\mu_n^{-\tau} 
		=
		O(x^{\sigma_b-\tau} (\log x)^{\ell-1}).
	\end{align} 
\end{lemma}

\begin{lemma}\cite[Eq. 2.2.4]{Hafner} \label{lem: bdd for error E_1}
		Suppose \ref{H3} holds for the partial sums of coefficients of the  Dirichlet series 
	$\psi(s)=\sum_{n=1}^{\infty} b_n \mu_n^{-s}$. 
	Let
	$\r \leq  (2\s_b-\d)A-\frac{1}{2}$. 
	For $u,x > 0$, let 
	\begin{align}
		F(u,x)
		=
		\sum_{\mu_n >  x} 
		b_n 
		\mu_n^{-\frac{\delta}{2}- \frac{\rho}{2A}- \frac{1}{4A}  }
		e^{-\frac{1}{2}(\mu_n u)^{\frac{1}{A}}}.
	\end{align}
Then for any
$\epsilon >0$, 
as 
$u \rightarrow 0$, we have
\begin{align}
	F(u,x) 
	\ll 
	u^{-\sigma_b+\frac{\delta}{2}+ \frac{\rho}{2A}+ \frac{1}{4A}  - \epsilon}
	e^{-\frac{1}{4}(u x)^{\frac{1}{A}}},
\end{align}
uniformly for 
$x > \frac{\mu_1}{2}$.

\end{lemma}

The following lemma gives an explicit formula for $r_Q(n)$, the number of representations of a positive integer $n$ by quadratic form $\mathcal{Q}=u^2+5v^2$.
\begin{lemma}\cite[Corollary 3.3]{BerkovichHamza2009}\label{lem:count r_Q(n)}
		 If the prime factorization of $n$ is given by 
	$$
	n=2^a5^b\prod_{i=1}^r p_i^{v_i}\prod_{j=1}^s q_j^{w_j},
	$$
	where
	$
	p_i\equiv1,3,7,9\pmod{20}$,
	and
	$q_j\equiv11,13,17,19\pmod{20},
	$
	then
	$$
	r_\mathcal Q(n)=\bigl(1+(-1)^{a+t}\bigr)
	\prod_{i=1}^r(1+v_i)
	\prod_{j=1}^s\frac{1+(-1)^{w_j}}{2},
	$$
	where $t$ is the number of prime factors of $n$, counting multiplicity, that are congruent to $3$ or $7$ modulo $20$.
\end{lemma}

The following lemma is a corollary to Theorem 2 of Girstmair, K\"uhleitner, M\"uller and Nowak \cite[p. 189]{Girstmair}.
\begin{lemma}\label{lem:girstmair size}
	Let $q\geq 1$, and let $a_1,\ldots,a_d$ be $d$ integers which are
	incongruent modulo $q$ and relatively prime to $q$. Denote by
	$\pi(x,r_1,\ldots,r_d)$ the number of integers $n\leq x$ which have
	$r_v$ distinct prime divisors $p\equiv a_v \pmod q$ for $1\leq v\leq d$,
	and no other ones. For given $B\geq 1$, and
	$
	r=\sum_{v=1}^d r_v,
	$
	we have, uniformly in $1\leq r\leq B\log_2 x$ and $x\geq 3$,
	\begin{align}
		\pi(x,r_1,\ldots,r_d)
		&=
		\frac{r}{r_1!\cdots r_d!\,\phi(q)^r}
		\frac{x}{\log x}
		(\log_2 x)^{r-1}
		\left(
		1+O\left(\frac{r}{\log_2 x}\right)
		\right),
		\\
		\pi(x,r_1,\ldots,r_d)
		&\asymp
		\frac{r}{r_1!\cdots r_d!\,\phi(q)^r}
		\frac{x}{\log x}
		(\log_2 x)^{r-1}.
	\end{align}
\end{lemma}
\normalcolor

\section{Proof of Theorem \ref{thm:main theorem}} \label{sec: results}
\subsection{A key proposition}
We begin with the following proposition which works as an essential tool for the proof of the theorem. 

\begin{proposition}\label{prop: key proposition}
	Consider 
	$\phi(s), \psi(s) \in \mathcal{C}$ 
as defined in \eqref{dirichlet series phi psi} satisfying  hypotheses \ref{H1} to \ref{H3}.
	For $x>0$ and $\rho \geq 0$, consider $A^{\rho}(x),Q_{\r}(x)$ as defined in \eqref{def:A_rho} and \eqref{def:Q_rho} respectively. 	Then, for 
	$\r < (2\s_b-\d)A-\frac{1}{2}$, we have
	\begin{align}
	&\frac{1}{\Gamma(x+1)}
	\int_{0}^{\infty}
	\left[
	A^{\r}(uv^A) - Q_{\r}(uv^A) 
	\right] e^{-v} v^x  dv
	\\
	&= 
	(1+o(1))
	A_0 	\omega' Q^{ \frac{1}{2A} + \frac{\rho}{A}}
	u^{\frac{\delta}{2} - \frac{\rho}{2A}- \frac{1}{4A} + \rho}
	x^{\frac{\delta A}{2} - \frac{\rho}{2}- \frac{1}{4} + A \rho} \times
	\\
	&
	\sum_{n=1}^{\infty} 
	b_n 
	\mu_n^{-\frac{\delta}{2}- \frac{\rho}{2A}- \frac{1}{4A} }
	e^{-D(\mu_n u)^{\frac{1}{A}}}
	\cos\left(\frac{(\mu_n u)^{\frac{1}{2A}} \sqrt{x}}{(hQ^2)^{\frac{1}{2A}}}+ \pi (A \gamma - \mu) \right)
	+ 
	O \left(x^{\frac{\delta A}{2} - \frac{\rho}{2}- \frac{1}{2} + A \rho + \epsilon}\right), \label{final estimation}
	\end{align}
for any sufficiently small $\epsilon > 0$ and $D=  \frac{1}{8}(hQ^2)^{-\frac{1}{A}},$ 
provided 
\begin{align}\label{condition on u in prop}
	x^{-\epsilon_0 A} \leq u \leq u_0,
\end{align}
for any
$\epsilon_0 >0$
and 
 $u_0 > 0$  sufficiently small.
 The parameters 
	$A_0, h, \mu, \gamma,$  
are	as defined in Lemma \ref{lem:asym for I_rho}.
\end{proposition}

\begin{proof}
Using Theorem \ref{thm:FE implies Rieszsum}, for 
$\rho \geq 0$ 
and 
$m \in \mathbb{N}$
sufficiently large, we write
\begin{align} \label{riesz sum 1}
	A^{\r+m}(v) - Q_{\r+m}(v)
	= 
	\omega' Q^{\delta+2 \rho +2m}
	\sum_{n=1}^{\infty} \frac{b_n}{\mu_n^{\d + \r +m}} 
	  I_{\r+m}^{(a)}\Big(\frac{\mu_n v}{Q^2}\Big),
\end{align}
where 
$I_{\r+m}^{(a)}(\ast)$ is defined in \eqref{def:I_rho}.
Note that in $I_{\r+m}^{(a)}(\ast)$, we can choose $a$ to be large enough which will enable us to make the parameter $U$ large in the  $I_{\r+m}(\ast)$ whenever needed.
\normalcolor
Let 
$x$ 
be a sufficiently large positive number. We replace 
$v$ 
by 
$uv^{A}$ in \eqref{riesz sum 1}, 
where $u$ is a parameter that will be chosen later and multiply throughout by 
$$
\frac{D^m}{(Dv)^m} (e^{-v} v^x),
$$
where the operator 
$\frac{D}{Dv}$
is given by
$ \frac{d}{dv} \frac{1}{v^{A-1}}$.
Integrating with respect to 
$v$ 
from 
$0$ to $\infty$,
using integration by parts,
the left-hand side of \eqref{riesz sum 1} becomes
\begin{align}
	\int_{0}^{\infty}
	&\left[
	A^{\r+m}(uv^A) -  Q_{\r+m}(uv^A)
	\right]
	\frac{D^m}{(Dv)^m}
	 (e^{-v} v^x)
     dv 
     \\
    & = 
     \left[
     \left[
     A^{\r+m}(uv^A) - Q_{\r+m}(uv^A)
     \right] 
     \frac{1}{v^{A-1}}
     \frac{D^{m-1}}{(Dv)^{m-1}}
     (e^{-v} v^x)
     \right]_0^{\infty}
     \\
     &-
     \int_{0}^{\infty} uAv^{A-1} 
     \left[
     A^{\r+m-1}(uv^A) - Q_{\r+m-1}(uv^A)
     \right]
     \frac{1}{v^{A-1}}
    \frac{D^{m-1}}{(Dv)^{m-1}}
    (e^{-v} v^x) 
     dv.
\end{align}
 The first term on the right-hand side vanishes due to the exponential factor and the second term can be evaluated by doing integration by parts again. After repeated integration by parts $m-1$ times, the left-hand side above is equal to
\begin{align}\label{lhs of borel mean}
	(-uA)^m 
	\int_{0}^{\infty}
	\left[
	A^{\r}(uv^A) - Q_{\r}(uv^A) 
	\right] e^{-v} v^x  dv.
\end{align} The right-hand side of \eqref{riesz sum 1} becomes 
\begin{align} \label{rhs of borel mean 1}
	\int_{0}^{\infty}  
	\omega' Q^{\delta+2 \rho +2m}
	\sum_{n=1}^{\infty} \frac{b_n}{\mu_n^{\d + \r +m}} 
	&I_{\r+m}^{(a)}\Big(\frac{\mu_n uv^A}{Q^2}\Big)
	\frac{D^m}{(Dv)^m}
	(e^{-v} v^x)
	 dv
	\\
	&= 
	\omega' Q^{\delta+2 \rho +2m}
	\sum_{n=1}^{\infty} \frac{b_n}{\mu_n^{\d + \r +m}} 
	\int_{0}^{\infty} 
	I_{\r+m}\Big(\frac{\mu_n uv^A}{Q^2}\Big)
	\frac{D^m}{(Dv)^m}
	(e^{-v} v^x)  
	dv,
\end{align} upon interchanging the order of integration and summation due to absolute convergence of the integral and  using Lemma \ref{lem:I_rho to I*_rho}.
Using integration by parts and Lemma \ref{lem:d/dx I_rho} we get,
 \begin{align}
 	\int_{0}^{\infty} 
 	I_{\r+m}\Big(\frac{\mu_n uv^A}{Q^2}\Big)
 	\frac{D^m}{(Dv)^m}
 	(e^{-v} v^x)  dv
 	&=
    \left[
    I_{\r+m}\Big(\frac{\mu_n uv^A}{Q^2}\Big)  \frac{1}{v^{A-1}}
    \frac{D^{m-1}}{(Dv)^{m-1}}
    (e^{-v} v^x) 
    \right]_0^\infty
    \\
    &- 
    \frac{\mu_n uA}{Q^2}
    \int_{0}^{\infty} 
    I_{\r+m-1}\Big(\frac{\mu_n uv^A}{Q^2}\Big) 
    \frac{D^{m-1}}{(Dv)^{m-1}}
    (e^{-v} v^x)  dv.
\end{align}
The first term on the right-hand side above vanishes due to the exponential factor and we do integration by parts to evaluate the second term. By repeating integration by parts another $m-1$ times, we obtain 
\begin{align} \label{rhs of borel mean 2}
	\int_{0}^{\infty} 
	I_{\r+m}\Big(\frac{\mu_n uv^A}{Q^2}\Big)
	\frac{D^m}{(Dv)^m}
	(e^{-v} v^x)  
	 dv
	=
	\Big(\frac{-\mu_n u A}{Q^2}\Big)^m
	\int_{0}^{\infty} 
	I_{\r}\Big(\frac{\mu_n uv^A}{Q^2}\Big) e^{-v} v^x  dv.
\end{align} Combining \eqref{rhs of borel mean 1} and \eqref{rhs of borel mean 2}, we get
\begin{align} \label{rhs of borel mean}
	\int_{0}^{\infty}  
	\omega' Q^{\delta+2 \rho +2m}
	\sum_{n=1}^{\infty} \frac{b_n}{\mu_n^{\d + \r +m}} 
	I_{\r+m}^{(a)}\Big(\frac{\mu_n uv^A}{Q^2}\Big)
	&\frac{D^m}{(Dv)^m}
	(e^{-v} v^x)   dv
	\\
	&=
	(- u A)^m
	\omega' Q^{\delta+2 \rho}
	\sum_{n=1}^{\infty} \frac{b_n}{\mu_n^{\d + \r }} 
	\int_{0}^{\infty} 
	I_{\r}\Big(\frac{\mu_n uv^A}{Q^2}\Big) e^{-v} v^x  dv.
\end{align}
From \eqref{lhs of borel mean} and \eqref{rhs of borel mean}, we obtain
\begin{align}\label{lhsrhs borel mean}
	\int_{0}^{\infty}
	\left[
	A^{\r}(uv^A) - Q_{\r}(uv^A) 
	\right] e^{-v} v^x  dv
	=
	\omega' Q^{\delta+2 \rho}
	\sum_{n=1}^{\infty} \frac{b_n}{\mu_n^{\d + \r }} 
	\int_{0}^{\infty} 
	I_{\r}\Big(\frac{\mu_n uv^A}{Q^2}\Big) e^{-v} v^x  dv.
\end{align}
Letting 
$\gamma_1 \in (0,1)$,
$\gamma_2 \in \left(0, \frac{1}{3}\right)$,
we split the sum over all 
$n$ 
into the parts
$\mu_n \in \left(0, \frac{1}{u} x^{-A \gamma_1}\right)$,
$\mu_n \in \left[\frac{1}{u} x^{-A \gamma_1}, \frac{1}{u} x^{A \gamma_2}\right]$ 
and 
$\mu_n > \frac{1}{u} x^{A \gamma_2} $, 
where $u$ and $x$ are parameters such that $x$ is sufficiently large and $u$ satisfies
\begin{align}\label{condition on u}
	x^{-\gamma_1 \epsilon} \leq u \leq u_0
\end{align}
for some 
$0 < \epsilon < A$
and 
$u_0$ 
arbitrarily small. For $x$ sufficiently large, the first interval
$\left(0, \frac{1}{u} x^{-A \gamma_1}\right)$ 
is empty since
$\frac{1}{u}x^{- A\gamma_1} \leq x^{-\gamma_1(A - \epsilon)} \rightarrow 0$ 
as 
$x \rightarrow \infty$.
We may thus write the sum in \eqref{lhsrhs borel mean} as

\begin{align}\label{partion into integrals}
	\bigg[
	 \sum_{\frac{1}{u} x^{- A\gamma_1} \leq \mu_n \leq    \frac{1}{u}x^{A\gamma_2}}
	+ \sum_{ \mu_n > \frac{1}{u}x^{A\gamma_2}}
	\bigg]	
	\frac{b_n}{\mu_n^{\d + \r }} 
	\int_{0}^{\infty} 
	I_{\r}\Big(\frac{\mu_n uv^A}{Q^2}\Big) e^{-v} v^x  dv
	= I_1+I_2 \quad \text{(say)}.
\end{align}
We use the asymptotic formula for 
$I_{\r}$ 
from Lemma \ref{lem:asym for I_rho} 
with 
$m=0$,
$\nu = 0$ (due to \ref{H2})
and 
$\omega, h, \mu, \gamma$ 
as defined in \eqref{notations in lem for I_rho}.
 In particular, we choose $U$ such that 
$c_\rho+U \geq d_0$ in Lemma \ref{lem:asym for I_rho} and use the asymptotic formula \eqref{formula for I_rho} without the second error term.
\normalcolor 
This gives
\begin{align}
	I_1
	= 
	&\sum_{0 < \mu_n \leq  \frac{1}{u}x^{A\gamma_2}} 
	\frac{b_n}{\mu_n^{\d + \r }} \times
	\\
	&\int_{0}^{\infty}
	\left\{
	A_0
	\Big(\frac{\mu_n uv^A}{Q^2}\Big)^{\frac{\omega  -1/2}{2A}} 
	\cos\left(\Big(\frac{\mu_n u v^A}{hQ^2}\Big)^{\frac{1}{2A}}+\pi (A \gamma - \mu) \right)
	+ O\left( \Big(\frac{\mu_n uv^A}{Q^2}\Big)^{\frac{\omega -1}{2A}+\xi} \right)
	\right\}
	e^{-v} v^x  dv
\end{align}
\begin{align}
	=
	\sum_{0 < \mu_n \leq  \frac{1}{u}x^{A\gamma_2}} 
	\frac{b_n}{\mu_n^{\d + \r }} 
	\bigg\{
	&A_0
	\Big(\frac{\mu_n u}{Q^2}\Big)^{\frac{\omega  -1/2}{2A}} 
	\int_{0}^{\infty}
	e^{-v} 
	v^{x+ \frac{\omega}{2} -\frac{1}{4}}
	\cos\left(\frac{(\mu_n u)^{\frac{1}{2A}} v^{\frac{1}{2}}}{(hQ^2)^{\frac{1}{2A}}}+\pi (A \gamma - \mu) \right) dv 
	\\
	& + O \left( 
	\left(\mu_nu \right)^{\frac{\omega -1}{2A}+\xi}
	\int_{0}^{\infty}
	v^{x+ \frac{\omega}{2}- \frac{1}{2}+\xi A}
	e^{-v} dv
	\right)
	\bigg\}.
\end{align}
We make a change of variable from $v$ to $xv$ in the first integral above and use the definition of the Gamma function for the second integral to obtain
\begin{align}
	I_1
	= \!\!\!\!\!\!
	\sum_{0 < \mu_n \leq  \frac{1}{u}x^{A\gamma_2}} 
	\frac{b_n}{\mu_n^{\d + \r }} 
	\!
	\bigg\{
	&A_0
	\Big(\frac{\mu_n u}{Q^2}\Big)^{\frac{\omega -1/2}{2A}} 
	x^{x+ \frac{\omega}{2} + \frac{3}{4}} e^{-x}
	\int_{0}^{\infty}
	e^{x-xv} 
	v^{x+ \frac{\omega}{2} - \frac{1}{4}}
	\cos\left(\frac{(\mu_n u)^{\frac{1}{2A}} (xv)^{\frac{1}{2}}}{(hQ^2)^{\frac{1}{2A}}}+\pi (A \gamma - \mu) \right) dv 
	\\
	&+ O \left( 
	\left(\mu_nu \right)^{\frac{\omega - 1}{2A}+\xi}
	\Gamma \left(
	x+ \frac{\omega}{2} +\frac{1}{2}+\xi A 
	\right)
	\right)
	\bigg\}.
\end{align}
For any 
$\ell > \frac{3 \gamma_2}{2}$, 
Lemma \ref{lem:szego's lem} with 
 $u=(\mu_n u)^{\frac{1}{A}}$ 
 and 
 $c= (hQ^2)^{-\frac{1}{2A}}$ 
 gives
\begin{align}
	I_1
	=
	\sum_{0 < \mu_n \leq  \frac{1}{u}x^{A\gamma_2}} 
	\frac{b_n}{\mu_n^{\d + \r }} 
	&\bigg\{
	A_0 \sqrt{2 \pi}
	\Big(\frac{\mu_n u}{Q^2}\Big)^{\frac{\omega -1/2}{2A}} 
	x^{x+ \frac{\omega}{2} + \frac{1}{4}} e^{-x}
	e^{-D(\mu_n u)^{\frac{1}{A}}}
	\cos\left(\frac{(\mu_n u)^{\frac{1}{2A}} x^{\frac{1}{2}}}{(hQ^2)^{\frac{1}{2A}}}+ \pi (A \gamma - \mu)\right)
	\\
	&+O\big(
	\left(\mu_nu\right)^{\frac{\omega -1/2}{2A}} 
	x^{x+ \frac{\omega}{2} + \frac{3}{4}} e^{-x}x^{-1+\ell}
	\big)
	+ O \left( 
	\left(\mu_nu \right)^{\frac{\omega -1}{2A}+\xi}
	\Gamma \Big(x  + \frac{\omega}{2} + \frac{1}{2}+\xi A \Big)
	\right)
	\bigg\}, \label{eq:I_1 after szego lem}
\end{align}
where $D:= \frac{1}{8}(hQ^2)^{-\frac{1}{A}}$.
\normalcolor
We write the main term of 
$I_1$ 
as the sum over all 
$\mu_n$ 
and subtract the tail corresponding to 
$\mu_n >  \frac{1}{u}x^{A\gamma_2}$.
This gives that the main term of 
$I_1$ 
equals
\begin{align}
	&A_0 \sqrt{2 \pi}
	\Big(\frac{u}{Q^2}\Big)^{\frac{\omega -1/2}{2A}}
	x^{x+ \frac{\omega}{2} + \frac{1}{4}} e^{-x}
	\sum_{\mu_n } 
	b_n 
	\mu_n^{\frac{\omega -1/2}{2A}-\d - \r }
	e^{-D(\mu_n u)^{\frac{1}{A}}}
	\cos\left(\frac{(\mu_n u)^{\frac{1}{2A}} x^{\frac{1}{2}}}{(hQ^2)^{\frac{1}{2A}}}+ \pi (A \gamma - \mu) \right)
	\\
	& +
	O \bigg(
	u^{\frac{\omega -1/2}{2A}}
	x^{x+ \frac{\omega}{2} + \frac{1}{4}} e^{-x}
	\sum_{\mu_n >  \frac{1}{u}x^{A\gamma_2}} 
	b_n 
	\mu_n^{\frac{\omega -1/2}{2A}-\d - \r }
	e^{-D(\mu_n u)^{\frac{1}{A}}}
	\bigg).
\end{align} \normalcolor
Substituting this into \eqref{eq:I_1 after szego lem}, we have obtained 
\begin{align} 
	I_1
	&=
	A_0 \sqrt{2 \pi}
	\Big(\frac{u}{Q^2}\Big)^{\frac{\omega -1/2}{2A}}
	x^{x+ \frac{\omega}{2} + \frac{1}{4}} e^{-x}
	\sum_{\mu_n } 
	b_n 
	\mu_n^{\frac{\omega -1/2}{2A}-\d - \r }
	e^{-D(\mu_n u)^{\frac{1}{A}}}
	\cos\left(\frac{(\mu_n u)^{\frac{1}{2A}} x^{\frac{1}{2}}}{(hQ^2)^{\frac{1}{2A}}}+ \pi (A \gamma - \mu) \right)
	\\
	&+
	O \bigg(
	u^{\frac{\omega -1/2}{2A}}
	x^{x+ \frac{\omega}{2} + \frac{1}{4}} e^{-x}
	\sum_{\mu_n >  \frac{1}{u}x^{A\gamma_2}} 
	b_n 
	\mu_n^{\frac{\omega -1/2}{2A}-\d - \r }
	e^{-D(\mu_n u)^{\frac{1}{A}}} \normalcolor
	\bigg) 
	\\
	&+ O \bigg(
	u^{\frac{\omega -1/2}{2A}}
	x^{x+ \frac{\omega}{2} - \frac{1}{4}+\ell} e^{-x}
	\sum_{\mu_n \leq  \frac{1}{u}x^{A\gamma_2}} 
	b_n 
	\mu_n^{\frac{\omega -1/2}{2A}-\d - \r }
	\bigg)
	\\
	& + O \bigg(
	u^{\frac{\omega -1}{2A}+\xi}
	\Gamma \Big(x+ \frac{\omega}{2} + \frac{1}{2} +\xi A \Big)
	\sum_{\mu_n \leq  \frac{1}{u}x^{A\gamma_2}} 
	b_n 
	\mu_n^{\frac{\omega -1}{2A}-\d - \r +\xi }
	\bigg).
	\end{align}
Dividing
$I_1$ 
by
$\Gamma(x+1)$, 
 using the Stirling formula \eqref{stirling for real} for the first three terms, and using \eqref{gamma of x plus alpha} for the last term above,  we get
\begin{align}
	&\frac{1}{\Gamma(x+1)}  I_1
	\\
	&=
	(1+o(1))
	A_0
	\Big(\frac{u}{Q^2}\Big)^{\frac{\omega -1/2}{2A}}
	x^{ \frac{\omega}{2} - \frac{1}{4}} 
	\sum_{\mu_n } 
	b_n 
	\mu_n^{\frac{\omega -1/2}{2A}-\d - \r }
	e^{-D(\mu_n u)^{\frac{1}{A}}}
	\cos\left(\frac{(\mu_n u)^{\frac{1}{2A}} x^{\frac{1}{2}}}{(hQ^2)^{\frac{1}{2A}}}+ \pi (A \gamma - \mu) \right)
	\\
	&  +
	O \bigg(
	u^{\frac{\omega -1/2}{2A}}
	x^{ \frac{\omega}{2} - \frac{1}{4}} 
	\sum_{\mu_n >  \frac{1}{u}x^{A\gamma_2}} 
	b_n 
	\mu_n^{\frac{\omega -1/2}{2A}-\d - \r }
	e^{-D(\mu_n u)^{\frac{1}{A}}}
	\bigg) 
	\\
	&+ O \bigg(
	u^{\frac{\omega -1/2}{2A}}
	x^{ \frac{\omega}{2} - \frac{3}{4}+\ell}
	\sum_{\mu_n \leq  \frac{1}{u}x^{A\gamma_2}} 
	b_n 
	\mu_n^{\frac{\omega -1/2}{2A}-\d - \r }
	\bigg)
	\\
	&+ O \bigg(
	u^{\frac{\omega -1}{2A}+\xi}
	x^{ \frac{\omega}{2} - \frac{1}{2} + \xi A}
	\sum_{\mu_n \leq  \frac{1}{u}x^{A\gamma_2}} 
	b_n 
	\mu_n^{\frac{\omega -1}{2A}-\d - \r + \xi }
	\bigg)
\end{align}
\begin{align}
	&=
	(1+o(1))
	A_0
	\Big(\frac{u}{Q^2}\Big)^{\frac{\omega -1/2}{2A}}
	x^{ \frac{\omega}{2} - \frac{1}{4}} 
	\sum_{\mu_n } 
	b_n 
	\mu_n^{\frac{\omega -1/2}{2A}-\d - \r }
	e^{-D(\mu_n u)^{\frac{1}{A}}}
	\cos\left(\frac{(\mu_n u)^{\frac{1}{2A}} x^{\frac{1}{2}}}{(hQ^2)^{\frac{1}{2A}}}+ \pi (A \gamma - \mu) \right)
	\\
	& + E_1 + E_2 +E_3. \label{error in I_1}	\quad \text{(say)}
\end{align}
Using Lemma \ref{lem: bdd for error E_1} with 
$u$ replaced by 
$(2D)^A u$, we see that as 
$u \rightarrow 0$,
\begin{align}	
	E_1
	\ll_{\epsilon'}
	u^{\frac{\omega -1/2}{2A}}
	x^{ \frac{\omega}{2} - \frac{1}{4}} 	
	u^{-\sigma_b- \frac{\omega -1/2}{2A}+ \d + \r - \epsilon'}
	e^{-\frac{D}{2}x^{\gamma_2}},
\end{align} \normalcolor
 for any 
 $\epsilon' >0$.
Since 
$u$ 
is sufficiently small, we see that 
\begin{align}
	\label{u theta}
u^{-\theta} \le u^{-|\theta|},
\end{align}
for any $\theta \in \mathbb R$. Using this and 
$u^{-1} \leq x^{\gamma_1 \epsilon}$ from \eqref{condition on u}, 
we obtain the estimate
\begin{align}\label{error E1}
	E_1 
	&\ll_{\epsilon'}
	x^{\gamma_1 \epsilon \left|\frac{\omega -1/2}{2A} \right|}
	x^{ \frac{\omega}{2} - \frac{1}{4}} 
	x^{\gamma_1 \epsilon \left| -\sigma_b- \frac{\omega -1/2}{2A}+ \d + \r - \epsilon' \right|}
	e^{-\frac{D}{2}x^{\gamma_2}}
	\\
	&\ll x^{-M},
\end{align} \normalcolor
for any $M>0$, as $x \rightarrow \infty$.
To estimate 
$E_2$ 
in \eqref{error in I_1}, we apply \eqref{eq:asymp of b_n} 
to the partial sum inside the error term with 
$\tau = -\frac{\omega -1/2}{2A}+\d + \r$
and 
$x=\frac{1}{u}x^{A\gamma_2}$, 
noting that the condition   
\begin{align}\label{condition for sum b_n}
	-\frac{\omega -1/2}{2A}+\d + \r \leq \sigma_b
\end{align}
follows from
\begin{align}
	\r \leq  (2\s_b-\d)A-\frac{1}{2}.
\end{align}
This gives
\begin{align} \label{E2}
	E_2
	\ll_{\epsilon'}
	u^{\frac{\omega -1/2}{2A}}
	x^{ \frac{\omega}{2} - \frac{3}{4}+\ell}
	\left( \frac{1}{u}x^{A\gamma_2}\right)^{\sigma_b + \frac{\omega -1/2}{2A}-\d - \r + \epsilon'},
\end{align}
 for any 
$\epsilon' >0$.
The conditions 
$\frac{1}{u} \leq x^{\gamma_1 \epsilon}$
and 
$0 < \epsilon < A$
give
\begin{align}\label{relation u and x}
	\frac{1}{u}x^{A\gamma_2} \leq x^{(\gamma_1 + \gamma_2)A}.
\end{align}
Substituting \eqref{relation u and x} into \eqref{E2}, using the bounds \eqref{u theta}, \eqref{condition on u},  and putting 
$\ell = \frac{3\gamma_2}{2}+\frac{1}{16}$, 
we obtain
\begin{align}
	E_2 
	\ll_{\epsilon'}
	x^{ \frac{\omega}{2} - \frac{1}{2}}
	x^{\gamma_1 \epsilon \left|  \frac{\omega -1/2}{2A}  \right|+(\gamma_1 + \gamma_2) A \left(\kappa + \epsilon'\right)+ \frac{3\gamma_2}{2} - \frac{3}{16} },
\end{align}
where $\kappa= \sigma_b + \frac{\omega -1/2}{2A}-\d - \r \geq 0$ from \eqref{condition for sum b_n}.
Fixing
$\epsilon' >0$ 
and taking 
$\gamma_1, \gamma_2$ 
to be sufficiently small positive numbers, we see that as $x \rightarrow \infty$
\begin{align}\label{error E2}
	E_2 
	\ll
	x^{ \frac{\omega}{2} - \frac{1}{2}}.
\end{align} 
For 
$E_3$, 
if 
$-\frac{\omega -1}{2A}+\d + \r -\xi > \sigma_b$,
 then the series inside the error term is convergent. Therefore, we have
 \begin{align}
 	E_3 
 	\ll
 	u^{\frac{\omega -1}{2A}+\xi}
 	x^{ \frac{\omega}{2} - \frac{1}{2}+ \xi A }.
 \end{align}
 Using \eqref{u theta} and \eqref{condition on u}, the above bound implies
\begin{align}
		E_3 
	&\ll
	x^{\gamma_1 \epsilon\left|\frac{\omega -1}{2A}+ \xi \right|+ \frac{\omega}{2} - \frac{1}{2} + \xi A}
	\\
	&\ll
	x^{ \frac{\omega}{2} - \frac{1}{2} + \epsilon}, \label{error E11}
\end{align} 
for any 
$\epsilon > 0$ upon taking 
$\gamma_1$ and $\xi$ 
to be sufficiently small, as 
$x \rightarrow \infty$. 
On the other hand, whenever 
$-\frac{\omega -1}{2A}+\d + \r - \xi \leq  \sigma_b$,
we apply \eqref{eq:asymp of b_n} 
to the partial sum inside the error term with 
$\tau = -\frac{\omega -1}{2A}+\d + \r - \xi$ 
and 
$x=\frac{1}{u}x^{A\gamma_2}$
to get
\begin{align}
	E_3
	\ll_{\epsilon'}
	u^{\frac{\omega -1}{2A}+\xi}
	x^{ \frac{\omega}{2} - \frac{1}{2}+\xi A}
	\Big(\frac{1}{u}x^{A\gamma_2}\Big)^{ \left(\sigma_b + \frac{\omega -1}{2A}-\d - \r + \xi + \epsilon'\right) },
\end{align}
for any 
$\epsilon' > 0$. 
 Using \eqref{relation u and x}, along with the bounds \eqref{u theta} and \eqref{condition on u},  we obtain
\begin{align}
		E_3
		\ll_{\epsilon'}
	 	x^{ \frac{\omega}{2} - \frac{1}{2}+ \xi A}
	 	x^{\gamma_1 \epsilon \left| \frac{\omega -1}{2A}+ \xi  \right|+(\gamma_1 + \gamma_2)A \left(\sigma_b + \frac{\omega - 1}{2A}-\d - \r + \xi + \epsilon'\right)}.
\end{align}
Fixing 
$\epsilon'>0$
and taking 
$\xi,$
$\gamma_1,$
and
$\gamma_2$ 
to be sufficiently small positive numbers, we get 
\begin{align}\label{error E12}
E_3
\ll
x^{ \frac{\omega}{2} - \frac{1}{2}+ \epsilon},
\end{align} 
for any
$\epsilon > 0$
 as $x \rightarrow \infty$. 
We substitute the estimates
\eqref{error E1}, \eqref{error E2}, \eqref{error E11} and \eqref{error E12}
for
$E_1,E_2$ 
and 
$E_3$
 into \eqref{error in I_1} to obtain
\begin{align}
	&\frac{1}{\Gamma(x+1)}  I_1
	\\
	&=
	(1+o(1))
	A_0 
	\Big(\frac{u}{Q^2}\Big)^{\frac{\omega -1/2}{2A}}
	x^{ \frac{\omega}{2} - \frac{1}{4}} 
	\sum_{\mu_n } 
	b_n 
	\mu_n^{\frac{\omega -1/2}{2A}-\d - \r }
	e^{-D(\mu_n u)^{\frac{1}{A}}}
	\cos\left(\frac{(\mu_n u)^{\frac{1}{2A}} x^{\frac{1}{2}}}{(hQ^2)^{\frac{1}{2A}}}+ \pi (A \gamma - \mu) \right)
	\\
	&+ 
	O \left(x^{ \frac{\omega}{2} - \frac{1}{2}+ \epsilon}\right). \label{final estimation I_1}
\end{align}
We now proceed to evaluate the term $I_2$ in \eqref{partion into integrals}. We
divide
$I_2$ 
by
$\Gamma(x+1)$
and without loss of generality use the asymptotic formula  \eqref{formula for I_rho} for 
$I_{\r}(x)$ 
from Lemma \ref{lem:asym for I_rho} with
$\nu = 0$ (due to \ref{H2}) and 
$\omega, h, \mu$ 
as defined in \eqref{notations in lem for I_rho},
to get
\begin{align}
	\frac{1}{\Gamma(x+1)}	I_2
     =
	&\sum_{ \mu_n > \frac{1}{u}x^{A\gamma_2}}	
	\frac{b_n}{\mu_n^{\d + \r }} 
		\frac{1}{\Gamma(x+1)}
	\times
	\\
	&
	\int_{0}^{\infty} 
	\Bigg\{
	\sum_{t=0}^{m} 
	A_t \Big(\frac{\mu_n uv^A}{Q^2}\Big)^{\frac{\omega-t-\frac{1}{2}}{2A}} 
	\cos\left(\Big(\frac{\mu_n uv^A}{hQ^2}\Big)^{\frac{1}{2A}}+\pi\left(\frac{t}{2}+ A \gamma- \mu\right)\right)
	\\
	&
	 + O\bigg( (\mu_n uv^A)^{\frac{\omega-m-1}{2A}+ \xi}\bigg) + O\big((\mu_n uv^A)^{\rho+
	\delta - c_\rho -U}\big)
	\Bigg\}
	e^{-v} v^x  dv = J_1+J_2+ J_3. \quad \text{(say)}	\label{estimation of I_2}
\end{align} 
To obtain a bound for the term $J_1$, it is enough to find a bound for the integral \normalcolor
\begin{align}
	\frac{1}{\Gamma(x+1)} A_t
	&\int_{0}^{\infty} 
	 \Big(\frac{\mu_n uv^A}{Q^2}\Big)^{\frac{\omega-t-\frac{1}{2}}{2A}} 
	e^{ i\left((h^{-1}Q^{-2}\mu_n uv^A)^{\frac{1}{2A}}+\pi\left(\frac{t}{2}+ A\gamma - \mu\right)\right)}
	e^{-v} v^x  dv
	\\
	&=
	A_t
	(-1)^{\frac{t}{2}+ A \gamma - \mu} 
	(\Big(\frac{\mu_n u}{Q^2}\Big)^{\frac{\omega-t-\frac{1}{2}}{2A}}
	\frac{1}{\Gamma(x+1)}
	\int_{0}^{\infty} 
	v^{x+\frac{\omega-t-\frac{1}{2}}{2}} 
	e^{ i\left((h^{-1}Q^{-2}\mu_n u)^{\frac{1}{2A}}\sqrt{v}\right)}
	e^{-v}
	dv,
\end{align}	
after rearranging and using the fact that 
$e^{i \pi}= -1$.
We make a change of variable from 
$v$ 
to 
$v^2$
and use \eqref{gamma of x plus alpha} to write
\begin{align}
		&\frac{1}{\Gamma(x+1)} A_t
	\int_{0}^{\infty} 
	\Big(\frac{\mu_n uv^A}{Q^2}\Big)^{\frac{\omega-t-\frac{1}{2}}{2A}} 
	e^{ i\left((h^{-1}Q^{-2}\mu_n uv^A)^{\frac{1}{2A}}+\pi\left(\frac{t}{2}+ A \gamma- \mu\right)\right)}
	e^{-v} v^x  dv
	\\
	&\sim
	2 A_t 
	(-1)^{\frac{t}{2}+ A \gamma - \mu} 
	\Big(\frac{\mu_n u}{Q^2}\Big)^{\frac{\omega-t-\frac{1}{2}}{2A}}
	\frac{x^{\frac{\omega-t-\frac{1}{2}}{2}+\frac{1}{2}}}{\Gamma\left(x+\frac{\omega-t-\frac{1}{2}}{2}+\frac{3}{2}\right)}
	\int_{0}^{\infty} 
	v^{2\left(x+\frac{\omega-t-\frac{1}{2}}{2}+\frac{1}{2}\right)} 
	e^{ i(h^{-1}Q^{-2}\mu_n u)^{\frac{1}{2A}}v}
	e^{-v^2}
	dv.\label{ist int in I_2}
\end{align}
To apply Lemma \ref{lem:bdd for integral I_2} to the integral on the right-hand side above, we must verify that
\begin{align}\label{condition to apply lem in I_2}
	\left(x+\frac{\omega-t-\frac{1}{2}}{2}+\frac{1}{2}\right)^\epsilon 
	< 
	(h^{-1}Q^{-2}\mu_n u)^{\frac{1}{2A}},
\end{align}
for some $\epsilon >0$,
as 
$x \rightarrow \infty$. 
Indeed, as 
$\omega, t$ 
are constants here, we have that the left-hand side of \eqref{condition to apply lem in I_2} is
$\ll_\epsilon x^{\epsilon}$.
We choose
$0 <\epsilon < \frac{\gamma_2}{2}$ (recall that $\gamma_2>0$),
and invoke the condition 
$\mu_n > \frac{1}{u}x^{A\gamma_2}$ 
imposed on the sum
$I_2$. 
This gives that the left-hand side of \eqref{condition to apply lem in I_2} is 
\begin{align}
	\ll_{\epsilon'} x^{\frac{\gamma_2}{2}- \epsilon'}
	\ll_{\epsilon'} (\mu_n u)^{\frac{1}{2A}- \epsilon'},
\end{align}
for any 
$\epsilon' > 0$.
As 
$h, Q, A$ 
are fixed,
the bound \eqref{condition to apply lem in I_2} is evident.
We are now ready to use Lemma \ref{lem:bdd for integral I_2} with 
$u$ 
and 
$x$
replaced by 
$(h^{-1}Q^{-2}\mu_n u)^{\frac{1}{2A}}$
and 
$x+\frac{\omega-t-\frac{1}{2}}{2}+\frac{1}{2}$ 
respectively. This gives that the right-hand side of \eqref{ist int in I_2} is
\begin{align}
	&\ll_{M_1}
	(\mu_n u)^{\frac{\omega-t-\frac{1}{2}}{2A}}
	x^{\frac{\omega-t-\frac{1}{2}}{2}+\frac{1}{2}}
	(\mu_n u)^{\frac{-M_1}{2A}},
\end{align}
for any 
$M_1 >0$.
Using 
$\mu_n > \frac{1}{u}x^{A\gamma_2}$,
this is
\begin{align}
	\ll_{M_1}
	(\mu_n u)^{\frac{\omega-t-\frac{1}{2}}{2A}+ \frac{1}{A \gamma_2}\left|\frac{\omega-t-\frac{1}{2}}{2}+\frac{1}{2}\right|- \frac{M_1}{2A} }.
\end{align}
 Choosing $M_1$  to be sufficiently large, we obtain
\begin{align}\label{bound of first term in I_2}
	\frac{1}{\Gamma(x+1)} A_t
	\int_{0}^{\infty} 
	\Big(\frac{\mu_n uv^A}{Q^2}\Big)^{\frac{\omega-t-\frac{1}{2}}{2A}} 
	e^{ i\left((h^{-1}Q^{-2}\mu_n uv^A)^{\frac{1}{2A}}+\pi\left(\frac{t}{2}+ A \gamma- \mu\right)\right)}
	e^{-v} v^x  dv
	\ll_{M}
	(\mu_n u)^{-M},
\end{align}
for any $M >0$.
To estimate $J_2$, we see that the integral is
\begin{align}
	\frac{1}{\Gamma(x+1)}
	 (\mu_n u)^{\frac{\omega-m-1}{2A}+\xi}
	 \int_{0}^{\infty} 	
	 e^{-v} v^{x+\frac{\omega-m-1}{2}+\xi A}  dv
	 =
	 \frac{1}{\Gamma(x+1)}
	 (\mu_n u)^{\frac{\omega-m-1}{2A}+ \xi}
	 \Gamma\left(x+\frac{\omega-m+1}{2}+\xi A\right),
\end{align}
using the definition of the Gamma function.
Using \eqref{gamma of x plus alpha}
and for large enough 
 $m$, we have
 \begin{align}
 	\frac{1}{\Gamma(x+1)}
 	(\mu_n u)^{\frac{\omega-m-1}{2A}+\xi}
 	\Gamma\left(x+\frac{\omega-m+1}{2}+ \xi A\right)
 	&\ll_{M_1}
 	(\mu_n u)^{-M_1}
 	x^{-M_1}
 	\\
 	&
 	\ll_{M_1} 
 	(\mu_n u)^{-M_1},\label{bound of error in I_2}
 \end{align}
for any $M_1>0$,
since $x^{-M_1}\leq 1$ for $x\geq 1$.
\normalcolor
Combining \eqref{estimation of I_2}, \eqref{bound of first term in I_2} and \eqref{bound of error in I_2}, we get
\begin{align}
	J_1+J_2
	&\ll_{M_1}
	\sum_{ \mu_n > \frac{1}{u}x^{A\gamma_2}}	
	\frac{b_n}{\mu_n^{\d + \r }} 
	(\mu_n u)^{-M_1}.	
	\end{align}
Taking 
$M_1>\sigma_b - \delta - \rho$,
we apply Lemma \ref{lem:tail of sum b_n} with
$x=\frac{1}{u}x^{A\gamma_2}$,
to obtain
	\begin{align}
	J_1+J_2
	&\ll_{M_1}
	u^{-M_1} 
	\left(\frac{x^{A\gamma_2}}{u}\right)^{\sigma_b-\delta-\rho -M_1} 
	\left(\log \frac{x^{A\gamma_2}}{u} \right)^{\ell-1}.
\end{align}
Since 
$u^{-1} \leq x^{\gamma_1 \epsilon}$,
we see that 
$\log \left( \frac{1}{u} x^{A\gamma_2} \right) \ll \log x$,
where the implied constant may depend on 
$\gamma_1$, $\gamma_2$, $A$ 
and 
$\epsilon$,
which are fixed quantities. This gives 
\begin{align}
	J_1+J_2
	\ll 
	u^{\delta+\rho -\sigma_b} 
	x^{A\gamma_2(\sigma_b-\delta-\rho -M_1)}
	\log x.
\end{align}
Using \eqref{u theta} and \eqref{condition on u}, 
we obtain
\begin{align}\label{final estimation I_2}
	J_1+J_2
	&\ll 
	x^{\gamma_1 \epsilon|\delta+\rho -\sigma_b|+ A\gamma_2(\sigma_b-\delta-\rho -M_1) }
	\log x
	\\
	&
	\ll_M x^{-M}
\end{align}
for any 
$M >0$
by taking 
$M_1$
sufficiently large in the penultimate step.
 For the term $J_3$, the integral is 
\begin{align}
	\frac{(\mu_n u)^{\delta + \rho -c_\rho -U}}{\Gamma(x+1)}
	\int_{0}^{\infty} 	
	e^{-v} v^{x+(\delta + \rho -c_\rho -U)A}  dv
	&
	=
	\frac{1}{\Gamma(x+1)}
	(\mu_n u)^{\delta + \rho -c_\rho -U}
	\Gamma\left(x+(\delta + \rho -c_\rho -U)A+ 1\right)
	\\
	&= 
	(\mu_n u)^{\delta + \rho -c_\rho -U}
	x^{(\delta + \rho -c_\rho -U)A},
\end{align}
upon using  \eqref{gamma of x plus alpha} and the definition of the Gamma function. Hence, we have 
\begin{align}
	J_3 
	= 
	u^{\delta + \rho -c_\rho -U}
	x^{(\delta + \rho -c_\rho -U)A}
	\sum_{ \mu_n > \frac{1}{u}x^{A\gamma_2}}
	b_n
	\mu_n^{-c_\rho -U}.
\end{align}
For sufficiently large  $U$, the series above is absolutely convergent and we have,
\begin{align}
	J_3 
	&\ll
	u^{\delta + \rho -c_\rho -U}
	x^{(\delta + \rho -c_\rho -U)A}
	\\
	&
	\ll
	x^{\gamma_1 \epsilon |\delta + \rho -c_\rho -U|+ (\delta + \rho -c_\rho -U)A},
\end{align} 
upon using \eqref{condition on u} and \eqref{u theta}. Further, taking sufficiently small $\gamma_1$ and large enough $U$, we have as $x\rightarrow \infty$,
 $$J_3 \ll x^{-M},$$ 
 for a suitable $M>0$. Hence, we have $x \rightarrow \infty$,
\begin{align}\label{term I_2}
	\frac{1}{\Gamma(x+1)}	I_2 \ll x^{-M},
\end{align}
for any $M >0$. \normalcolor
 Combining 
 \eqref{lhsrhs borel mean}, 
 \eqref{partion into integrals},
 \eqref{final estimation I_1},
 and 
\eqref{term I_2},
 recalling 
$\omega = A\d+(2A-1)\r$,
  we get,
\begin{align}
	&\frac{1}{\Gamma(x+1)}
	\int_{0}^{\infty}
	\left[
	A^{\r}(uv^A) - Q_{\r}(uv^A) 
	\right] e^{-v} v^x  dv
	\\ 
	&= 
	(1+o(1))
	A_0 	\omega' Q^{ \frac{1}{2A} + \frac{\rho}{A}}
	u^{\frac{\delta}{2} - \frac{\rho}{2A}- \frac{1}{4A} + \rho}
	x^{\frac{\delta A}{2} - \frac{\rho}{2}- \frac{1}{4} + A \rho} \times
	\\
	&\sum_{n=1}^{\infty} 
	b_n 
	\mu_n^{-\frac{\delta}{2}- \frac{\rho}{2A}- \frac{1}{4A} }
	e^{-D(\mu_n u)^{\frac{1}{A}}}
	\cos\left(\frac{(\mu_n u)^{\frac{1}{2A}} \sqrt{x}}{(hQ^2)^{\frac{1}{2A}}}+ \pi (A \gamma - \mu) \right)
	+  
	O \left(x^{\frac{\delta A}{2} - \frac{\rho}{2}- \frac{1}{2} + A \rho + \epsilon }\right). 
\end{align}

\end{proof}

\subsection{Deducing Theorem \ref{thm:main theorem} from Proposition \ref{prop: key proposition}}
Consider the equality \eqref{final estimation} for
 $X \leq x \leq X^{3}$ (say), 
for some large $X$. 
 Denoting the right hand side of \eqref{final estimation}
as 
$H(x)$,
 let
 \begin{align}\label{def:F(k)}
 	F(\sqrt{x}) \coloneqq 
 \sum_{n=1}^{\infty} 
 b_n 
 \mu_n^{-\frac{\delta}{2}- \frac{\rho}{2A}- \frac{1}{4A} }
 e^{-D(\mu_n u)^{\frac{1}{A}}}
 \cos\left(\frac{(\mu_n u)^{\frac{1}{2A}} \sqrt{x}}{(hQ^2)^{\frac{1}{2A}}}+ \pi (A \gamma - \mu) \right). 
 \end{align}
  Then for 
  $\sqrt{X} \leq x \leq X^{\frac{3}{2}},$ 
  we have 
  \begin{align}\label{def:H(x)}
  	H(x^2)
  	= (1+o(1))
  	A_0 	\omega' Q^{ \frac{1}{2A} + \frac{\rho}{A}}
  	u^{\frac{\delta}{2} - \frac{\rho}{2A}- \frac{1}{4A} + \rho}
  	x^{ A\d+(2A-1)\r - \frac{1}{2}} 
  	F(x) 
  	+ 
  	O(x^{A\d+(2A-1)\r-1+2\epsilon}).
  \end{align}\normalcolor
We would like  to apply Lemma \ref{lem:sourmelidis theorem} to $F(x)$ with
\begin{align}
	f(n)
	=
	b_n 
	\mu_n^{-\frac{\delta}{2}- \frac{\rho}{2A}- \frac{1}{4A} } 
	e^{-D(\mu_n u)^{\frac{1}{A}}}, \quad
	\lambda_n
	=
	\frac{(\mu_n u)^{\frac{1}{2A}}}{ (hQ^2)^{\frac{1}{2A}} 2\pi},\hspace{1mm}
	\text{and} \quad
	\beta=\pi (A \gamma - \mu). \label{notations for sound lem}
\end{align}
Let 
$N \geq 2$ 
be a parameter depending on $x$,  as stated in  Theorem \ref{thm:main theorem}.
By the choice of 
$\mathcal{M}_N$,
we see that
\begin{align}
	\mathcal{M}_N 
	\subseteq
	\left\{ m \in \mathbb{N} : \mu_m \in \left[\left(\frac{1}{2}\right)^{2A} \mu_N, \left(\frac{3}{2}\right)^{2A} \mu_N \right] \right\},
\end{align}
and there exists a function $G(N)$ such that  
\begin{align} \label{derivation of H4}
\frac{1}{\mu_{N}} 	\sum_{\substack{m \in \mathbb{N} \\ \mu_m \in \left[\left(\frac{1}{2}\right)^{2A} \mu_N, \left(\frac{3}{2}\right)^{2A} \mu_N \right]}} b_m
\ll
G({N})	
\ll
\frac{1}{|\mathcal{M}_N|} \sum_{m \in \mathcal{M}_N} b_m. 
\end{align}
\normalcolor
By our  choice of $\lambda_n$ in \eqref{notations for sound lem}, it follows that 
$\lambda_m \in \left[\frac{1}{2}\lambda_N, \frac{3}{2}\lambda_N\right]$  
for  each
$m \in \mathcal{M}_N$.
The parameter $Y$ in Lemma \ref{lem:sourmelidis theorem} will be chosen later to satisfy $Y\geq 1$ and 
$\frac{Y}{2}\geq \sqrt{X}$.

 We choose 
$u = \frac{1}{\mu_N}$
in \eqref{def:F(k)}. Note that this satisfies \eqref{condition on u in prop} due to the condition \ref{assumption 2}. 
We apply Lemma \ref{lem:sourmelidis theorem} with $\mathcal{M}$= $\mathcal{M}_N, M= |\mathcal{M}_N|$. 
Due to the choice of $Y$ and condition \ref{assumption 4} on $|\mathcal{M}_N|$, Lemma \ref{lem:sourmelidis theorem} yields that
there exists a point
$x \in  \big[\sqrt{X},  2X \big]$
\normalcolor 
such that  
\begin{align}
	|F(x)| &\geq 
	C \normalcolor
	\sum_{m \in \mathcal{M}_N} b_m 
	\mu_m^{-\frac{\delta}{2}- \frac{\rho}{2A}- \frac{1}{4A}}
	e^{-D\big(\frac{\mu_m}{\mu_N}\big)^{\frac{1}{A}}}
	+ O \Bigg( \frac{2^{|\mathcal{M}_N|}Y \log X}{X} \normalcolor \sum_{\substack{n \\ \lambda_n \leq 2 \lambda_N}} b_n 
	\mu_n^{-\frac{\delta}{2}- \frac{\rho}{2A}- \frac{1}{4A}}
	e^{-D\big(\frac{\mu_n}{\mu_N}\big)^{\frac{1}{A}}}\Bigg)
	\\
	&-\frac{8}{\pi Y (hQ^2)^{-\frac{1}{2A}}} \sum_{n} b_n 
	\mu_n^{-\frac{\delta}{2}- \frac{\rho}{2A}- \frac{1}{4A}}
	e^{-D\big(\frac{\mu_n}{\mu_N}\big)^{\frac{1}{A}}}=
	 T_1+T_2-T_3 \quad \text{(say)}. \label{after sounds lem}
\end{align}
We proceed to show that 
$T_2, T_3$ 
are error terms.
From the choice of  
$\lambda_n$
in \eqref{notations for sound lem}, the condition 
$\lambda_n \leq 2 \lambda_N$ 
in $T_2$
above can be replaced by $\mu_n \leq 2^{2A} \mu_N$.
Moreover, the decreasing function 
$e^{-D\big(\frac{\mu_n}{\mu_N}\big)^{\frac{1}{A}}}$
can be bounded above by 
$e^{-D\big(\frac{\mu_1}{\mu_N}\big)^{\frac{1}{A}}}$. 
We have
\begin{align}
	T_2 
	&\ll
 \frac{2^{|\mathcal{M}_N|}Y \log X}{X} \normalcolor \sum_{ \mu_n \leq 2^{2A} \mu_N} b_n 
	\mu_n^{-\frac{\delta}{2}- \frac{\rho}{2A}- \frac{1}{4A} }
	e^{-D\left(\frac{\mu_1}{\mu_N}\right)^{\frac{1}{A}}}
	\\
	&\ll
 \frac{2^{|\mathcal{M}_N|}Y \log X}{X} \normalcolor
	\sum_{ \mu_n \leq 2^{2A} \mu_N} b_n 
	\mu_n^{-\frac{\delta}{2}- \frac{\rho}{2A}- \frac{1}{4A} },
\end{align}
as $N \rightarrow \infty$.
The condition 
$\r <  (2\s_b-\d)A-\frac{1}{2}$ 
is equivalent to
$\frac{\delta}{2}+ \frac{\rho}{2A}+ \frac{1}{4A} < \sigma_b$,
so that we can apply Lemma \ref{lem:asymp of b_n with log} with 
$x=2^{2A} \mu_N $ to get
\begin{align}
	T_2
	&\ll
 \frac{2^{|\mathcal{M}_N|}Y \log X}{X} \normalcolor \mu_N^{ \sigma_b -\frac{\delta}{2}- \frac{\rho}{2A}- \frac{1}{4A}}
	(\log \mu_{N})^{\ell-1}, 
	\label{estimation of T_2}
\end{align}
as $N \rightarrow \infty$.
 Further, using Hypothesis \ref{H3} to represent 
 $\mu_N^{ \sigma_b}
 (\log \mu_{N})^{l-1}$  as a sum, we write
 \begin{align}
 	T_2
 	&\ll
 \frac{2^{|\mathcal{M}_N|}Y \log X}{X} \normalcolor \mu_N^{ -\frac{\delta}{2}- \frac{\rho}{2A}- \frac{1}{4A}}
 	 \sum_{\substack{m \in \mathbb{N} \\ \mu_m \in \left[\left(\frac{1}{2}\right)^{2A} \mu_N, \left(\frac{3}{2}\right)^{2A} \mu_N \right]}} b_m
 	 \\
 	 &\ll
 \frac{2^{|\mathcal{M}_N|}Y \log X}{ X} \normalcolor \mu_N^{ -\frac{\delta}{2}- \frac{\rho}{2A}- \frac{1}{4A}+1} G(N) ,
 	 \label{estimation of T_2 after H3}
 \end{align}
 upon using \eqref{derivation of H4}. We now turn our attention to $T_3$.
 We apply Lemma \ref{lem: bdd for error E_1} with 
 $u=\frac{(2D)^A}{\mu_{N}}$ 
 and a fixed
 $x_0 \in (\frac{\mu_1}{2}, \mu_1)$ and obtain for any $\epsilon >0$,
\begin{align}\label{estimation of T_3}
	T_3 
	&\ll 
	\frac{\mu_{N}^{\sigma_b-\frac{\delta}{2}- \frac{\rho}{2A}- \frac{1}{4A}+\epsilon}}{Y}
	e^{-\frac{D}{2}\big(\frac{x_0}{\mu_N}\big)^{\frac{1}{A}}}
	\\
	&\ll
	\frac{\mu_{N}^{\sigma_b-\frac{\delta}{2}- \frac{\rho}{2A}- \frac{1}{4A}+\epsilon}}{Y},
\end{align}
as $N \rightarrow \infty$.
\normalcolor
We will proceed to obtain a lower bound for  
$T_1$. 
 We have for
$\mu_m \in \left[\left(\frac{1}{2}\right)^{2A} \mu_N, \left(\frac{3}{2}\right)^{2A} \mu_N \right]$, 
\begin{align}
	T_1 
	&=
	C 
	\sum_{m \in \mathcal{M}_N} b_m 
	\mu_m^{-\frac{\delta}{2}- \frac{\rho}{2A}- \frac{1}{4A} }
	e^{-D\big(\frac{\mu_m}{\mu_N}\big)^{\frac{1}{A}}}
	\\
	&\geq 
	C
	\left(\left(\frac{3}{2}\right)^{2A} \mu_N\right)^{-\frac{\delta}{2}- \frac{\rho}{2A}- \frac{1}{4A} }
	e^{-\frac{9D}{4}}
	\sum_{m \in \mathcal{M}_N} b_m 
	\\
	&
	\gg
	\mu_N^{-\frac{\delta}{2}- \frac{\rho}{2A}- \frac{1}{4A}}
	|\mathcal{M}_N| G({N}), \label{estimation of T_1}
\end{align}
upon using \eqref{derivation of H4}
as $N \rightarrow \infty$.
From \eqref{after sounds lem}, \eqref{estimation of T_2}, \eqref{estimation of T_3}, and \eqref{estimation of T_1},
 we obtain that there exists 
$x \in  \big[\sqrt{X},  2X \big]$
\normalcolor 
such that
\begin{align}\label{lower bdd after Ts' case 2}
	|F(x)| 
	&\gg
	\mu_N^{-\frac{\delta}{2}- \frac{\rho}{2A}- \frac{1}{4A}}
		|\mathcal{M}_N| G({N}) 
	+ 
	O\left(\frac{2^{|\mathcal{M}_N|}Y \log X}{ |\mathcal{M}_N| X} \normalcolor \mu_N^{ -\frac{\delta}{2}- \frac{\rho}{2A}- \frac{1}{4A}+1} |\mathcal{M}_N| G(N)\right)
	\\
	&+ 
	O \left( \frac{\mu_{N}^{\sigma_b-\frac{\delta}{2}- \frac{\rho}{2A}- \frac{1}{4A}+\epsilon}}{Y}
	\right).\normalcolor
\end{align}
 We choose $Y= 2 \sqrt{X}$ and 
 note that this choice of $Y$ satisfies 
 $Y \geq 1$ and   $\frac{Y}{2} \geq \sqrt{X}$. 
By conditions \ref{assumption 2} and \ref{assumption 4}, \normalcolor
the second term on the right-hand side of \eqref{lower bdd after Ts' case 2} is 
\begin{align}
	&=
	 o\left(\mu_N^{ -\frac{\delta}{2}- \frac{\rho}{2A}- \frac{1}{4A}}
	 |\mathcal{M}_N|G(N)\right),
	 \label{upper bdd of 2nd term}
\end{align}
as  $x \to \infty$.
Substituting this into \eqref{lower bdd after Ts' case 2},  we obtain that
there exists
$x \in  \big[\sqrt{X},  2X \big]$
\normalcolor 
 such that
\begin{align}
	|F(x)| 
	\gg
	\mu_N^{-\frac{\delta}{2}- \frac{\rho}{2A}- \frac{1}{4A}}
		|\mathcal{M}_N| G({N})
	+
	O \left( \frac{\mu_{N}^{\sigma_b-\frac{\delta}{2}- \frac{\rho}{2A}- \frac{1}{4A}+\epsilon}}{\sqrt{X}}\right). \label{lower bdd of F after L} \normalcolor
\end{align}
 From \eqref{def:g(x)}, the first term above is 
$g(x)$ and  by condition \ref{assumption 2}, the second term is $o(1)$\normalcolor.  Since $g(x) \asymp h(x)$, where $h(x)$ has the form \eqref{conditions on g(x)} and $h(x) \to \infty$ as $x \to \infty$, we   obtain that 
there exists a point 
$x \in  \big[\sqrt{X},  2X \big]$
\normalcolor  
such that 
\begin{align}\label{final bdd F(k^2)}
	|F(x)| \gg g(x). 
\end{align}
 We substitute the bound from \eqref{final bdd F(k^2)} into \eqref{def:H(x)} to see that
there exists a point 
$x \in  \big[\sqrt{X},  2X \big]$ 
\normalcolor 
such that
\begin{align}
	|H(x^2)| \gg 
	u^{\frac{\delta}{2} - \frac{\rho}{2A}- \frac{1}{4A} + \rho}
	x^{ A\d+(2A-1)\r - \frac{1}{2}} 
	g(x)
	+ 
	O(x^{A\d+(2A-1)\r-1+ 2 \epsilon}).
\end{align}
Since the growth of $g(x)$ is of the form \eqref{conditions on g(x)}, we note that 
$g(\sqrt{x}) \gg g(x)$, which means that  there exists
$x \in \left[X,  4X^{2}\right]$ 
\normalcolor
such that 
\begin{align}
	|H(x)| \coloneqq
	 \frac{1}{\Gamma(x+1)} 
\left|	\int_{0}^{\infty}
	\left[
	A^{\r}(uv^A) - Q_{\r}(uv^A) 
	\right] e^{-v} v^x  dv
	\right|
	&\gg 
	u^{\frac{\delta}{2} - \frac{\rho}{2A}- \frac{1}{4A}+ \rho}
	x^{ \frac{A \delta}{2}- \frac{\rho}{2} - \frac{1}{4}+A \rho} g(x)
	\\
	&+ O \left(x^{\frac{A \delta}{2}- \frac{\rho}{2} - \frac{1}{2}+A \rho + \epsilon}\right).
	\label{Hx}
\end{align}
As $u$ is sufficiently small, we have 
$u^{|\theta|} \leq u^{\theta}$ 
for any
$\theta \in \mathbb{R}$. Using this and \eqref{condition on u}, we have 
 \begin{align}
 	u^{\frac{\delta}{2}+ \rho - \frac{\rho}{2A}- \frac{1}{4A}}
 	x^{ \frac{A \delta}{2}- \frac{\rho}{2} - \frac{1}{4}+A \rho} g(x) 
 	&\gg 
 	x^{-\gamma_1 \epsilon \left|\frac{\delta}{2} - \frac{\rho}{2A}- \frac{1}{4A} + \rho\right| +  \frac{A \delta}{2}- \frac{\rho}{2} - \frac{1}{4}+A \rho  } 	g(x)
 	\\
 	&
 	\gg 
 	x^{  \frac{A \delta}{2}- \frac{\rho}{2} - \frac{1}{4}+A \rho- \epsilon'},
 \end{align}
for any $\epsilon' > 0$, 
by choosing 
$\gamma_1 >0$ to be
sufficiently small and using 
\eqref{conditions on g(x)} to estimate  the growth of 
$g(x)$.  The error term in \eqref{Hx} is thus negligible and we see that for some
$x \in \left[X,  4X^{2}\right]$,
\normalcolor
we have
\begin{align}
		\frac{1}{\Gamma(x+1)} 
	\left|\int_{0}^{\infty}
	\left[
	A^{\r}(uv^A) - Q_{\r}(uv^A) 
	\right] e^{-v} v^x  dv
	\right|
	\gg
	u^{\frac{\delta}{2} - \frac{\rho}{2A}- \frac{1}{4A} + \rho}
	x^{  \frac{A \delta}{2}- \frac{\rho}{2} - \frac{1}{4}+A \rho} g(x).  \label{final inequality 2}
\end{align}
We will now show that
\begin{align} \label{omega bound}
	A^{\r}(x) - Q_{\r}(x) 
	=
	\Omega
	\left(
	x^{\frac{\delta}{2} - \frac{\rho}{2A}- \frac{1}{4A}+ \rho} g(x)
	\right).
\end{align}
Suppose this is not true, then it means that 
\begin{align}
	A^{\r}(t) - Q_{\r}(t) 
	=
	o
	\left(
	t^{\frac{\delta}{2} - \frac{\rho}{2A}- \frac{1}{4A} + \rho} g(t)
	\right),
\end{align}
as $t \rightarrow \infty$. 
In particular, there exists a constant 
$c>0$
such that
\begin{align}\label{upper bdd assumption}
		|A^{\r}(t) - Q_{\r}(t)| 
		<
		c t^{\frac{\delta}{2} - \frac{\rho}{2A}- \frac{1}{4A}+ \rho} g(t),
\end{align}
for all
$t$
sufficiently large.
We will be done if we show that 
$c$
cannot be taken arbitrarily small. In order to do this, we will demonstrate a lower bound for 
$c$.
The bound \eqref{upper bdd assumption} gives
\begin{align}
		\frac{1}{\Gamma(t+1)} 
	\left|\int_{0}^{\infty}
	\left[
	A^{\r}(uv^A) - Q_{\r}(uv^A) 
	\right] e^{-v} v^t  dv \right|
	&<
	\frac{1}{\Gamma(t+1)} 
	\int_{0}^{\infty}
	c
	(uv^A)^{\frac{\delta}{2} - \frac{\rho}{2A}- \frac{1}{4A} + \rho} g(uv^A)
	 e^{-v} v^t  dv
	 \\
	 &
	 \leq
	 c c_1 
	 u^{\frac{\delta}{2} - \frac{\rho}{2A}- \frac{1}{4A} + \rho}
	 t^{ \frac{A \delta}{2}- \frac{\rho}{2} - \frac{1}{4}+A \rho}
	 g(ut^A), \label{final inequality}
\end{align}
for some constant 
$c_1 >0$,
as
$t \rightarrow \infty$, 
upon using Lemma \ref{lem: asymp for I_k,u}  and the fact that $g(x) \asymp h(x)$. 
We now let 
$t=x$ 
and observe that, due to the range \eqref{condition on u in prop} of 
$u$, we have 
\begin{align}
	x^{(1- \epsilon_0)A} \leq ux^A \leq x^A,
\end{align}
which implies
$g(ux^A) \asymp g(x)$, 
since 
$g(x) \asymp h(x)$, which
is of the form \eqref{conditions on g(x)}. 
This allows us to compare \eqref{final inequality} with \eqref{final inequality 2}.
 Denoting the implied constant in \eqref{final inequality 2} by 
$c'$,
we obtain
\begin{align}
	c' u^{\frac{\delta}{2} - \frac{\rho}{2A}- \frac{1}{4A} + \rho}
	x^{ \frac{A \delta}{2}- \frac{\rho}{2} - \frac{1}{4}+A \rho} g(x)
	<
	c c_2
	  u^{\frac{\delta}{2} - \frac{\rho}{2A}- \frac{1}{4A} + \rho}
	x^{  \frac{A \delta}{2}- \frac{\rho}{2} - \frac{1}{4}+A \rho}
	g(x),
\end{align}
for some constant $c_2>0$,
which further implies
\begin{align}\label{existence of lbdd of c}
	c> 
\frac{c'}{c_2} > 0.
\end{align}
This proves the omega bound \eqref{omega bound}.

It now remains to prove the 
$\Omega_\pm$-bounds 
\eqref{special case 1 thm} and \eqref{special case 2 thm}. In the asymptotic formula \eqref{final estimation}, since the integral is absolutely convergent, \eqref{final estimation} yields
\begin{align}
	\frac{1}{\Gamma(x+1)}
	\int_{0}^{\infty}
	&  
	\Re\left[\overline{\omega'}\big(
	A^{\r}(uv^A) - Q_{\r}(uv^A)\big)
	\right]  \normalcolor
	e^{-v} v^x  dv
	\\
	&=
	(1+o(1))
	A_0  Q^{ \frac{1}{2A} + \frac{\rho}{A}}
	u^{\frac{\delta}{2} - \frac{\rho}{2A}- \frac{1}{4A} + \rho}
	x^{  \frac{A \delta}{2}- \frac{\rho}{2} - \frac{1}{4}+A \rho} 
	F(\sqrt{x})
	+ 
	O \left(x^{ \frac{A \delta}{2}- \frac{\rho}{2} - \frac{1}{2}+A \rho + \epsilon}\right),
\end{align}
where 
$F(\sqrt{x})$ is as defined in \eqref{def:F(k)}.
Note that,  under hypothesis \ref{H2}, the explicit expression of $A_0$ in \cite[pg 23]{rsv} shows that $A_0>0$.
We choose $u = \frac{1}{\mu_N}$ and apply
Lemma \ref{lem:sourmelidis theorem} after \eqref{notations for sound lem} to obtain that
there exists a point 
$x \in \left[\sqrt{X}, X\right]$ 
\normalcolor
such that 
\begin{align}
	F(x) &\geq 
	2C
	\sum_{m \in \mathcal{M}_N} b_m 
	\mu_m^{-\frac{\delta}{2}- \frac{\rho}{2A}- \frac{1}{4A}}
	e^{-D\big(\frac{\mu_m}{\mu_N}\big)^{\frac{1}{A}}}
	+ O \Bigg(  \frac{2^{|\mathcal{M}_N|}Y \log X}{X} \normalcolor \sum_{n \geq 1} b_n 
	\mu_n^{-\frac{\delta}{2}- \frac{\rho}{2A}- \frac{1}{4A}}
	e^{-D\big(\frac{\mu_n}{\mu_N}\big)^{\frac{1}{A}}}\Bigg),
\end{align}
whenever
$A \gamma - \mu \equiv 0 \pmod{2}$. 
Following an argument similar to that discussed above from \eqref{after sounds lem} up to \eqref{final inequality 2}, we see that for 
$X \geq 2$ 
there exists a point 
$x \in \left[X, X^{2}\right]$ \normalcolor
such that 
\begin{align}
	\frac{1}{\Gamma(x+1)} 
	\int_{0}^{\infty}
\Re\left[\overline{\omega'}\big(
	A^{\r}(uv^A) - Q_{\r}(uv^A)\big)
	\right]  \normalcolor e^{-v} v^x  dv
	\geq 
	c_1
	u^{\frac{\delta}{2} - \frac{\rho}{2A}- \frac{1}{4A} + \rho}
	x^{  \frac{A \delta}{2}- \frac{\rho}{2} - \frac{1}{4}+A \rho}
	g(x), 
\end{align}
for some constant $c_1$.
In  order to show
\begin{align}\label{omega plus result}
	\Re \bigg(\overline{\omega'}\big(A^{\rho}(x)-Q_{\r}(x)\big)\bigg) \normalcolor
= \Omega_+\left(x^{\frac{\delta}{2} - \frac{\rho}{2A}- \frac{1}{4A} + \rho  } g(x)\right),
\end{align}
we again assume that there is a constant 
$c>0$ such that
\begin{align}
	\Re \bigg(\overline{\omega'}\big(A^{\rho}(x)-Q_{\r}(x)\big)\bigg) \normalcolor
	<
	c x^{\frac{\delta}{2} - \frac{\rho}{2A}- \frac{1}{4A} + \rho} g(x), 
\end{align}
for all sufficiently large values of $x$ and follow the method discussed above from \eqref{upper bdd assumption} to \eqref{existence of lbdd of c},  to demonstrate a lower bound for 
$c$.
Similarly, whenever 
$A\gamma - \mu \equiv 1 \pmod {2 } $,
we apply the relevant case of Lemma \ref{lem:sourmelidis theorem} and follow similar arguments  to obtain 
\begin{align}
	\Re \bigg(\overline{\omega'}\big(A^{\rho}(x)-Q_{\r}(x)\big)\bigg) \normalcolor= \Omega_-\left(x^{\frac{\delta}{2} - \frac{\rho}{2A}- \frac{1}{4A} + \rho  } g(x)\right). 
\end{align}

\qed

In the following sections, we give applications of Theorem \ref{thm:main theorem} to specific arithmetical functions.

\section{The Piltz Divisor Function: Proof of Theorem \ref{thm: piltz}}\label{sec:piltz div function}
Let $\mathbb{K} $ be an algebraic number field of degree $m=r_1+2r_2$, where $r_1$ is the number of real embeddings and $2r_2$ is the number of complex embeddings of $\mathbb{K}$. Let $\mathcal{O_{\mathbb{K}}}$ denote the ring of algebraic integers of $\mathbb{K}$. The Dedekind zeta function associated with $\mathbb{K}$ is defined by
\begin{equation}
	\zeta_{\mathbb{K}}(s) 
	:= 
	\sum_{\mathfrak{a} \subseteq \mathcal{O}_{\mathbb{K}} }
	\frac{1}{\mathcal{N}(\mathfrak{a})^s} \quad (\Re(s) > 1),
\end{equation}
where the sum is taken over non-zero integral ideals $\mathfrak{a}$ of $\mathbb{K}$ and $\mathcal{N}(\mathfrak{a})$ denotes the norm of $\mathfrak{a}$. Denoting $d_{\mathbb{K}}(n)$ to be  the number of non-zero ideals in $\mathcal{O}_{\mathbb{K}}$ with norm $n$, $\zeta_{\mathbb{K}}(s)$ can also be expressed as 
\begin{equation}
	\zeta_{\mathbb{K}}(s) 
	=
	\sum_{n=1}^{\infty}
	\frac{d_{\mathbb{K}}(n)}{n^s} \quad (\Re(s) > 1).	
\end{equation}
For $k \in \mathbb N$, let the arithmetical function $d_{\mathbb{K},k}(n) $ be defined by the relation 
\begin{equation}
\zeta^k_{\mathbb{K}}(s) 
=
\sum_{n=1}^{\infty}
\frac{d_{\mathbb{K},k}(n)}{n^s} \quad (\Re(s) > 1).
\end{equation}
Here $d_{\mathbb{K},k}(n)$ counts the number of $k$-tuples $(\mathfrak{a}_1, \hdots, \mathfrak{a}_k)$ of $\mathcal{O}_{\mathbb{K}}$ ideals with $\mathcal{N}(\mathfrak{a}_1 \cdots \mathfrak{a}_k)=n$. 
Note that when $\mathbb{K}=\mathbb Q$ (so, $r_1=1$, $r_2=0$), the function $d_{\mathbb{K},k}(n)$ is simply $d_k(n)$, that is, the number of ways to write $n$ as a product of $k$ positive integers. 

The Dirichlet series 
$ \zeta^k_{\mathbb{K}}(s)$
satisfies the functional equation 
\begin{align}
	c^{-s} 
	\Gamma\left(\frac{s}{2}\right)^{kr_1} 
	\Gamma(s)^{r_2}
	\zeta^k_{\mathbb{K}}(s)
	=
	c^{-(1-s)} 
	\Gamma\left(\frac{1-s}{2}\right)^{kr_1} 
	\Gamma(1-s)^{kr_2}
	\zeta^k_{\mathbb{K}}(1-s),
\end{align}
where $c$ is a constant depending on $\mathbb{K}$ and $k$.
One arrives at this functional equation by letting  
 $a_n=b_n= d_{\mathbb{K},k}(n)$,  $\lambda_n=\mu_n=n$,  
 $\d=1$, $Q=c^{-1}$, $\omega'=1$, 
 and 
 $r=kr_1+kr_2$, 
 with 
 $\a_i=\frac{1}{2}$  for $i=1, \hdots, kr_1$ 
 and 
 $\a_i =1$  for $i= kr_1+1, \hdots , kr_1+kr_2$, and $\b_i=0$ for $i= 1, \hdots, kr_1+kr_2 $ in  \eqref{generalfneqn}. 
Recalling the notation from Theorem \ref{thm:main theorem}, we have
\begin{align}
 A= \frac{kr_1}{2}+ kr_2 = \frac{km }{2}, \: \:
	 \mu = \frac{1- kr_1 - kr_2}{2} \: \text{ and }
	 \: \gamma=  -\left(\frac{1}{2}+ \frac{\rho}{km} + \frac{1}{2km}\right).
\end{align}


 Before proving  Theorem \ref{thm: piltz}, we give  the following lemma which  shows that the exponent of $\log_2 x$ in the omega bound of Theorem \ref{thm: piltz} is an improvement over \eqref{hafner piltz bdd}. 
 \begin{lemma}\label{lem: tau < tau'}
	Fix 
	$m,k \in \mathbb{N}$ with $km>1$.  
	Let 
	$\tau$ 
	and 
	$\tau'$
	be as defined in \eqref{tau in piltz} and \eqref{tau' in piltz} respectively. Then, for $0 \leq \rho < \frac{km-1}{2}$,  we have $\tau'>\tau$. 
\end{lemma}
\begin{proof}
Recall that 
	$I$ 
	is the set 
	$ \{ v \in \{1, 2, \cdots, m\}  : \d_v > 0\}$. Thus \eqref{tau in piltz} may be written as
	\begin{align}\label{tau in piltzproof}
		\tau = \frac{km-1- 2 \rho}{2km}\big(k \log k + k \sum_{v \in I} \delta_v v \log v-k +1\big)+k-1.
	\end{align}
	Using the values of 
	$\tau$ 
	and 
	$\tau'$
given in \eqref{tau in piltzproof} and \eqref{tau' in piltz}, we obtain after some simplification 
	\begin{align}
		\tau'-\tau
		&=
		\frac{km}{2km}
		\left(
		\sum_{v \in I} \delta_v v 
		\Big(v^{\frac{km-1-2 \rho}{km+1+2 \rho}} k^{\frac{2km}{km+1+2 \rho}} - k \log v
		\Big) - k \log k +k -2\right)
		\\
		&+
		\frac{1+2 \rho}{2km} 
		\left(
		\sum_{v \in I} \delta_v v 
		\Big(
		v^{\frac{km-1-2 \rho}{km+1+2 \rho}} k^{\frac{2km}{km+1+2 \rho}} + k \log v
		\Big) + k \log k -k 
		\right) -k +1.  
		\label{tau-tau' first}
	\end{align}
Since
	$0 \leq \rho < \frac{km-1}{2}$, 
we have
	$0 < km+1 + 2 \rho < 2km$,
	 we may then write
	\begin{align}\label{notation epsilon in piltz lemma}
		\frac{2km}{km+1+2 \rho} = 1+ \epsilon,
	\end{align}
	for some $\e >0$. In particular, 
	\begin{align}\label{notation2 epsilon in piltz lemma}
		\frac{km-1-2 \rho}{km+1+2 \rho}= \epsilon,
		\quad
		\frac{1+2 \rho}{2km} =\frac{1}{1+\e} - \frac{1}{2} = \frac{1- \epsilon}{2(1+\epsilon)}.
	\end{align}
	Using \eqref{notation epsilon in piltz lemma}, \eqref{notation2 epsilon in piltz lemma} and the property
	$\sum_{v \in I} \delta_v v = 1$, 
	 \eqref{tau-tau' first} can be rewritten as 
	\begin{align}
		\tau'-\tau
		&=
		\frac{1}{2}
		\sum_{v \in I} \delta_v v 
		\Big(v^{\epsilon} k^{1+\epsilon} - k \log v
		- k \log k +k -2\Big)
		\\
		&+
		\frac{1- \epsilon}{2(1+\epsilon)}
		\sum_{v \in I} \delta_v v 
		\Big(
		v^{\epsilon} k^{1+\epsilon} + k \log v
		+ k \log k -k \Big)
		+\sum_{v \in I} \delta_v v (1-k). 
		\label{tau-tau' second}
	\end{align}
	Since 
	$\frac{1- \epsilon}{2(1+\epsilon)} = \frac{1}{2}- \frac{\epsilon}{1+ \epsilon}$, we rearrange and simplify the right hand side of \eqref{tau-tau' second} to obtain
	\begin{align}
		\tau'-\tau 
		&= 
		\sum_{v \in I} \delta_v 
		v^{1+\epsilon} k^{1+\epsilon}
		-\frac{ \epsilon}{1+\epsilon}
		\sum_{v \in I} \delta_v v 
		\Big(
		v^{\epsilon} k^{1+\epsilon} + k \log v
		+ k \log k -k\Big)
		-\sum_{v \in I} \delta_v v k
		\\
		&=
		\sum_{v \in I}
		\frac{\delta_v v k}{1+\epsilon}
		\Big( (vk)^{\epsilon}  - \epsilon  \log (kv)-1\Big).
		\label{tau-tau' final}
	\end{align}
Using the fact that
	$e^x > 1+x$ for all $x >0$, we deduce that 
$
		(vk)^{\epsilon} - \epsilon  \log (kv)-1  
		> 0,
$
	whenever 
	$kv >1$. 
	Thus, from \eqref{tau-tau' final} we conclude
	$\tau'-\tau >0 $ 
	whenever 
	$k \geq 2$. If 
	$k=1$ and $m>1$, we have 
	$(vk)^{\epsilon}  - \epsilon  \log (kv)-1 > 0$
	for $v\ge 2$, thus ensuring  $\tau'-\tau >0 $ .
\end{proof}


\subsection{Proof of Theorem \ref{thm: piltz}}
 We proceed to verify the hypotheses of Theorem \ref{thm:main theorem}. Hypotheses \ref{H1} and \ref{H2} are easily seen to be true.
For the Piltz divisor function, it is  known that 
\begin{align}\label{asym for piltz div fn}
	\sum_{n \leq x} d_{\mathbb{K},k}(n)
     \sim
     c x (\log x)^{k-1},
\end{align}
where $c$ is a constant depending on $\mathbb{K}$ and $k$.
Thus, for $b_n= d_{\mathbb K, k}(n)$,  \eqref{H3} holds with 
$\sigma_b = 1$
and 
$\ell = k$. 
Let $N$ be a parameter depending on $x$, to be chosen later such that the conditions \ref{assumption 1} to \ref{assumption 4} of Theorem \ref{thm:main theorem} are met. 
 We will choose the set of distinguished integers $\mathcal M_N$ following the arguments of Girstmair, K\"uhleitner, M\"uller and Nowak  \cite[p. 192]{Girstmair}. 
   We choose 
 $\mathcal{M}_N$
 to be the set of 
 integers 
 $n \in \left[\frac{1}{2^{km}}N, \left(\frac{3}{2}\right)^{km} N\right]$
such that 
$n$ 
has 
$\lfloor \lambda_v \log_2N \rfloor$ 
distinct prime factors from 
$P_v$ 
for all 
$v \in I$ 
and no others. 
 Here $\lambda_v$ is a positive real number that  will be chosen later.   
 From (9) of   \cite{Girstmair}, we have 
\begin{equation}\label{size of d(n) in piltz}
	d_{\mathbb{K},k}(n) \geq \prod_{v \in I}  (vk)^{\lfloor \lambda_v \log_2N \rfloor} \gg (\log N)^{\sum_{v \in I} \lambda_v \log(vk)}
\end{equation}
for all 
$n \in \mathcal{M}_N$. Moreover, from  \cite[eq. (10)]{Girstmair}, the size of the set $\mathcal{M}_N$ is given by 
\begin{equation}\label{size of W(N) in piltz}
	|\mathcal{M}_N| \asymp N (\log N)^{\a} (\log_2N)^{-\frac{R}{2}},
\end{equation}
where 
$\a = -1 + \sum_{v \in I} \lambda_v (1+ \log \d_v - \log \lambda_v)$ and $R= |I|$.
We proceed to verify that condition \eqref{condition of H4} (with 
$Y= \frac{N}{2^{km}}$ and $c= 3^{km}$) holds. 
From \eqref{size of d(n) in piltz}, we see that 
\begin{align}\label{sum d(n) in piltz}
	\sum_{n \in \mathcal{M}_N} d_{\mathbb{K},k}(n) 
	\gg
	 |\mathcal{M}_N| (\log N)^{\sum_{v \in I} \lambda_v \log(vk)}.
\end{align}
On the other hand,  from \eqref{asym for piltz div fn}, we deduce that for fixed $k,m$, 
\begin{align}\label{sum d(n) in piltz in interval}
	\sum_{n \in \left[\frac{N}{2^{km}}, \left(\frac{3}{2}\right)^{km} N\right] } d_{\mathbb{K},k}(n)
	\sim
	c_1 
	N (\log N)^{k-1},
\end{align}
where $c_1$ is a constant depending on $ k$ and the field $\mathbb K$.
Taking 
$G(N)=c_2 (\log N)^{\sum_{v \in I} \lambda_v \log(vk)}$ and
combining \eqref{sum d(n) in piltz} and \eqref{sum d(n) in piltz in interval}, we see  that condition \eqref{condition of H4} 
 holds if we have 
\begin{align} 
	\label{need from lambda 2}
	 (\log N)^{k-1}
	  \ll
	  (\log N)^{\sum_{v \in I} \lambda_v \log(vk)}.
\end{align}
In what follows, we will choose $\lambda_v$ such that
$k-1 \leq \sum_{v \in I} \lambda_v \log(vk)$, so that  the above bound holds. The hypotheses of Theorem \ref{thm:main theorem} are thus satisfied and  applying it for $0 \leq \rho <  \frac{km-1}{2}$, we obtain 
 \begin{align}\label{omega for piltz with N}
	\Delta_{\mathbb{K}, k,\rho}(x)
	&= 
	\Omega
	\left(
	x^{\frac{1}{2}+ \rho -\frac{\rho}{km}- \frac{1}{2km}}
	N^{-\frac{1}{2}  -\frac{\rho}{km} - \frac{1}{2km}}
|\mathcal{M}_N|G(N)
	\right)
\\
&= 	\Omega
\left(
x^{\frac{1}{2}+ \rho -\frac{\rho}{km}- \frac{1}{2km}}
N^{-\frac{\r}{km}+\frac{1}{2}-\frac{1}{2km}}
(\log N)^{\a+\sum_{v \in I} \lambda_v \log(vk)}
(\log_2N)^{-\frac{R}{2}}
\right), 
\label{omega for piltz with N 2nd}
\end{align}
 using  
 $G(N)= (\log N)^{\sum_{v \in I} \lambda_v \log(vk)}$ and  \eqref{size of W(N) in piltz}. 
We consider  condition \ref{assumption 4} and use   \eqref{size of W(N) in piltz} to simplify the desired inequality 
\begin{align}\label{choice of N ineq in Piltz}
2^{|\mathcal{M}_N|} \leq x^{\epsilon},
\end{align}
\normalcolor
to  see that $N$ must satisfy 
$
	N \le c \log x (\log_2 x)^{-\a} (\log_3 x)^{\frac{R}{2}},
$
for a suitable constant $c>0$. In particular, we  make the choice 
\begin{align}
	N= c \log x (\log_2 x)^{-\a} (\log_3 x)^{\frac{R}{2}},
\end{align} 
which satisfies \ref{assumption 1}-\ref{assumption 4} of Theorem \ref{thm:main theorem}. 
 Substituting this value of $N$ into \eqref{omega for piltz with N 2nd}, we see that 	$\Delta :=\Delta_{\mathbb{K}, k,\rho}(x)$ satisfies 
\begin{align} 	 
\Delta =	\Omega&
	\Bigg( x^{\theta} (\log x)^{\kappa}
	(\log_2 x)^{- \a \kappa}
	(\log_3 x)^{\frac{R}{2} \kappa}
 \left\lbrace\log_2 x - \a \log_3 x + \frac{R}{2} \log_4 x\right\rbrace^{\a+\sum_{v \in I} \lambda_v \log(vk)}	
	\\
	\times&  \left\lbrace 
	\log \Big(\log_2 x - \a \log_3 x + \frac{R}{2} \log_4 x\Big) \right\rbrace^{-\frac{R}{2}}\Bigg), \label{piltz bdd after notation change}
\end{align}
where
\begin{align} \label{theta kappa}
\theta = \frac{(km-1)(2 \rho +1)}{2km}\:  \text{ and }
\: 
\kappa = \frac{km-1- 2 \rho}{2km}.
\end{align}
Upon simplification, we obtain  
\begin{align} \label{lower bdd after notation change}
\Delta_{\mathbb{K}, k,\rho}(x) = \Omega \left(
	x^{\theta}
(\log x)^{\kappa}
(\log_2 x)^{- \a \kappa +\a+\sum_{v \in I} \lambda_v \log(vk)}
(\log_3 x)^{\frac{R}{2} (\kappa-1)}
 \right). 
\end{align}
 We choose $\lambda$ suitably such that the exponent of $\log_2 x$ is maximum. It is easy to see that this occurs for 
 \begin{align}
 	\lambda_v= \d_v \cdot (vk)^{\left(\frac{\r}{km}+\frac{1}{2}+\frac{1}{2km}\right)^{-1}} =  k^{\frac{2km}{km+ 2\rho+1}} \cdot
 	\delta_v \cdot  v^{\frac{2km}{km+ 2\rho+1}}, 
  \end{align}
which satisfies \eqref{need from lambda 2}  for $0 \leq \rho <  \frac{km-1}{2}$ since $2km/ (km+2\rho+1)>1$ and $\sum_{v \in I} \d_v v=1$.  
 With this choice of 
$\lambda$ 
we   obtain the desired result. 


We  discuss some explicit examples in the following sections. 

\subsection{The divisor function $d_k(n)$} \label{sec:gen div fn}
      Let 
$k \geq 2$ 
be an integer  and 
$d_k(n)$ 
denote the number of ways of expressing $n$ as a product of $k$ positive integers. 
Putting 
$\mathbb{K}=\mathbb{Q}, m=1, I=\{1\}$, and $\delta_1=1$
in Theorem \ref{thm: piltz} yields Corollary \ref{cor: d_k}.


For the interested reader, we briefly delineate some of  the salient  features of the $d_k(n)$ case, even though they follow from the  proof of the more general Theorem \ref{thm: piltz}. 
The set  of distinguished integers $\mathcal{M}_N$ in this case is    the set of
integers in the interval
$\left[\frac{1}{2^k}N, \left(\frac{3}{2}\right)^k N\right]$
having exactly 
$\lfloor\lambda \log_2 N\rfloor$ 
distinct prime factors, where $\lambda$ is a positive real number, whose optimal value turns out to be $ k^{\frac{2k}{k+2\rho+ 1}}$. In this case, it is possible to estimate the size of $\mathcal{M}_N$ using 
 Lemma \ref{lem: size having exactly k fact}, to get \begin{align}
	|\mathcal{M}_N|
	\asymp
	\frac{N (\log_2 N)^{\lfloor\lambda \log_2 N\rfloor-1}}{\log N (\lfloor\lambda \log_2 N\rfloor-1)!}.
\end{align}
Using Stirling's formula \eqref{stirling for real}, we simplify the above expression to get
\begin{align}\label{size of W(N) in d_k}
	|\mathcal{M}_N|
	\asymp
	\frac{N}{\sqrt{\log_2 N}}(\log N)^{\lambda-1-\lambda \log \lambda}.
\end{align}  
Since each 
$n \in \mathcal{M}_N$ 
has  exactly 
$\lfloor \lambda \log_2 N \rfloor$ 
distinct  prime factors, for 
$n \in \mathcal{M}_N$,
we have the inequality
$d_k(n) \geq k^{ \lfloor \lambda \log_2 N \rfloor} \gg  (\log N)^{\lambda \log k}$, which yields 
\begin{align}
	\sum_{n \in \mathcal{M}_N} d_k(n)
	&\gg 
	|\mathcal{M}_N|  (\log N)^{\lambda \log k} \label{bdd for sum d_k(n)}.
\end{align}
On the other hand, it is  known that 
\begin{align}\label{sum d_k(n) in interval}
	\sum_{n \in \left[\frac{1}{2^k}N, \left(\frac{3}{2}\right)^k N\right] } d_k(n) \sim 
	c_1N (\log N)^{k-1},
\end{align}
where 
$c_1 >0$ is a constant depending upon $k$. 
One can see that   condition \eqref{condition of H4} 
is satisfied for 
$G(N)= (\log N)^{\lambda \log k}$ (with 
$Y= \frac{N}{2^k}$ and $c= 3^k$). 
The other aspects of the computation are as in the proof of Theorem \ref{thm: piltz}.

 \subsection{Other special cases}

\begin{example}
	Let 
$\mathbb{K}$ 
be a normal extension of
$\mathbb{Q}$ 
of degree
$m \geq 2$.
Then the Dirichlet density
$\delta_v = 0$
for 
$1 \leq v <m$
and 
$\delta_m= \frac{1}{m}$ (cf. W. Narkiewicz \cite[p.357, Corollary 4]{Narkiewicz}). In this case 
$R=|I|=1$, since 
$\delta_v >0$
only for $v=m$. 
Hence, by using \eqref{final piltz div bdd}, for 
$0 \leq \rho <  \frac{km-1}{2}$, we have
\begin{align}\label{bdd in piltz normal ext}
	\Delta_{\mathbb{K}, k,\rho}(x)
	= 
	\Omega \Bigg(
	&x^{\theta}
	(\log x)^{\kappa}
	(\log_2 x)^{\tau'}
	(\log_3 x)^{\frac{\kappa-1}{2}}\Bigg),
\end{align}
where
$\theta, \kappa$
are as in \eqref{theta kappa} and 
\begin{align}
	\tau' = \frac{km+1+2 \rho}{2km} \Big( \frac{1}{m}(mk)^{\frac{2km}{km+1+2 \rho}}-1\Big).
\end{align}
In this case, the  bound \eqref{hafner piltz bdd} obtained by Hafner gives an omega bound with   the smaller exponent 
\begin{align}
	\tau = \frac{km-1- 2 \rho}{2km}\big(k \log k + k  \log m-k +1\big)+k-1.
\end{align}
for $\log_2 x$. 
\end{example}

 \begin{example}
Consider a number field 
	$\mathbb{K}$ 
	of degree 
	$m=3$
	whose discriminant
	$D$
	is not a perfect square. This is an example of a number field which is not a normal extension. In this case, we have (cf. Girstmair, K\"uhleitner, M\"uller and Nowak \cite[p. 204, Subsection 4.3]{Girstmair})
$
		\delta_1= \frac{1}{2}, \quad \delta_2= 0, \quad \delta_3= \frac{1}{6}
$. 
Thus, we have 
	$R=|I|=2$ 
	since both $\delta_1, \delta_3 >0$.
	By using \eqref{final piltz div bdd}, for 
	$0 \leq \rho <  \frac{3k-1}{2}$, we obtain
	\begin{align}\label{bdd in piltz non normal ext}
		\Delta_{\mathbb{K}, k,\rho}(x)
		= 
		\Omega \Bigg(
		&x^{\theta}
		(\log x)^{\kappa}
		(\log_2 x)^{\tau'}
		(\log_3 x)^{\kappa-1}\Bigg),
	\end{align}
	where 
	\begin{align}
		\theta = \frac{(3k-1)(2 \rho +1)}{6k},
\:\: 		
\kappa = \frac{3k-1- 2 \rho}{6k} 
\: \text{ and }\:
	\tau' = \frac{3k+1+2 \rho}{6k} \Big( \frac{1}{2} k^{\frac{6k}{3k+1+2 \rho}}+ \frac{1}{6} (3k)^{\frac{6k}{3k+1+2 \rho}}-1\Big).
	\end{align}
	In this case, the exponent of 
	$\log_2 x$ in Hafner's bound \eqref{hafner piltz bdd} is given by
	\begin{align}
		\tau= \frac{3k-1- 2 \rho}{6k}\big(k \log k + \frac{k  \log 3}{2}-k +1\big)+k-1. 
	\end{align}
\end{example}

\begin{example}
	Consider 
	$\mathbb{K} = \mathbb{Q}(\alpha)$, 
	where
	$\alpha$ is a root of an irreducible polynomial of degree $7$ over $\mathbb{Q}$ with 
	Galois group  
	PSL$(3,2)$ (this is the projective special linear group of degree $3$ over the field of two elements, and is a simple group of order 168). Here, the degree of extension is $m=7$. In this case, only the sets $P_1,P_3$ and  $P_7$ are nonempty. The corresponding Dirichlet densities are given by (cf. Girstmair, K\"uhleitner, M\"uller and Nowak \cite[p. 204, Subsection 4.4]{Girstmair})
$
		\delta_1= \frac{7}{12}, \quad \delta_3= \frac{1}{8}, \quad \delta_7= \frac{1}{168}
$. 
	So, we have 
	$R=|I|=3$ 
	since  $\delta_1, \delta_3, \delta_7  >0$.
	By using \eqref{final piltz div bdd}, for 
	$0 \leq \rho <  \frac{7k-1}{2}$, we obtain
	\begin{align}\label{bdd in piltz 7 deg ext}
		\Delta_{\mathbb{K}, k,\rho}(x)
		= 
		\Omega \Bigg(
		&x^{\theta}
		(\log x)^{\kappa}
		(\log_2 x)^{\tau'}
		(\log_3 x)^{\frac{3}{2}(\kappa-1)}\Bigg),
	\end{align}
	where 
	\begin{align}
		&\theta = \frac{(7k-1)(2 \rho +1)}{14k},
\: \: 
	\kappa = \frac{7k-1- 2 \rho}{14k} \: \: \text{ and }
		\\
		&	\tau' = \frac{7k+1+2 \rho}{14k} \Big( \frac{7}{12} k^{\frac{14k}{7k+1+2 \rho}}+ \frac{1}{8} (3k)^{\frac{14k}{7k+1+2 \rho}}+\frac{1}{168} (7k)^{\frac{14k}{7k+1+2 \rho}}-1\Big).
	\end{align}
	If we compare with Hafner's bound given in \eqref{hafner piltz bdd}, the exponent of 
	$\log_2 x$ in this case 
	is given by
	\begin{align}
		\tau=  \frac{7k-1- 2 \rho}{14k}\big(k \log k + \frac{  3k\log 3}{8}+ \frac{k \log 7}{24}-k +1\big)+k-1. 
	\end{align}
\end{example}

\section{The Circle problem: Proof of Theorem \ref{thm: r(n)}} \label{sec:circle prob}
      Let 
$r(n)$
denote the number of ways of expressing $n$ as a sum of two squares. The associated Dirichlet series  
$ \phi(s):=\sum_{n=1}^{\infty}r(n) n^{-s}$ 
satisfies the functional equation 
\begin{align}\label{fn eqn of r(n)}
	\pi^{-s}
	\Gamma(s)
	\phi(s)
	=
	\pi^{-(1-s)}
	\Gamma(1-s)
	\phi(1-s).
\end{align}       
 This functional equation follows from \eqref{generalfneqn} by letting 
$a_n=b_n=r(n), \lambda_n= \mu_n= n, \d=1, Q= \pi^{-1}, \omega'=1, r=1,$ 
and 
$\a_1=1, \b_1 =0$. 
Recalling the notation from Theorem \ref{thm:main theorem}, we have   
$A= 1$, 
$\mu= 0$,
and 
$\gamma=  -\left(\frac{3}{4}+ \frac{\rho}{2}\right)$.


We proceed to verify the hypotheses of Theorem \ref{thm:main theorem}. Equation \eqref{fn eqn of r(n)} shows that Hypotheses \ref{H1} and \ref{H2} hold.  For $b_n=r(n)$, \ref{H3} is known to hold with 
$\sigma_b = 1$, 
$\ell = 1$ and $c_0= \pi$. 
Let $N$ be a parameter depending on $x$, to be chosen later such that the conditions \ref{assumption 1}-\ref{assumption 4} of Theorem \ref{thm:main theorem} are met. 
We choose 
$\mathcal{M}_N$
to be the set of
integers in the interval
$\left[\frac{N}{4}, \frac{9N}{4}\right]$
having exactly 
$\lfloor\lambda \log_2 N\rfloor$ 
distinct prime factors all of which are $ \equiv 1 \pmod{4}$, where $\lambda$ is a positive real number to be chosen later. 
 The size of the set $\mathcal{M}_N$  is given by (see \cite[Equation (3.7)]{Sound})
\begin{align}\label{size of W(N) in r(n)}
	|\mathcal{M}_N|
	\asymp
	\frac{N}{\sqrt{\log_2 N}}(\log N)^{\lambda-1-\lambda \log \lambda - \lambda \log 2}.
\end{align}  
We proceed to check that condition \eqref{condition of H4} (with 
$Y= \frac{N}{4}$ and $c= 9$) holds for
 $G(N)=c_1 (\log N)^{\lambda \log 2}$.
Since each 
$n \in \mathcal{M}_N$ 
has  exactly 
$\lfloor \lambda \log_2 N \rfloor$ 
distinct  prime factors, 
we have the inequality
$r(n)  \gg  (\log N)^{\lambda \log 2}$ for 
$n \in \mathcal{M}_N$.
Using the inequality for $r(n)$, we have the estimate
\begin{align}
	\sum_{n \in \mathcal{M}_N} r(n)
	&\gg
	|\mathcal{M}_N|  (\log N)^{\lambda \log 2} \label{bdd for sum r(n)}.
\end{align}
It is also known that 
\begin{align}\label{sum r(n) in interval}
	\sum_{n \in \left[\frac{N}{4}, \frac{9N}{4}\right] }r(n) \sim 
	c_2N,
\end{align}
where 
$c_2 >0$ is a constant.
Taking 
$G(N)=c_1 (\log N)^{\lambda \log 2}$
and
combining  \eqref{bdd for sum r(n)} and \eqref{sum r(n) in interval}, we see that condition \eqref{condition of H4} 
holds since
$
	c_2 \leq  (\log N)^{\lambda \log 2}.
$
Hence, by Theorem \ref{thm:main theorem}, for 
$0 \leq \rho <  \frac{1}{2}$, we obtain 
\begin{align}\label{omega for r(n) with N}
	P_{\rho}(x)
	&= 
	\Omega
	\left(
	x^{\frac{1}{4} +\frac{\rho}{2}}
	N^{-\frac{3}{4}-\frac{\rho}{2}} 
	|\mathcal{M}_N|G(N)
	\right) 
	\\
	&= 
	\Omega
	\left(
	x^{\frac{1}{4} +\frac{\rho}{2}}
	\frac{N^{\frac{1}{4} -\frac{\rho}{2}}}{\sqrt{\log_2 N}}
	(\log N)^{\lambda-1-\lambda \log \lambda }
	\right), 
	\label{omega bdd final r(n)}
\end{align}
using   $G(N)= (\log N)^{\lambda \log 2}$
and 
\eqref{size of W(N) in r(n)}.

We consider  condition \ref{assumption 4} and use   \eqref{size of W(N) in r(n)} to simplify the desired inequality  
$2^{|\mathcal{M}_N|} \leq x^{\epsilon}$, \normalcolor
to  see that $N$ must satisfy 
\begin{align}
	N \le c \log x (\log_2 x)^{1-\lambda+\lambda \log \lambda+ \lambda \log 2} (\log_3 x)^{\frac{1}{2}},
\end{align}
for a suitable constant $c>0$. In particular, we  choose 
\begin{align}
	N= c \log x (\log_2 x)^{1-\lambda+ \lambda \log \lambda+ \lambda \log 2} (\log_3 x)^{\frac{1}{2}},
\end{align}
which satisfies conditions \ref{assumption 1}-\ref{assumption 4} of Theorem \ref{thm:main theorem}. 
Substituting this value of 
$N$ into \eqref{omega bdd final r(n)},  we obtain
\begin{align}
	P_{\rho}(x)
	=
	\Omega \left(
	x^{\frac{1}{4} +\frac{\rho}{2}}
	(\log x)^{\frac{1}{4} -\frac{\rho}{2}} 
	(\log_2 x)^{\left(\frac{\r}{2}+\frac{3}{4}\right)(\lambda-1 - \lambda \log \lambda) + \left( \frac{1}{4} -\frac{\rho}{2}\right) \lambda \log 2 }
	(\log_3 x)^{-\frac{\r}{4} - \frac{3}{8}} \right),
\end{align}
upon simplification.
We choose $\lambda$ suitably such that the exponent of $\log_2 x$ is maximum. It is easy to see that this occurs for
\begin{align}
	\label{choice lambda}
	\lambda = 2^{\frac{1-2 \rho}{2 \rho+3}}.
\end{align} 
With this choice of 
$\lambda$,
we   obtain the desired result. 

We do not have any $\Omega_\pm$  results in this case since the conditions in   
\eqref{special case 1 thm} and  \eqref{special case 2 thm} are not satisfied for any $0 \leq \rho < \frac{1}{2}$.

\section{The Sphere Problem: Proof of Theorem \ref{thm: r_3(n)}}\label{sec:sphere prob}
Let 
$r_3(n)$
denote the number of ways of expressing $n$ as a sum of three squares. The associated Dirichlet series  
$ \phi(s):=\sum_{n=1}^{\infty}r_3(n) n^{-s}$ 
satisfies the functional equation 
\begin{align}\label{fn eqn of r_3(n)}
	\pi^{-s}
	\Gamma(s)
	\phi(s)
	=
	\pi^{-\left(\frac{3}{2}-s\right)}
	\Gamma\left(\frac{3}{2}-s\right)
	\phi\left(\frac{3}{2}-s\right).
\end{align}       
This functional equation follows from \eqref{generalfneqn} by letting 
$a_n=b_n=r_3(n), \lambda_n= \mu_n= n, \d=\frac{3}{2}, Q= \pi^{-1}, \omega'=1, r=1,$ 
and 
$\a_1=1, \b_1 =0$. 
Recalling the notation from Theorem \ref{thm:main theorem}, we have   
$A= 1$, 
$\mu= 0$,
and 
$\gamma=  -\left(1+ \frac{\rho}{2}\right)$.

We proceed to verify the hypotheses of Theorem \ref{thm:main theorem}. Equation \eqref{fn eqn of r_3(n)} shows that Hypotheses \ref{H1} and \ref{H2} hold.  For $b_n=r_3(n)$, \ref{H3} is known to hold with 
$\sigma_b = \frac{3}{2}$, 
$\ell = 1$ and $c_0= \frac{4 \pi}{3}$. 
Let $N$ be a parameter depending on $x$, to be chosen later such that the conditions \ref{assumption 1}-\ref{assumption 4} of Theorem \ref{thm:main theorem} are met. 
We choose 
$\mathcal{M}_N$
to be the set of
all integers in the interval
$\left[\frac{N}{4}, \frac{9N}{4}\right]$.
The size of the set $\mathcal{M}_N$  is given by 
$
	|\mathcal{M}_N|
	\asymp
	N.
$
 Since
\begin{align}\label{sum r_3(n) in interval}
	\sum_{n \in \left[\frac{N}{4}, \frac{9N}{4}\right] }r_3(n) \sim 
	c_1N^{\frac{3}{2}},
\end{align}
for a positive  constant $c_1$, condition \eqref{condition of H4} (with 
$Y= \frac{N}{4}$ and $c= 9$) trivially holds for
$G(N)= c_1 N^{\frac{1}{2}}$.
Hence, by Theorem \ref{thm:main theorem}, for 
$0 \leq \rho <  1$, we obtain 
\begin{align}
	P_{3,\rho}(x)
	&= 
	\Omega
	\left(
	x^{\frac{1}{2} +\frac{\rho}{2}}
	N^{-1-\frac{\rho}{2}} 
	|\mathcal{M}_N|G(N)
	\right) 
	= 
	\Omega
	\left(
	x^{\frac{1}{2} +\frac{\rho}{2}}
	N^{\frac{1}{2} -\frac{\rho}{2}}
	\right), 
	\label{omega bdd final r_3(n)}
\end{align}
using   $G(N)=  c_1 N^{\frac{1}{2}}$
and 
$ |\mathcal{M}_N|
\asymp
N$.
We choose $N$ such that 
$2^{|\mathcal{M}_N|} \leq x^{\epsilon}$, \normalcolor
as needed  in condition \ref{assumption 4}.
This means that $N$ must satisfy 
$N\le c \log x$ 
for some suitable constant 
$c>0$. In particular, we choose 
$
	N =  c \log x,
$
which satisfies conditions \ref{assumption 1}-\ref{assumption 4} of Theorem \ref{thm:main theorem}.
Substituting this value of 
$N$ into \eqref{omega bdd final r_3(n)},  we obtain
\begin{align}
	P_{3,\rho}(x)=
	\Omega
	\left(
	x^{\frac{1}{2} +\frac{\rho}{2}}
	(\log x)^{\frac{1}{2} -\frac{\rho}{2}}
	\right).
\end{align}
We do not have any $\Omega_+$  result for this problem since the condition in   
\eqref{special case 1 thm} is not satisfied for any $0 \leq \rho < 1$. By \eqref{special case 2 thm}, the above bound holds for $\Re(P_{3,\rho}(x))$ with $\Omega$ replaced by 
$\Omega_-$ whenever $\rho=0$.

\section{Lattice points in an Ellipse: Proof of Theorem \ref{thm:r_Q(n)}}

 Let
$\mathcal{Q}(u,v):=u^2+5v^2$
be a positive definite binary quadratic form, and define
\begin{align}
	r_\mathcal{Q}(n):=\#\{(u,v)\in\mathbb Z^2:\mathcal{Q}(u,v)=n\}.
\end{align}
The associated Dirichlet series
\begin{align}
	\phi(s):=\sum_{n=1}^{\infty}\frac{r_\mathcal{Q}(n)}{n^s}
	=\sum_{(u,v)\in\mathbb Z^2\setminus\{(0,0)\}}\frac{1}{\mathcal{Q}(u,v)^s},
	\qquad \Re(s)>1,
\end{align}
 is the Epstein zeta function associated with the quadratic form $\mathcal {Q}$. 
 The Dirichlet series $\phi(s)$ satisfies the functional equation 
 \begin{align}\label{fn eqn of r_Q(n)}
 	\left(\frac{\sqrt 5}{\pi}\right)^s\Gamma(s)\phi(s)
 	&=
 	\left(\frac{\sqrt 5}{\pi}\right)^{1-s}\Gamma(1-s)\phi(1-s).
 \end{align}
 This functional equation can be obtained from \eqref{generalfneqn} by letting 
 $ a_n=b_n=r_\mathcal Q(n), \lambda_n=\mu_n=n,
 \delta=1, Q= \sqrt{5}{\pi}^{-1}, \omega'=1,r=1,$ and $\alpha_1=1, \beta_1=0$.
 Recalling the notation from Theorem \ref{thm:main theorem}, we have   
 $A= 1$, 
 $\mu= 0$,
 and 
 $\gamma=  -\left(\frac{3}{4}+ \frac{\rho}{2}\right)$.

We proceed to verify the hypotheses of Theorem \ref{thm:main theorem}. Equation \eqref{fn eqn of r_Q(n)} shows that Hypotheses \ref{H1} and \ref{H2} hold.
Further Hypothesis \ref{H3} holds with
$
\sigma_b=1, \ell=1, c_0=\frac{\pi}{\sqrt5}.
$

We now choose the set of distinguished integers $\mathcal{M}_N$ as follows. 
Consider the set of primes
$$
\mathcal P:=\{p\text{ prime}:p\equiv1,3,7,9\pmod{20}\}.
$$
Let $B_N$ be the set of integers $n\in[N/4,3N/4]$ having exactly $\lfloor\lambda\log_2N\rfloor$ distinct prime factors, all of which belong to $\mathcal P$, where $\lambda$ is a positive real number to be chosen later.
 That is,
\begin{align}
	B_N:=
	\left\{n\in\left[\frac N4,\frac{3N}{4}\right):
	\omega(n)=\lfloor\lambda\log_2N\rfloor,\, p\mid n\Longrightarrow p\equiv1,3,7,9\pmod{20}\right\},
	\label{ellipse-BN}
\end{align}
where $\omega(n)$ denotes the number of distinct prime divisors of $n$.
For an integer $n\in B_N$, let $t(n)$ be the  number of prime factors of $n$, counted with multiplicity, which are congruent to $3$ or $7$ modulo $20$.
By Lemma \ref{lem:count r_Q(n)}, we have for $n\in B_N$
\begin{align}
	r_\mathcal Q(n)
	&=
	\begin{cases}
		2d(n),&\text{if} \,\, t(n) \,\, \text{is even},\\
		0,&\text{if} \,\, t(n) \,\, \text{is odd}.
	\end{cases}
	\label{ellipse-representation}
\end{align}
To manage the case where $r_\mathcal Q(n)=0$, we define 
$T:[N/4,3N/4) \rightarrow [N/4,9N/4]$ by
\begin{align}
	T(n):=
	\begin{cases}
		n,&\text{if} \,\, t(n) \,\, \text{is even},\\
		3n,&\text{if} \,\, t(n) \,\, \text{is odd},
	\end{cases}
	\label{ellipse-map-T}
\end{align}
and take the set of distinguished integers $\mathcal M_N$ to be $T(B_N)$. Note that we have 
$r_\mathcal Q(n)=2d(n)$
 for all 
 $n \in \mathcal M_N$. 
 Since $T$ is injective, we have
\begin{align}\label{size of W(N) in r_Q(n)}
	|B_N|=|\mathcal M_N|
	\asymp
	\frac{N}{\sqrt{\log_2N}}
	(\log N)^{\lambda-1-\lambda\log\lambda-\lambda\log2},
\end{align}
by using Lemma \ref{lem:girstmair size}. 

We proceed to check that condition \eqref{condition of H4} (with 
$Y= \frac{N}{4}$ and $c= 9$) holds for
$G(N)=c_1 (\log N)^{\lambda \log 2}$.
Since each 
$n \in \mathcal{M}_N$ 
has  exactly 
$\lfloor \lambda \log_2 N \rfloor$ 
or 
$\lfloor \lambda \log_2 N \rfloor+1$ 
distinct  prime factors, 
we have the inequality
$r_\mathcal Q(n)  \gg  (\log N)^{\lambda \log 2}$ for 
$n \in \mathcal{M}_N$.
Using the inequality for $r_\mathcal{Q}(n)$, we have the estimate
\begin{align}
	\sum_{n \in \mathcal{M}_N} r_\mathcal{Q}(n)
	&\gg
	|\mathcal{M}_N|  (\log N)^{\lambda \log 2} \label{bdd for sum r_Q(n)}.
\end{align}
It is also known that 
\begin{align}\label{sum r_Q(n) in interval}
	\sum_{n \in \left[\frac{N}{4}, \frac{9N}{4}\right] }r_\mathcal{Q}(n) \sim 
	c_2N,
\end{align}
where 
$c_2 >0$ is a constant.
Taking 
$G(N)=c_1 (\log N)^{\lambda \log 2}$
and
combining  \eqref{bdd for sum r_Q(n)} and \eqref{sum r_Q(n) in interval}, we see that condition \eqref{condition of H4} 
holds since
$
c_2 \leq  (\log N)^{\lambda \log 2}.
$
Hence, by Theorem \ref{thm:main theorem}, for 
$0 \leq \rho <  \frac{1}{2}$, we obtain 
\begin{align}\label{omega for r(n) with N}
		P_\rho(x,\mathcal {Q})
	&= 
	\Omega
	\left(
	x^{\frac{1}{4} +\frac{\rho}{2}}
	N^{-\frac{3}{4}-\frac{\rho}{2}} 
	|\mathcal{M}_N|G(N)
	\right) 
	\\
	&= 
	\Omega
	\left(
	x^{\frac{1}{4} +\frac{\rho}{2}}
	\frac{N^{\frac{1}{4} -\frac{\rho}{2}}}{\sqrt{\log_2 N}}
	(\log N)^{\lambda-1-\lambda \log \lambda }
	\right), 
	\label{omega bdd final r_Q(n)}
\end{align}
using   $G(N)= (\log N)^{\lambda \log 2}$
and 
\eqref{size of W(N) in r_Q(n)}. 
The rest of the proof is exactly the same as in the circle problem after \eqref{omega bdd final r(n)}. Following that, we obtain Theorem \ref{thm:r_Q(n)}.

We do not have any $\Omega_\pm$  results in this case since the conditions in   
\eqref{special case 1 thm} and  \eqref{special case 2 thm} are not satisfied for any $0 \leq \rho < \frac{1}{2}$.

\normalcolor

\section{Concluding remarks}
Recently, Lamzouri \cite[Theorem~1.1]{lamzouri2026} improved the omega bounds in the divisor problem and the circle problem by removing the negative power of $\log_3 x$ in the earlier bounds of Soundararajan \cite{Sound}, Sourmelidis \cite{sourmelidis2024}, and Mahatab \cite{mahatab2025},  while also determining the signs of the
resulting large values. Subsequently, Karak and Mahatab \cite[Theorem~1.1]{karakmahatab2026} adapted Lamzouri's method to obtain analogous improvements for the Piltz divisor problem over number fields. Lamzouri's improvement is based on an
estimate for the  trigonometric series
$
 F(x)=\sum_{n=1}^{\infty}f(n) \cos (2 \pi \lambda_n x + \beta),
$
as defined in \eqref{def of F(x)}. 
In particular, Lamzouri's theorem states that for a finite set 
$\mathcal{M} \subset \mathbb{N}$ 
with cardinality 
$|\mathcal{M}|=M$ 
and any real numbers 
$3< Y < X$, we have
\begin{align}
	\max_{x \in [Y,X]}{} F(x) 
	\geq
	\delta r
	\sum_{m\in\mathcal M}f(m)
	+
	O\!\left(
	\left(\frac{1+r}{1-r}\right)^M
	\left(\frac{Y\log X}{X}\right)^{\alpha}
		\sum_{n \geq 1} f(n)
	\right),
\end{align} 
where
$0<r<1$
and
$
0<\delta<\cos\beta
$
are fixed numbers,
$
\alpha
=
\frac{2}{\pi}
\arctan\!\left(
\frac{\cos\beta-\delta}{|\sin\beta|}
\right)
$,
and
$
\beta\in\left(-\frac{\pi}{2},\frac{\pi}{2}\right)\setminus\{0\}
$.
This result suggests that we can replace Lemma \ref{lem:sourmelidis theorem} in the proof of Theorem \ref{thm:main theorem} and obtain an axiomatic bound for the error term 
$A^{\rho}(x)-Q_{\r}(x)$. In particular, by applying this bound for the trigonometric series 
$F(x)$ in \eqref{def:F(k)}, we obtain
\begin{align}\label{main term in lamz}
	\max_{x \in [Y,X]}{} F(x) 
 &\geq 
	\delta r
	\sum_{m \in \mathcal{M}} b_m 
	\mu_m^{-\frac{\delta}{2}- \frac{\rho}{2A}- \frac{1}{4A}}
	e^{-D(\mu_m u)^{\frac{1}{A}}}
	\\
	&+ O \Bigg( \left(\frac{1+r}{1-r}\right)^M
	\left(\frac{Y\log X}{X}\right)^{\alpha}
	 \sum_{n \geq 1} b_n 
	\mu_n^{-\frac{\delta}{2}- \frac{\rho}{2A}- \frac{1}{4A}}
	e^{-D(\mu_n u)^{\frac{1}{A}}}
	\Bigg). 
\end{align}
 To obtain a lower bound for $F(x)$ similar to 
\eqref{final bdd F(k^2)}, we require a suitable choice of the set $\mathcal{M}$ to make the error term smaller.  In \cite{lamzouri2026}, for the divisor problem, this  improvement  is obtained by choosing \(\mathcal{M}\) to consist of the integers corresponding to the \(|\mathcal{M}|\) largest values of the sequence \(d(n)n^{-3/4}\). 
\normalcolor
While the above approach can be applied to specific examples, the
corresponding estimates depend on the arithmetic structure of the
particular coefficient sequence. Such estimates require information
beyond the hypotheses of Theorem~\ref{thm:main theorem}, whose
formulation provides a unified treatment of the omega bounds
established in this paper.

\normalcolor

\section*{Acknowledgment}
The second author was  partially supported by the Vikram Sarabhai Research Fellowship provided by IIT Gandhinagar, under project number MIS/END-VSRF/MA/AV/2024-25/084. The third author is   grateful to the Max Planck Institute for Mathematics  in Bonn for its hospitality and financial support.

\section*{Conflict of interest statement and Data availability}
The authors confirm that there is no conflict of interest. 
	Data sharing is not applicable to this article as no datasets were generated or analysed during the
	current study.
\allowdisplaybreaks

\bibliographystyle{abbrv}
\bibliography{rsv_omega}

\end{document}